\documentclass[axioms,article,submit,moreauthors,pdftex]{Definitions/mdpi} 

\firstpage{1} 
\pubvolume{0}
\issuenum{0}
\articlenumber{0}
\pubyear{2021}
\copyrightyear{2021}
\datereceived{--} 
\dateaccepted{--} 
\datepublished{--} 
\hreflink{https://doi.org/} 

\usepackage{overpic}
\usepackage{subfig}
\usepackage{empheq}

\newcommand{\rev}[1]{\textcolor{black}{#1}}

\newcommand{\R}{\mathbb R}
\newcommand{\subfig}[1]{(\textbf{#1})}

\newcommand{\dx}{\Delta x}
\newcommand{\dt}{\Delta t}
\newcommand{\leg}{\textsc{l}} 
\newcommand{\pes}{\textsc{h}} 
\newcommand{\rhoL}{\rho_\textsc{l}}
\newcommand{\rhoLstar}{\rho^*_\textsc{l}}
\newcommand{\VLstar}{V^*_\textsc{l}}
\newcommand{\rhoPstar}{\rho^*_\textsc{h}}
\newcommand{\rhoP}{\rho_\textsc{h}}
\newcommand{\SL}{S_\textsc{l}}
\newcommand{\SP}{S_\textsc{h}}
\newcommand{\RL}{R_\textsc{l}}
\newcommand{\RP}{R_\textsc{h}}
\newcommand{\FL}{\mathcal{F}_\textsc{l}}
\newcommand{\FP}{\mathcal{F}_\textsc{h}}

\newcommand{\lL}{\ell_\textsc{l}}
\newcommand{\lP}{\ell_\textsc{h}}
\newcommand{\mumax}{\rhoP^{\text{max}}}
\newcommand{\rhomax}{\rhoL^{\text{max}}}
\newcommand{\Domain}{\mathcal{D}}

\Title{
\rev{Macroscopic} and multi-scale models for multi-class vehicular dynamics with uneven space occupancy: a case study
}

\TitleCitation{
Macroscopic and multi-scale models for multi-class vehicular dynamics with uneven space occupancy: a case study
}

\Author{
	Maya Briani $^{1}$\orcidA{}, 
	Emiliano Cristiani $^{1}$\orcidB{},
	Paolo Ranut $^{2}$
}

\AuthorNames{Maya Briani, Emiliano Cristiani, Paolo Ranut}

\AuthorCitation{Briani, M.; Cristiani, E.; Ranut, P.}

\address{
	\hspace{-0.8mm}$^{1}$ \quad Istituto per le Applicazioni del Calcolo, Consiglio Nazionale delle Ricerche, Rome, Italy\\
	$^{2}$ \quad Autovie Venete S.p.A. 
}

\corres{Correspondence: m.briani@iac.cnr.it; e.cristiani@iac.cnr.it}

\abstract{In this paper we propose two models describing the dynamics of heavy and light vehicles on a road network, taking into account the interactions between the two classes.
The models are tailored for two-lane highways where heavy vehicles cannot overtake. This means that heavy vehicles cannot  saturate the whole road space, while light vehicles can. In these conditions the creeping phenomenon can appear, i.e.\ one class of vehicles can proceed even if the other class has reached the maximal density.
The first model we propose couples two first-order macroscopic LWR models, while the second model couples a second-order microscopic Follow-the-Leader model with a first-order macroscopic LWR model.
Numerical results show that both models are able to catch some second-order (inertial) phenomena like stop \& go waves.  
Models are calibrated by means of real data measured by fixed sensors placed along the A4  Italian highway Trieste-Venice and its branches, provided by Autovie Venete S.p.A.}

\keyword{LWR model; Follow-the-Leader model; phase transition; creeping; seepage; fundamental diagram; lane discipline; networks} 

\MSC{35L65, 35F25, 90B20, 76T99}

\begin{document}
	
	\section{Introduction}\label{sec:intro}
	
	In this paper we deal with macroscopic and multi-scale modeling of traffic flow on a road network, focusing on multi-class dynamics which couple light and heavy vehicles (in the following, cars and trucks). The proposed models are characterized by the fact that cars and trucks interact with each other and that trucks are confined in a part of the road space (slow lane) and cannot overtake. As a consequence, when trucks saturate the space and form a queue, cars can still move, although at reduced speed.
	
	\paragraph{State of the art.} The literature about traffic flow is very large and many different aspects of traffic dynamics were described through mathematical models. 
	Let us start from classic approaches:
	in a single-lane \textit{microscopic} (agent-based) framework with $N$ vehicles and no overtaking, each vehicle $k\in\{1,\ldots,N\}$ is singularly identified by its position $X_k(t)$ and its velocity $V_k(t)$. 
	By assumption, the $(k+1)$-th vehicle is always in front of the $k$-th one.
	Also, each vehicle is assumed to adjust its acceleration based on the difference in positions and velocities between the vehicle itself and the vehicle in front of it. 
	This approach leads to the following system of ordinary differential equations:
	\begin{equation}\label{FtL}
	\left\{
	\begin{array}{ll}
	\dot X_k=V_k \\ [2mm]
	\dot V_k=A(X_k,X_{k+1},V_k,V_{k+1})
	\end{array}, \quad k=1,\ldots,N-1
	\right.
	\end{equation}
	where $A$ is a given acceleration function. The first vehicle in the row ($k=N$), called \emph{leader}, has an independent dynamics. 
	Since the whole dynamics is determined by the leader's one in a domino effect, these kind of models are known as Follow-the-Leader.
	
	Adopting instead a \textit{macroscopic} (fluid-dynamics) point of view, we describe the mass of vehicles by means of its density $\rho(x,t)$ only. The celebrated LWR model \cite{lighthill1955PRSLA, richards1956OR} is based on the observation that the density $\rho$ evolves in time ruled by the following conservation law 
	\begin{equation}\label{LWR}
	\partial_t\rho + \partial_x f(\rho)=0,\qquad x\in\R, \quad t>0
	\end{equation}
	where the function $f(\rho)$, called \emph{fundamental diagram}, is given and represents the flux of vehicles as a function of the density itself. 
	The velocity of the vehicles can be recovered from $\rho$ thanks to the relation 
	\begin{equation}\label{f=rhov}
	v(\rho)=\frac{f(\rho)}{\rho}, \qquad (\rho\neq 0).
	\end{equation}
	 
	It is important to note that the microscopic model \eqref{FtL} is second-order, i.e.\ acceleration based, while the macroscopic model \eqref{LWR} is first-order, i.e.\ velocity based. 
	The difference is important because velocity based models, allowing nonphysical instantaneous accelerations, are not able to catch effects caused by inertia, like stop \& go waves.
	Due to this difference, it is plain that the model \eqref{LWR} is \emph{not} the many-particle limit of the model \eqref{FtL}. 
	We refer the interested reader to the paper \cite{cristiani2019DCDS-B} for a review of various many-particle limits (i.e.\ micro-to-macro correspondences) and the existing multi-scale models. 
	
	\medskip
	
	A first generalization of the models \eqref{FtL} and \eqref{LWR} is that of \emph{road networks}. While managing junctions in microscopic models is relatively easy, doing the same in a macroscopic setting is more challenging. The reason is that, in general, the conservation of the mass alone is not sufficient to characterize a unique solution at junctions. 
	We refer the reader to the book by Garavello and Piccoli \cite{piccolibook} for more details about the ill-posedness of the problem at junctions. Multiple workarounds for such ill-posedness have been suggested in the literature: 
	(i) maximization of the flux across junctions and introduction of priorities among the incoming roads \cite{coclite2005SIMA, piccolibook, holden1995SIMA}; 
	(ii) introduction of a \emph{buffer} to model the junctions by means of additional ordinary differential equation coupled with \eqref{LWR} \cite{bressan2015NHM, garavello2012DCDS-A, herty2009NHM}; 
	(iii) reformulation of the problem on all possible paths on the network rather than on roads and junctions. The last approach has both a global formulation \cite{bretti2014DCDS-S, briani2014NHM, hilliges1995TRB} and a more manageable local formulation, described in \cite{briani2018CMS}, which is the one we will adopt in this paper.
	All these approaches allow to determine a unique solution for the traffic evolution on the network, but the solutions might be different.
	
	A second generalization of our interest is that of \emph{multi-class} dynamics. ``Multi-class'' is a very generic term used in literature to refer to the case in which the road is populated by different groups of vehicles/drivers, and tracking each group separately is desired. 
	Again, doing this in a microscopic framework is easy since it is sufficient to label each vehicle on the basis of the class it belongs to. 
	In the macroscopic setting, instead, we need to introduce as many density functions as there are classes, and then establish the interactions between classes. This leads to a system of conservation laws of the form
	\begin{equation}\label{LWRmc}
	\partial_t \rho_c+ \partial_x f_c(\rho_1,\ldots,\rho_C)=0,\qquad c=1,\ldots,C, \quad x\in\R, \quad t>0
	\end{equation}
	where $C$ is the number of classes, $c\in\{1,\ldots,C\}$, $\rho_c$ is the density of class $c$, and $f_c$ is the flux of class $c$ which depends of all densities (usually the dependence is on the sum of all densities $\sum_c\rho_c$).
	Multi-class models are used to describe very different situations, like the co-presence of vehicles with
	\begin{itemize}
		\item different driving modes (e.g.\ autonomous vs.\ classic);
		\item different origins and destinations;
		\item different length (i.e.\ space occupied);
		\item different velocity/flux function;
		\item reserved roads or reserved entry/exit lanes.
	\end{itemize}
	A complete review of multi-class models is out of the scope of this paper. We refer to the papers \cite{fan2015SIAP, kessels2016TRR, qian2017TRB} and to the recent books \cite{ferrarabook, kesselsbook} for an overview of the most used multi-class models.
	
	Before introducing our contributions, let us introduce the \emph{creeping or seepage effect} \cite{fan2015SIAP, agarwal2016TDE} which will be useful to describe the features of the proposed models. This term denotes the situations where the road space is shared by small and large vehicles, and small vehicles are able to move (at reduced velocity) even if large vehicles have reached the maximal density.
	This is in contrast with classical models like the one proposed by Benzoni-Gavage and Colombo \cite{benzoni2003EJAM}, in which the saturation of a class of vehicles stops immediately all the other classes. 
	It is useful to note that the creeping phenomenon is typically considered in a context of \emph{disordered traffic}, i.e.\ a traffic with no lane discipline: smaller vehicles (e.g., two wheels) slip into the empty spaces left by large vehicles, similarly to a motion through porous media.
	This is not the case considered here, since we assume a strict lane discipline.
	
	Finally let us recall some important contributions about the fundamental diagram and its properties. It is well known that a single function $f=f(\rho)$ is not able alone to describe real data correctly. Indeed, by \eqref{f=rhov} we deduce that for any given density value $\rho$ only one velocity $v(\rho)$ is possible. 
	This is not what happens in reality, where a scattered fundamental diagram is observed instead, due to the fact that different drivers respond in a different way to the same traffic conditions. 
	Many papers investigated this phenomenon from different points of view, trying to explain its features, including instabilities; see, e.g., \cite{balzotti2021NHM, fan2013TRR, fan2014NHM, klar2004SIAP, herty2012KRM, ni2018AMM, paipuri2020TRB, puppo2016CMS, tosin2017MMS, wang2013JAT, CGARZ}.
	
\begin{figure}[h!]
		\widefigure
		\includegraphics[width=12cm]{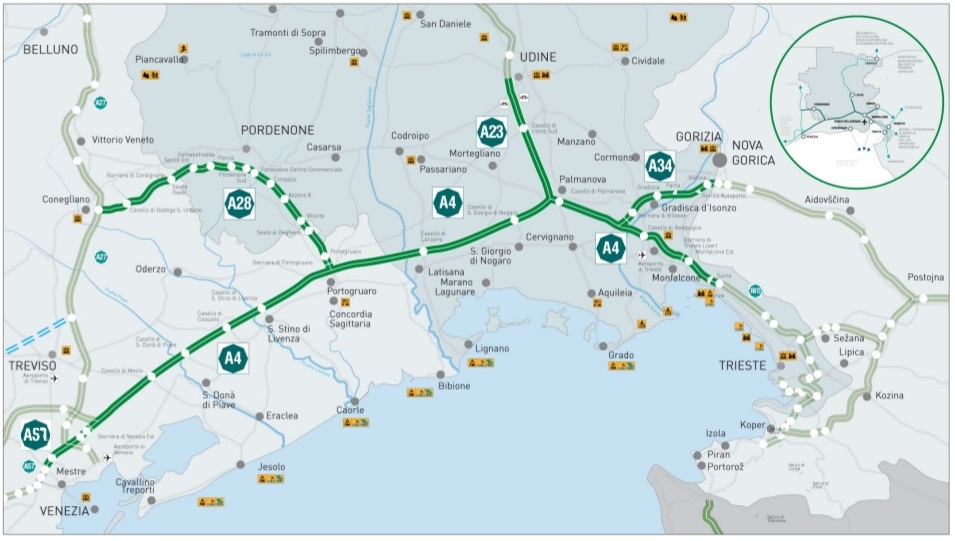}
		\caption{The Italian motorway A4 Trieste-Venice and its branches to/from Udine, Pordenone and Gorizia, managed by Autovie Venete S.p.A.}
		\label{fig:reteAV}	
	\end{figure}

	
	\paragraph{Case study.}
	In this paper we consider the Italian motorway A4 Trieste-Venice and its branches to/from Udine, Pordenone and Gorizia, managed by Autovie Venete S.p.A., see Fig.\ \ref{fig:reteAV}.
	\rev{At the time of the present study (2019), the motorway had two lanes per direction, except for the leftmost segment near Venice (Venice -- San Don\`a). To avoid heterogeneous conditions, we have dropped the three-lane segment of the road to focus exclusively on the parts with \emph{two lanes per direction}.} 
	In those segments, cars can use both lanes at any time, while trucks can use only the slow lane and cannot overtake. Due to the large flow of heavy vehicles, it happens some times that a queue of trucks is formed. 
	In this case, cars move into the fast lane and keep going, although at moderate speed. 
	When traffic conditions are sustained, the two classes of vehicles interact with each other: on the one hand, trucks act as moving bottlenecks for cars (cf.\ \cite{bretti2018MiE}), which are forced to slow down due to the restricted space; 
	on the other hand, trucks must slow down when cars find it convenient to occupy part of the slow lane.

	\paragraph{Our contribution.}
	In this paper we propose two models for describing multi-class traffic flow on networks in which vehicles belonging to different classes share the road space only partially. More precisely, light vehicles can occupy the whole road, while heavy vehicles only a part of it.
	To align with the case study, we will assume that the road has two lanes in total and trucks can occupy only the slow one, without overtaking.  
	\begin{enumerate}
		\item The first model is purely macroscopic. Both cars and truck are described by two coupled first-order LWR-based models. Fundamental diagrams are shaped in order to allow cars to move even in presence of fully congested trucks.
		Considering that the fundamental diagram of each class is influenced by the presence of the other class, in case of unstable (rapidly varying) traffic conditions of one class we observe a scattered behavior in the fundamental diagram of the other class. 
		Numerical results will show that this feature allows the model to catch, at least in part, some second-order (inertial) phenomena in traffic behavior, like stop \& go waves.
		\item The second model is multi-scale. Cars are described by a first-order LWR-based model, while trucks are described by a second-order microscopic Follow-the-Leader model. 
		For trucks, we consider the microscopic model used in \cite{cristiani2019DCDS-B}, inspired, in turn, by a model originally proposed in \cite{zhao2017TRB} and specifically designed to reproduce  stop \& go waves.
		The choice of second-order model for trucks is crucial, since inertia effects are not at all negligible for those vehicles, while they are less important in car dynamics. 
		Finally note that, since trucks are confined in only one lane and cannot overtake, their dynamics perfectly matches the constituting assumptions of the Follow-the-Leader model.
	\end{enumerate} 
		Let us finally mention that the idea of coupling first- and second-order models was already exploited in \cite{cristiani2019DCDS-B} in a single-class scenario.
	
	\begin{Remark} 
	Both models distinguish classes, but not lanes. The fact that trucks cannot use the fast lane while cars can occupy both slow and fast lanes is encapsulated in the choice of the fundamental diagrams.
	\end{Remark}

	\section{Dataset}\label{sec:data}
\begin{figure}[h!]
	 	\subfig{a}\includegraphics[width=12cm]{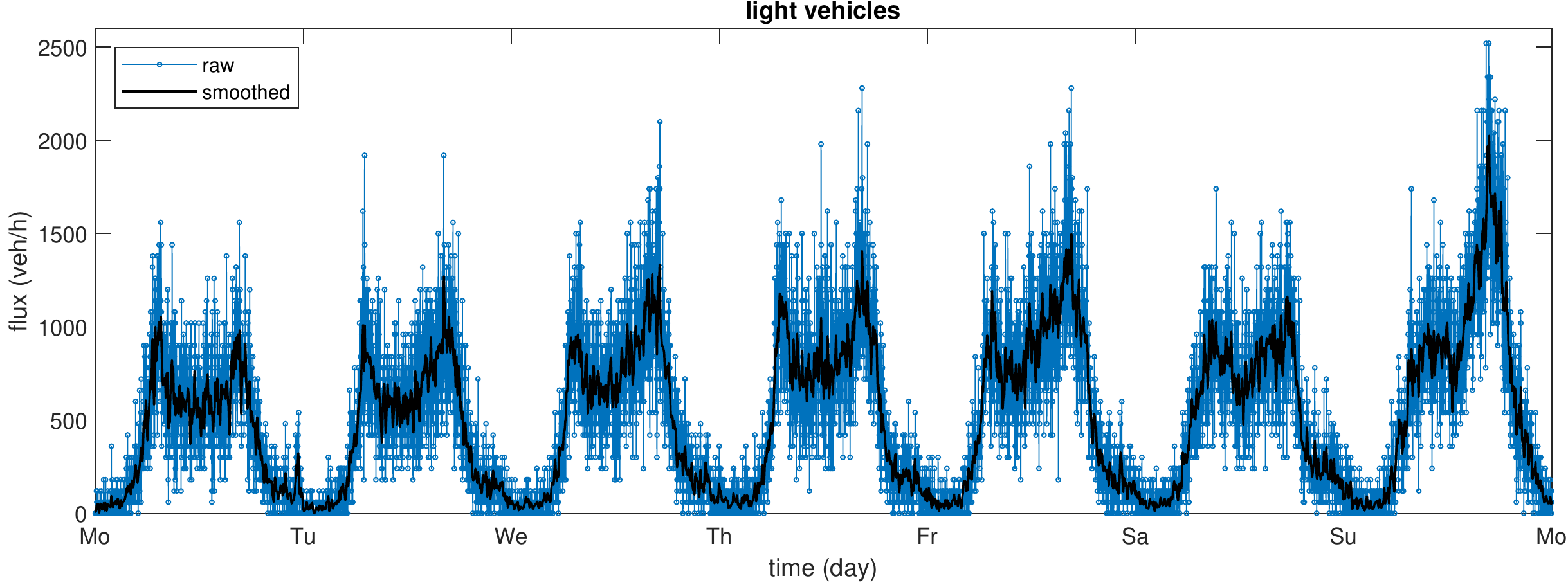}\\  [2mm]
	 	\subfig{b}\includegraphics[width=12cm]{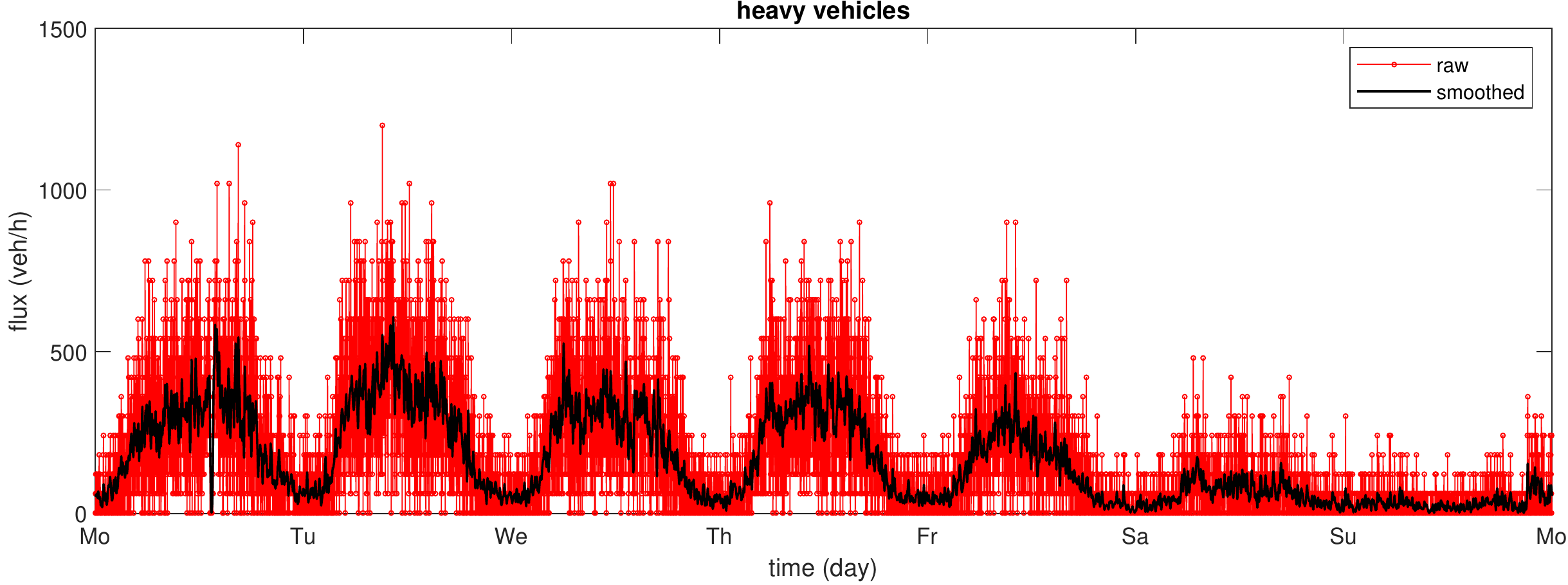}
	 	\caption{Typical weekly (from Monday to Sunday) flux data on the A4 motorway of \subfig{a} light and \subfig{b} heavy vehicles \rev{collected on March 2019 near Redipuglia}. Smoothed data are plotted in black. Note the flux drop of cars in the middle of the day and of trucks in the weekend.}
	 	\label{fig:weeklyflux}	
	 \end{figure}

\begin{figure}[h!]
 		\subfig{a}\includegraphics[width=6cm]{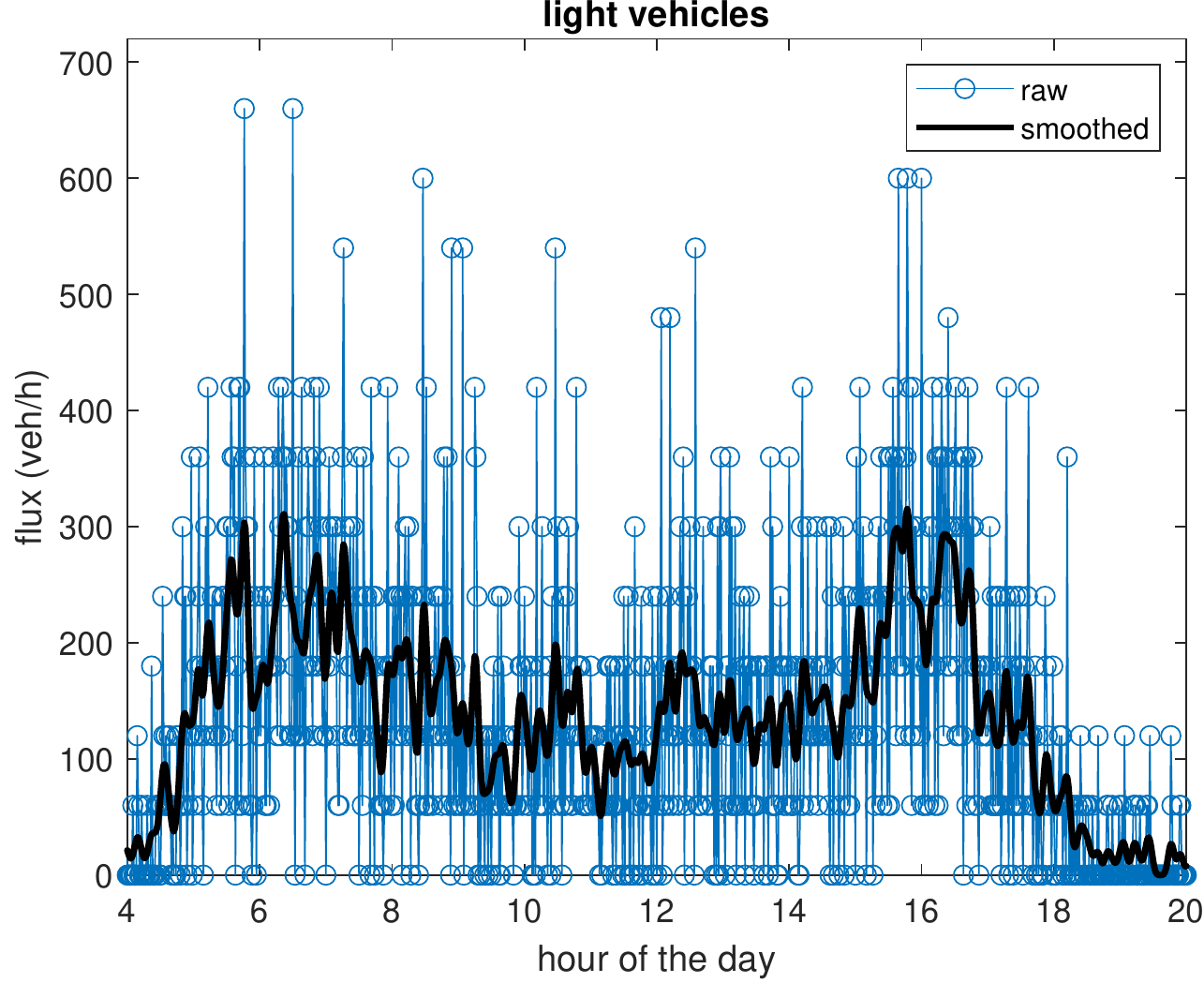}
 		\subfig{b}\includegraphics[width=6cm]{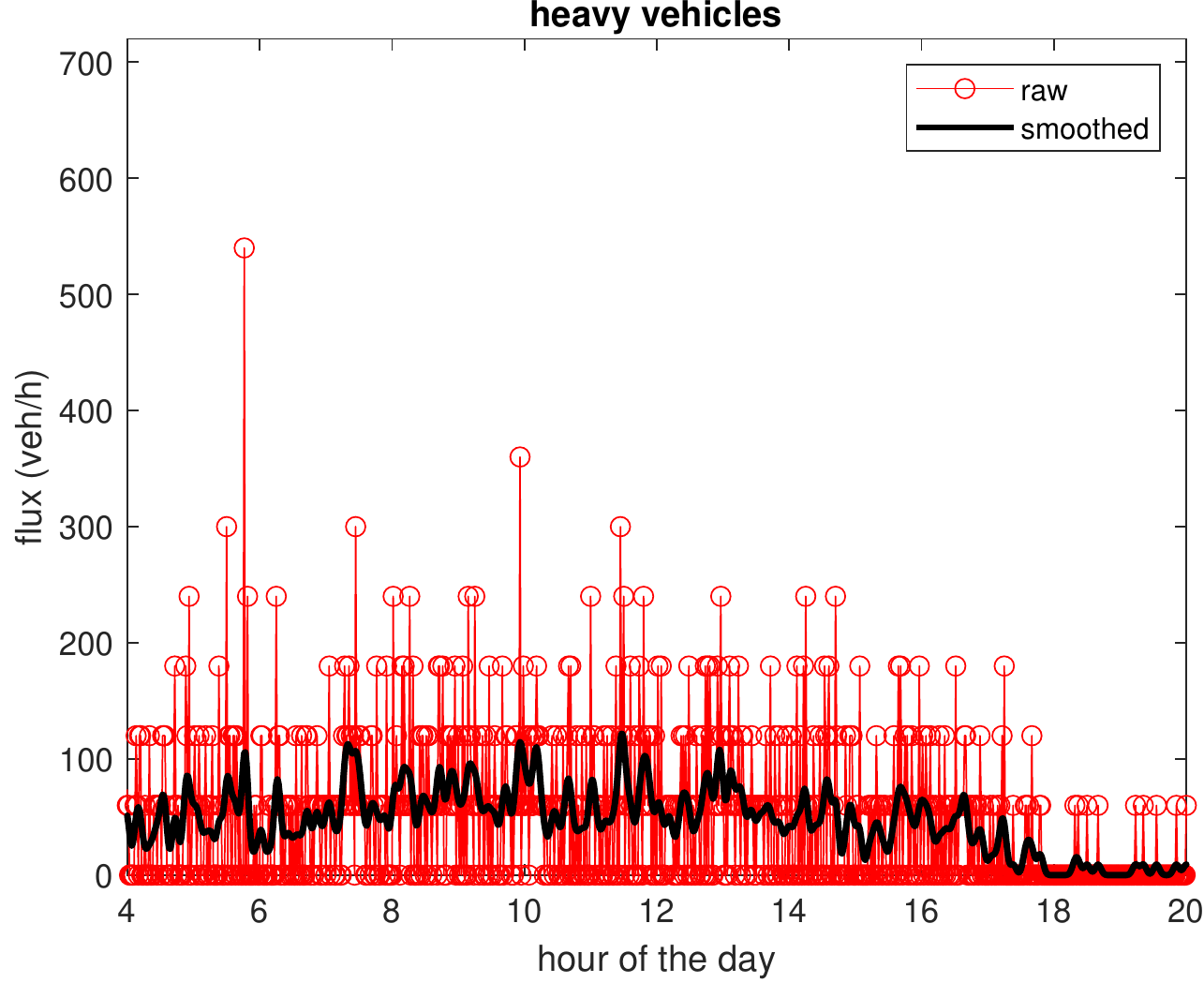}\\  [2mm]
 		\subfig{c}\includegraphics[width=6cm]{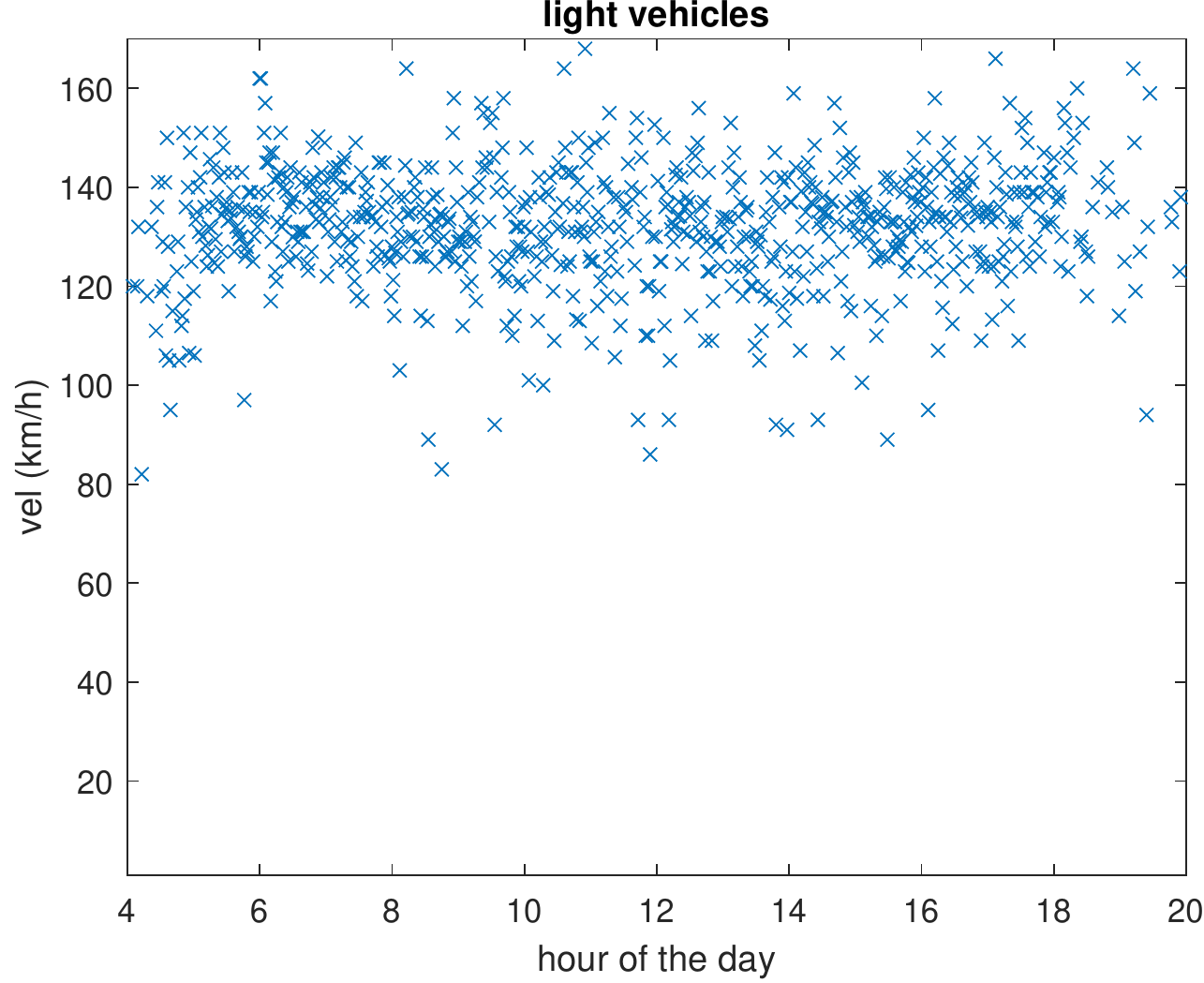}
 		\subfig{d}\includegraphics[width=6cm]{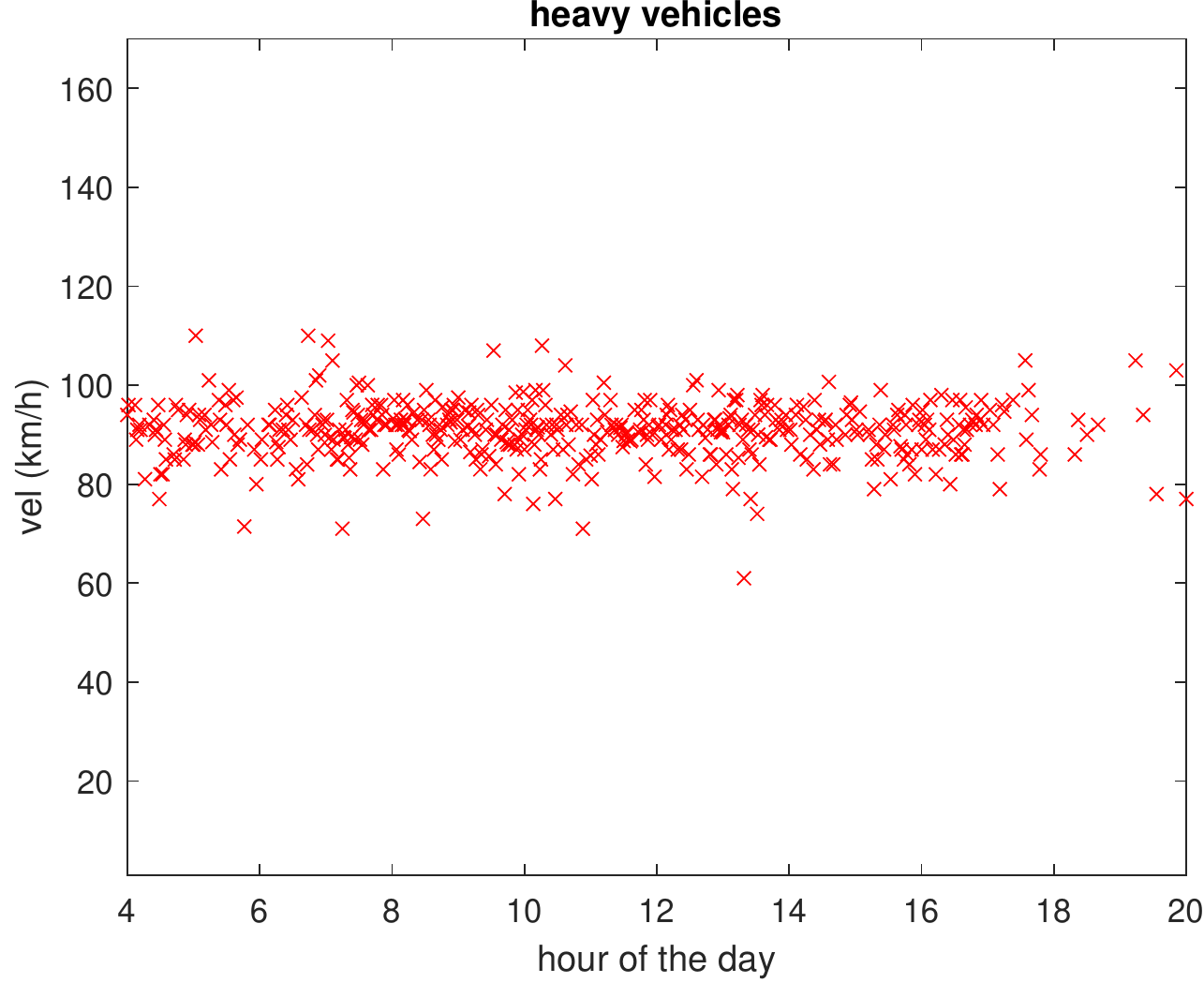}
 		\caption{Typical daily (Thursday) flux and velocity data on the A28 motorway of \subfig{a}-\subfig{c} light and \subfig{b}-\subfig{d} heavy vehicles \rev{collected on May 2019 near Sesto al Reghena}.}
 		\label{fig:normalfluxvel}	
 	\end{figure}

\begin{figure}[h!]
		\subfig{a}\includegraphics[width=6cm]{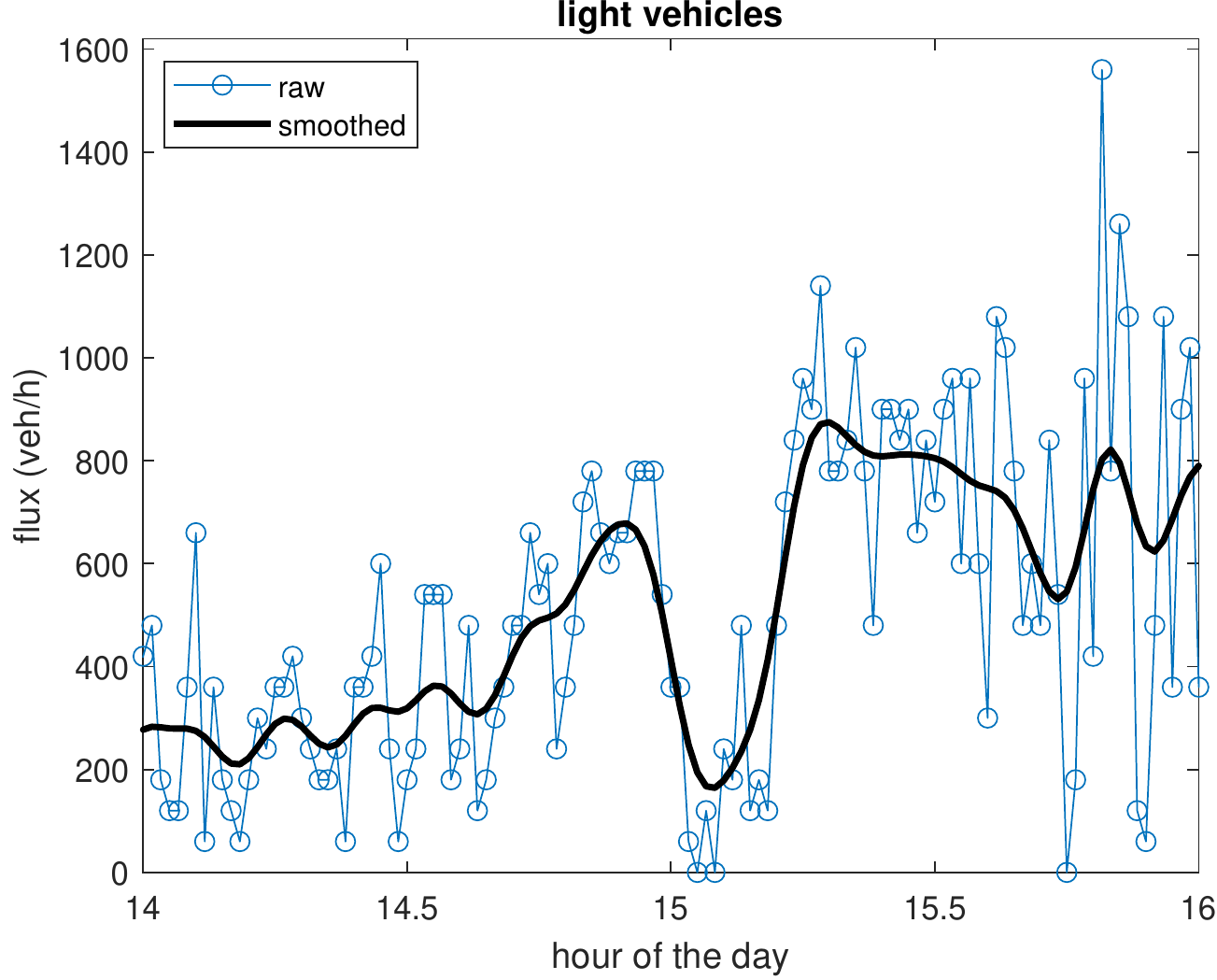}
		\subfig{b}\includegraphics[width=6cm]{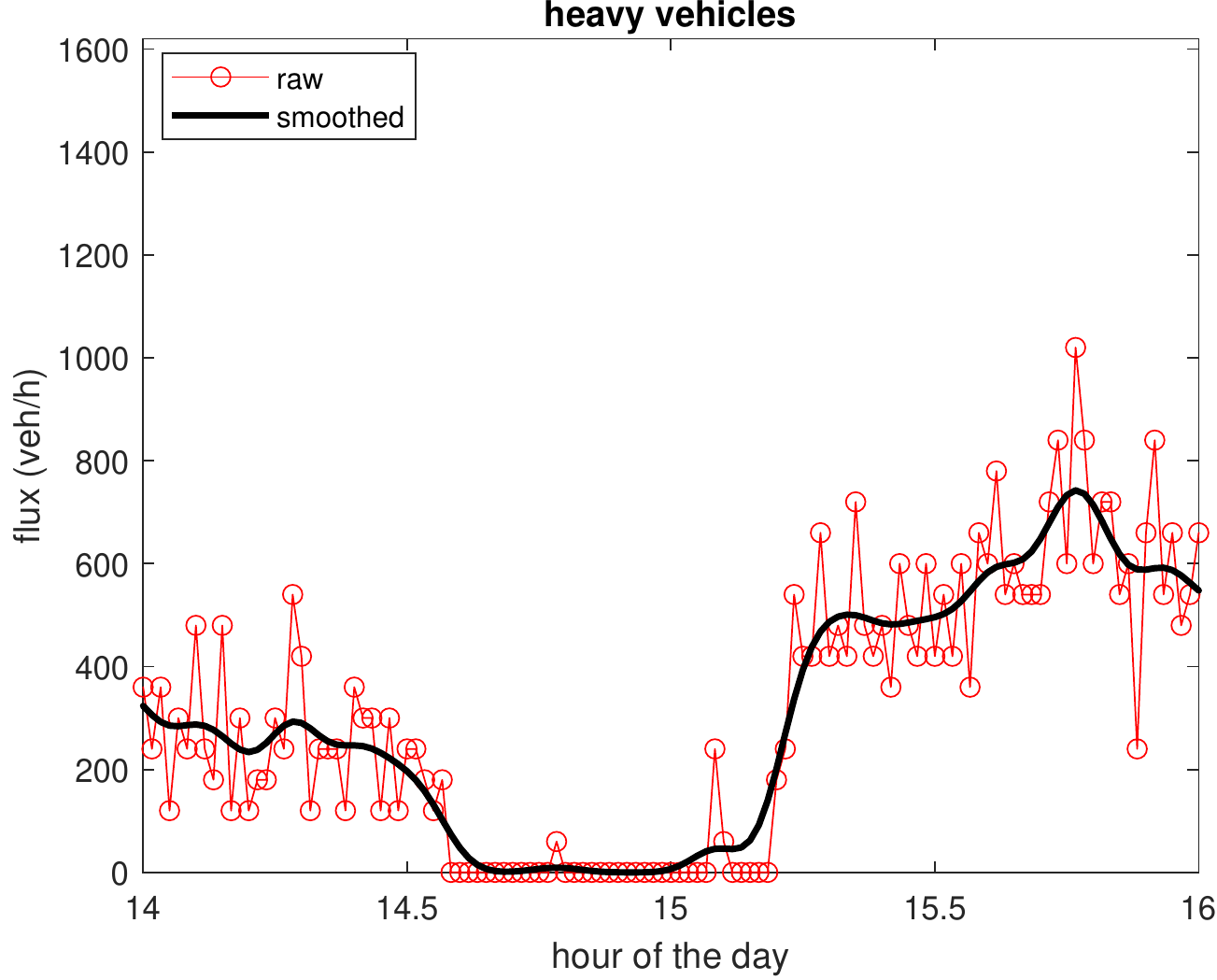}\\ [2mm]
		\subfig{c}\includegraphics[width=6cm]{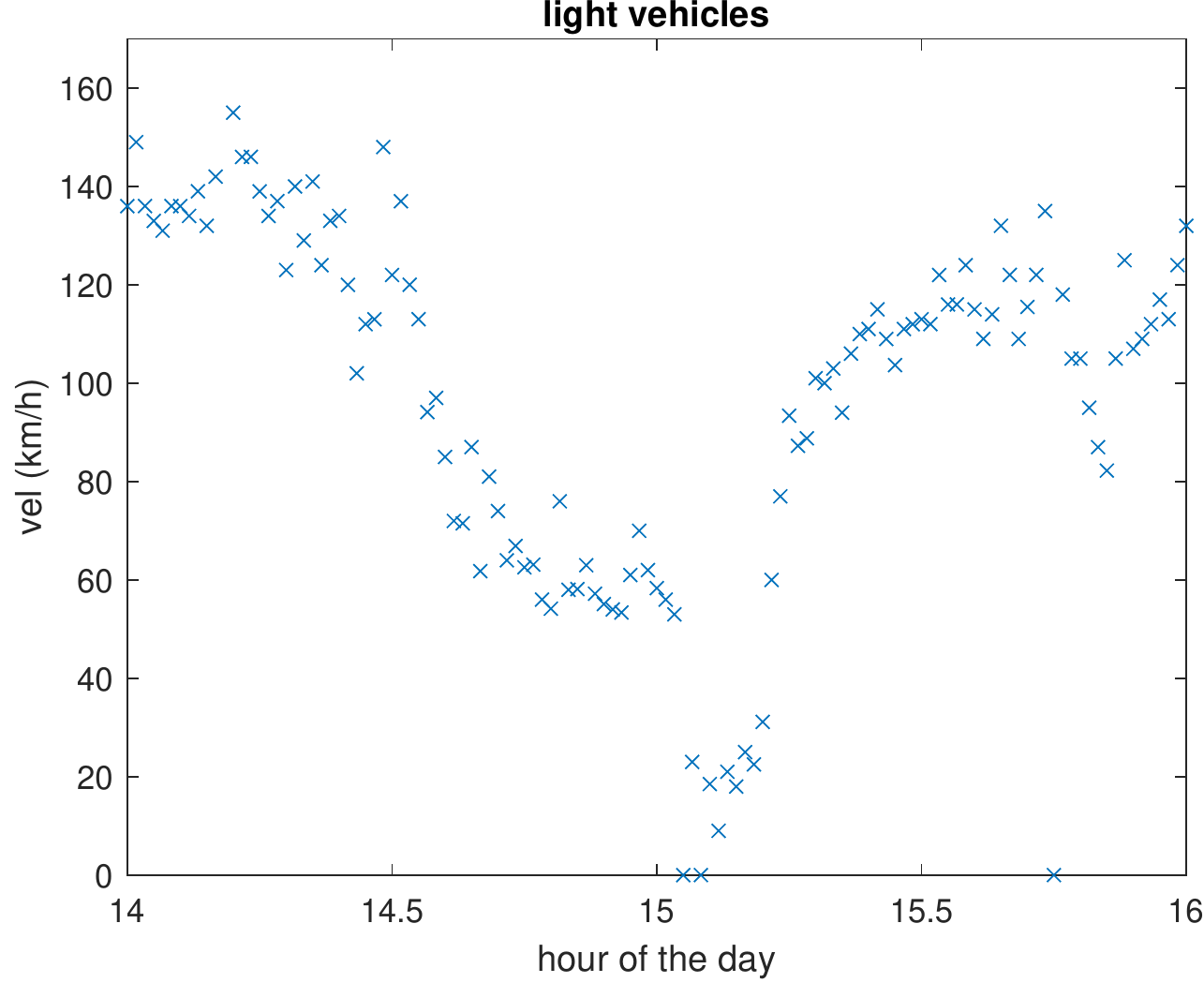}
		\subfig{d}\includegraphics[width=6cm]{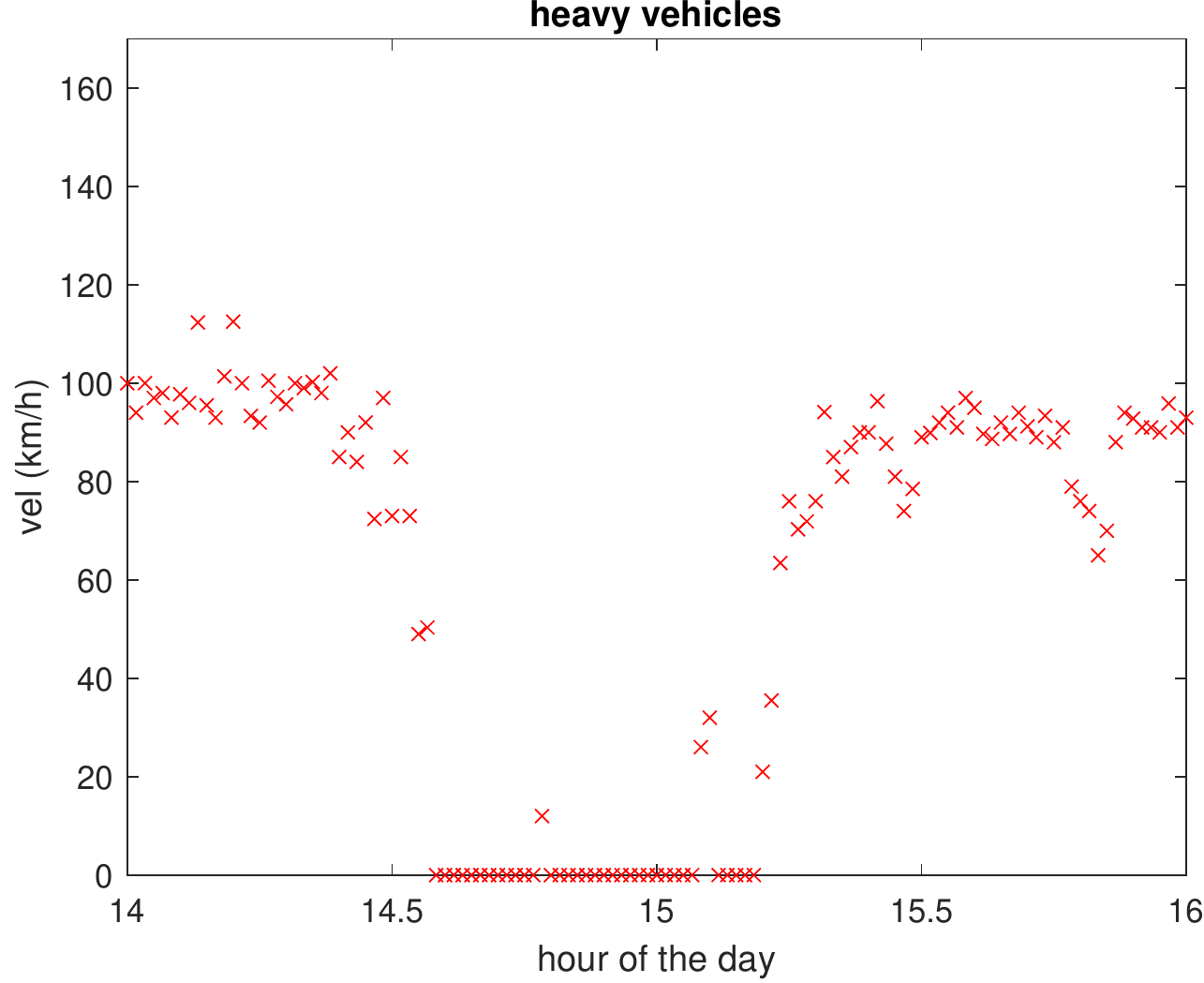}
		\caption{Creeping phenomenon \rev{registered on May 2019 near Portogruaro}: \subfig{a} Light vehicles move in the fast lane even if \subfig{b} heavy vehicles queue in the slow lane. \subfig{c} Light vehicles' velocity drops from $\sim$140 km/h to $\sim$60 km/h and then to $\sim$20 km/h while \subfig{d} heavy vehicles are completely stopped.}
		\label{fig:creepingfluxvel}
	\end{figure}

	 Autovie Venete constantly monitors traffic conditions by means of video cameras, mobile sensors, and fixed sensors. \rev{In this paper we focus on the latest kind of data}. 
	 Fixed sensors are located along the motorway, on each lane, and measure flux and velocity of all vehicles passing in front of them, distinguishing also the class of vehicles. 
	 Data are aggregated per minute and are stored in a data base for later analysis.
	 \rev{For light vehicles, we have further aggregated data coming from slow and fast lanes. 
	 For heavy vehicles, instead, we have considered the slow lane only.}
	 In Figs.\ \ref{fig:weeklyflux}--\ref{fig:creepingfluxvel} we show some flux and velocity data coming from some fixed sensors, used to conceive and calibrate the models  presented in this paper. For better readability, flux data are plotted both raw (as is) and smoothed by a Gaussian filter. \rev{Note that the flux data is always a multiple of 60 since it is evaluated every minute but it is expressed in terms of vehicles per hour.}


	\section{Models}\label{sec:models}
	
	In this section we present the two models. As already stated in the Introduction, the models are not meant to give the same results or to be one the many-particle limit of the other. 
	Nevertheless, they share the most important constitutive assumptions and for this reason they are expected to give the same qualitative results. 
	The most important common modeling assumption is that the car dynamics is influenced at any time by the presence of trucks, while the truck dynamics is affected by cars only if the density of cars exceed a certain threshold, which corresponds to the fact that cars cannot be confined in the fast lane any longer and must invade the slow lane where trucks live. 
	This assumption comes from an important evidence: cars tend to avoid to be trapped between two trucks in the slow lane, and prefer moving to the fast lane. Doing this, cars move to the side of trucks (overtaking them if possible) and do not affect their dynamics, unless the density of cars is so high that they must necessarily occupy the slow lane too.

	\subsection{Macroscopic model}\label{sec:mama_model}
	We denote by $\lL$ and by $\lP$ the average length of  light vehicles (cars) and heavy vehicles (trucks), respectively, and we define  
	\begin{equation}\label{eq:beta}
	\beta := \frac{\lL}{\lP}<1.
	\end{equation}
	We also denote by $\rhoL$ the density of cars and by $\rhoP$ the density of trucks. 
	Similarly, we denote by $\rhomax$ and $\mumax$ the maximal densities for cars and trucks, respectively. 
	\rev{They are defined as
		\begin{equation}\label{def_rhomax}
			\rhomax = \frac{2}{\lL}
			\quad\mbox{ and }\quad 
			\mumax = \frac{1}{\lP},
		\end{equation}
		having assumed that there are two available lanes for cars and only one for trucks.
	Note that density values are expressed in terms of number of vehicles per unit of space}. Considering that trucks occupy more space than cars, a direct comparison of the two densities is not meaningful. For this reason the two classes are typically compared in terms of occupied space.
	
	The two-class dynamics is physically admissible if the two densities fall in the set
	\begin{equation}\label{eq:dominio}
	\Domain:=\left\{(\rhoL,\rhoP) :  0\leq \rhoL\leq \rhomax,\ 0\leq \rhoP\leq \mumax,\  0\leq \rhoL+\frac{\rhoP}{\beta}\leq \rhomax\right\},
	\end{equation}
	which is well defined if $\rhomax-\frac{\mumax}{\beta}\geq 0$.
	In the following, in order to cope with the uneven space occupancy, we assume that the last condition is verified with the strict inequality,
	\begin{equation}\label{eq:DensCondstretta}
	\rhomax-\frac{\mumax}{\beta} > 0. 
	\end{equation}
	
	We consider the following two-class model for
	$(\rhoL,\rhoP)\in\Domain$, 
	\begin{equation}\label{eq:model}
	\begin{cases}
	\partial_t\rhoL + \partial_x f_{\leg}(\rhoL,\rhoP) =\, 0
	\smallskip\\
	\partial_t \rhoP + \partial_x f_{\pes}(\rhoL,\rhoP) =\, 0
	\end{cases}
	\qquad x\in\R,\quad t>0,
	\end{equation}
	where
	$$
	f_\leg(\rhoL,\rhoP):=\rhoL v_	\leg(\rhoL,\rhoP),
	\qquad
	f_\pes(\rhoL,\rhoP):=\rhoP v_\pes(\rhoL,\rhoP)
	$$
	define the two fundamental diagrams and $v_\leg$, $v_\pes$ are the speed functions for light and heavy vehicles, respectively.
	We then have a family of flow-density curves $\rhoL\mapsto f_\leg(\rhoL,\rhoP)$ for cars, parameterized by the trucks density $\rhoP$, and analogously a family of flow-density curves $\rhoP\mapsto f_\pes(\rhoL,\rhoP)$ for trucks, parameterized by $\rhoL$.
	
	We assume that the flux and speed functions satisfy the following properties:
	\begin{enumerate}[label=(L\arabic*)]
		\item $v_{\leg}(\rhoL,\rhoP) \geq 0$ for all $(\rhoL,\rhoP)\in\mathcal{D}$ and  $v_\leg(\rhoL,\rhoP) = 0$ iff $\rhoL = \rhoLstar(\rhoP)$, where 
		\begin{equation}\label{eq:RhoMax}
		\rhoLstar(\rhoP) := \rhomax - \rhoP/\beta
		\end{equation}
		is the maximum admissible cars density given the trucks density $\rhoP$;		
		\item $v_\leg(\rhoL,\rhoP)$ is a decreasing function with respect to $\rhoL$ and $\rhoP$;
		%
		\item $f_\leg(0,\rhoP) = 0$ and $f_\leg(\rhoLstar(\rhoP),\rhoP) = 0$ for all $\rhoP\in[0,\mumax]$;
		\item $f_\leg(\rhoL,\rhoP)$ is concave with respect to $\rhoL$ for any $\rhoP$. 
		We define 
		\begin{equation}
		\sigma_\leg(\rhoP):=\arg\max_{\rhoL} f_\leg(\rhoL,\rhoP) 
		\end{equation} 
		which represents as usual the interface between freeflow and congested regimes;
		\item $f_\leg(\rhoL,\rhoP)$ is a decreasing function with respect to $\rhoP$ for any $\rhoL$.
	\end{enumerate}

	Similarly,
	\begin{enumerate}[label=(H\arabic*)]
		\item $v_{\pes}(\rhoL,\rhoP) \geq 0$ for all $(\rhoL,\rhoP)\in\mathcal{D}$ and $v_\pes(\rhoL,\rhoP) = 0$ iff $\rhoP = \rhoPstar(\rhoL)$, where 
		\begin{equation}\label{eq:MuMax}
		\rhoPstar(\rhoL) := \min\left\{\mumax, \beta(\rhomax - \rhoL) \right\}
		\end{equation}
		is the maximum admissible trucks density given the cars density $\rhoL$;
		\item $v_\pes(\rhoL,\rhoP)$ is a decreasing function with respect to $\rhoL$ and $\rhoP$;
		\item $f_\pes(\rhoL,0) = 0$ and $f_\pes(\rhoL,\rhoPstar(\rhoL)) = 0$ for all $\rhoL\in[0,\rhomax]$;
		\item $f_\pes(\rhoL,\rhoP)$ is concave with respect to $\rhoP$ for any $\rhoL$. We define
		\begin{equation}
		\sigma_\pes(\rhoL):=\arg\max_{\rhoP} f_\pes(\rhoL,\rhoP) 
		\end{equation} 
		which represents as usual the interface between freeflow and congested regimes;
		\item $f_\pes(\rhoL,\rhoP)$ is a decreasing function with respect to $\rhoL$ for any $\rhoP$.
	\end{enumerate}
	
	\medskip
	
	To cope with the peculiarities of the dynamics, we consider a \emph{phase transition} (\rev{cf.\ }\cite{colombo2002SIAP, colombo2010JHDE, dellemonache2021AX}) caused by the presence of two states of the system:
	\begin{itemize}
		\item The \emph{partial-coupling phase} is in place when 
		\begin{equation}
		(\rhoL,\rhoP)\in\Domain_1:=\Big\{0\leq\rhoP\leq\mumax, \ 0\leq\rhoL\leq \rhomax-\mumax/\beta\Big\},
		\end{equation}
		see Fig.\ \ref{fig:dominio_MaMa}.
		\begin{figure}[h!]
			\begin{overpic}
				[width=0.7\columnwidth]{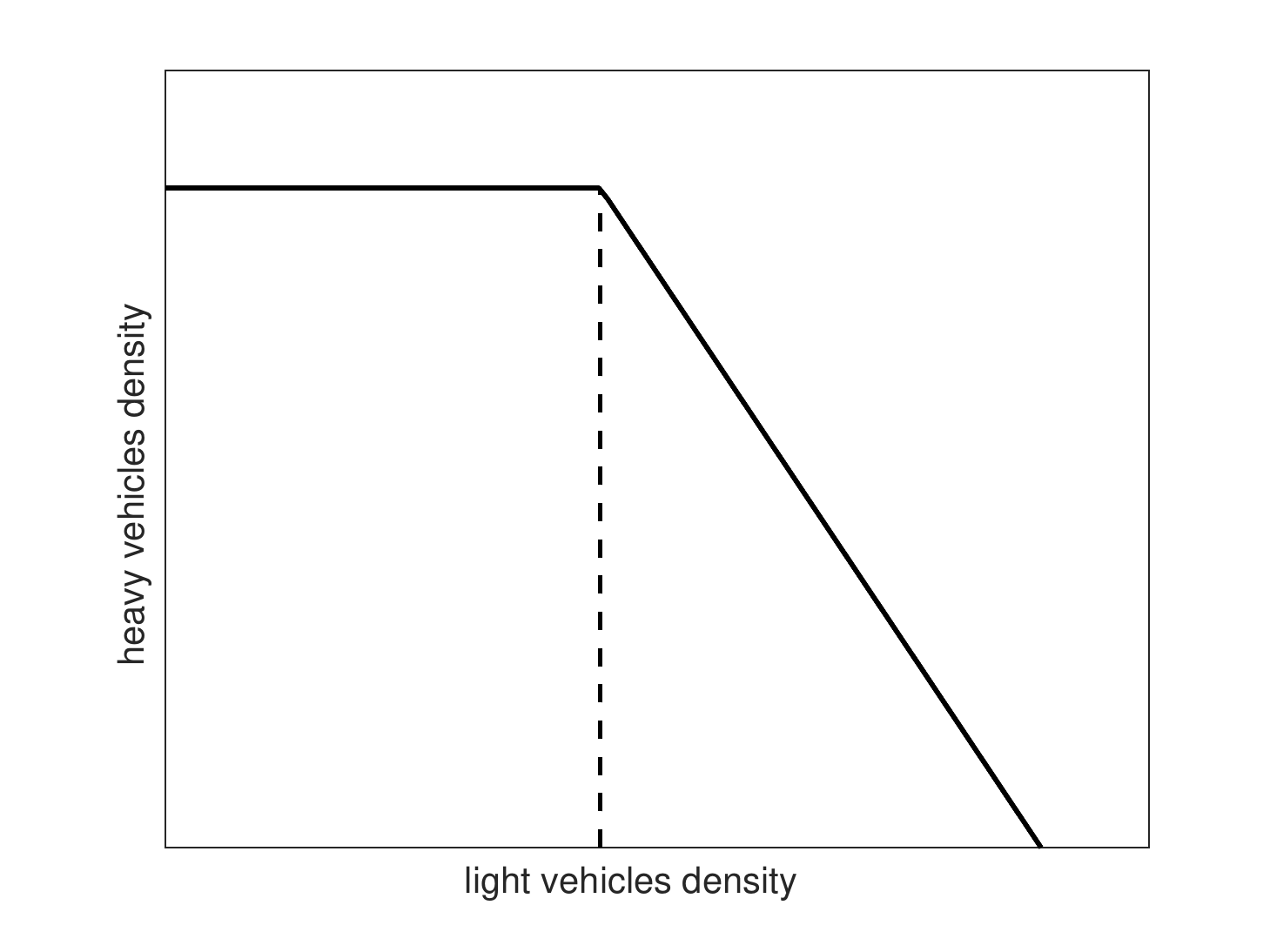}
				\put(17,30){partial}  \put(17,25){coupling}
				\put(27,15){$\Domain_1$}   \put(57,15){$\Domain_2$}
				\put(52,30){full}  \put(52,25){coupling}
				\put(44,20){\rotatebox{90}{\small transition level}}
				\put(80,5){$\rhomax$}
				\put(5,60){$\mumax$}
			\end{overpic}
			\caption{Domains $\Domain_1$ and $\Domain_2$ of the macroscopic model \eqref{eq:model}.}
			\label{fig:dominio_MaMa}
		\end{figure}
		In this phase we assume that cars are mainly in the fast lane and do not affect the trucks dynamics. Trucks are then independent from cars.
		 
		 For trucks we choose a triangular fundamental diagram with
		\begin{equation}\label{eq:HeavyVel_PC}
		v_\pes(\rhoP) = V_\pes^\text{max} \quad \mbox{ for all } \quad \rhoP\leq \sigma_\pes,
		\end{equation}
		 where 
		 $V_\pes^\text{max}$ is the maximum speed of trucks, see Fig.\ \ref{fig:Vel_PCphase}\subfig{b}.
		 
		 Cars do not interfere with trucks but adapt their dynamics to the presence of them. 
		 Also for cars we choose a (family of) triangular fundamental diagrams, see Fig.\ \ref{fig:Vel_PCphase}\subfig{a}.
		 \rev{Specifically, we set
	\begin{equation}\label{eq:LightVel}
		v_\leg(\rhoL,\rhoP) = 
		\left\{\begin{array}{lcl}
			 \VLstar(\rhoP) & \mbox{ if } & \rhoL\leq\sigma_\leg(\rhoP),
			\bigskip\\
			\displaystyle 
			\frac{\VLstar(\rhoP)\ \sigma_\leg(\rhoP)} {\rhoLstar(\rhoP)-\sigma_\leg(\rhoP)}
			\left(\frac{\rhoLstar(\rhoP)}{\rhoL}-1\right)& \mbox{ if } & \sigma_\leg(\rhoP) < \rhoL\leq\rhomax-\mumax/\beta,
		\end{array}\right.
	\end{equation}
	where $\VLstar(\rhoP)$ is the maximum speed of cars given the truck density. We also define $\VLstar(0)=V_\leg^\text{max}$ as the maximum speed of cars in absence of trucks. 
	Then, $\VLstar(\rhoP)\geq 0$ and $\sigma_\leg(\rhoP)\geq 0$ are continuous linear decreasing functions of $\rhoP$.	
		 }
		 	
		For $(\rhoL,\rhoP) \in \mathcal{D}_1$ the model \eqref{eq:model} then becomes
		\begin{equation}\label{eq:model_PC}
		\begin{cases}
		\partial_t\rhoL + \partial_x f_\leg(\rhoL,\rhoP)=\, 0 
		\smallskip\\
		\partial_t \rhoP + \partial_xf_\pes(\rhoP)=\, 0
		\end{cases}
		\end{equation}
		where $f_\leg(\rhoL,\rhoP)=\rhoL v_\leg(\rhoL,\rhoP)$ and $f_\pes(\rhoP)=\rhoP v_\pes(\rhoP)$ as described in Fig.\ \ref{fig:Vel_PCphase}. 
		\begin{figure}[h!]
			\subfig{a}
			\begin{overpic}
			[width=6cm]{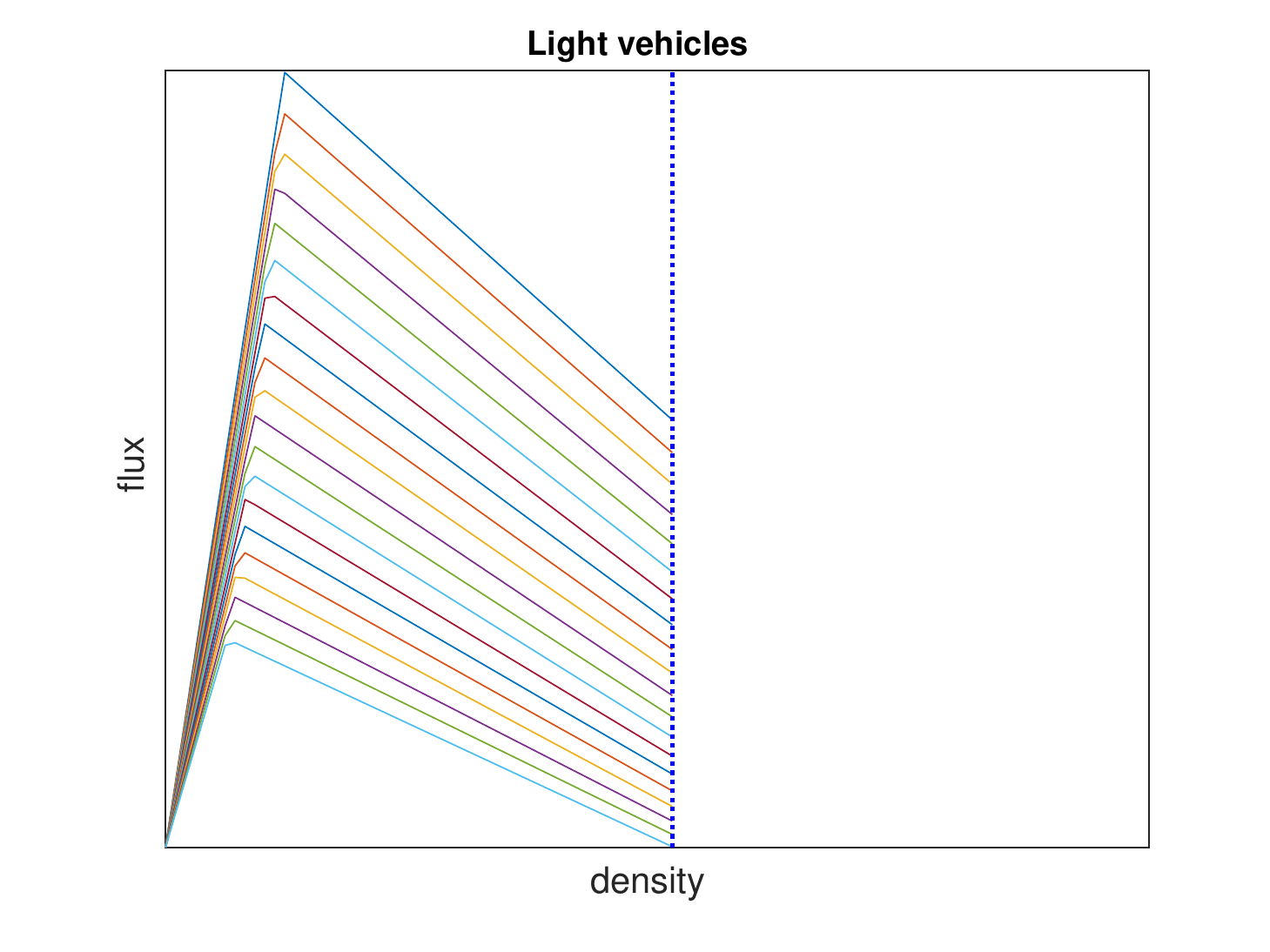}
			\put(50,47){\rotatebox{90}{\tiny transition level}}
			\end{overpic}
			\subfig{b}\includegraphics[width=6cm]{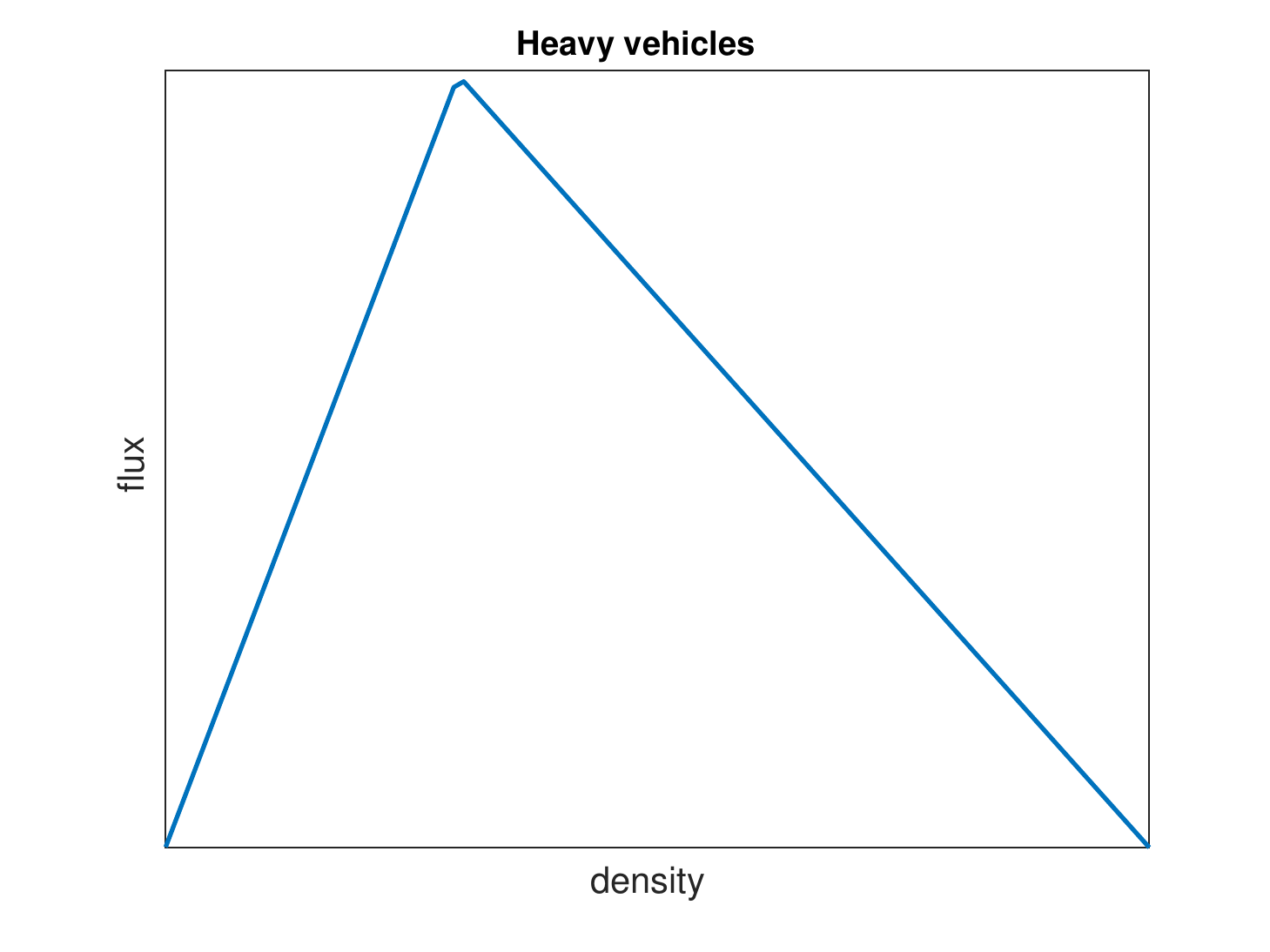}
			\caption{Fundamental diagrams of the macroscopic model in the partial-coupling phase, i.e.\ $(\rhoL,\rhoP)\in\Domain_1$.}
			\label{fig:Vel_PCphase}
		\end{figure} 
		\item The \emph{full-coupling phase} is in place when $(\rhoL,\rhoP)\in\Domain_2:=\Domain\backslash\Domain_1$,
		see Fig.\ \ref{fig:dominio_MaMa}.
		In this case, we assume that cars are too much to find it convenient to be confined in the fast lane. For this reason they invade the slow lane, thus influencing the dynamics of trucks. The two equations in system \eqref{eq:model} are then fully coupled. 
		
		As before, we choose for both classes a family of triangular fundamental diagrams which extend by continuity those defined in $\Domain_1$, as shown in Fig.\ \ref{fig:Vel_FCphase}.
		\begin{figure}[h!]
			\subfig{a}
			\begin{overpic}
				[width=6cm]{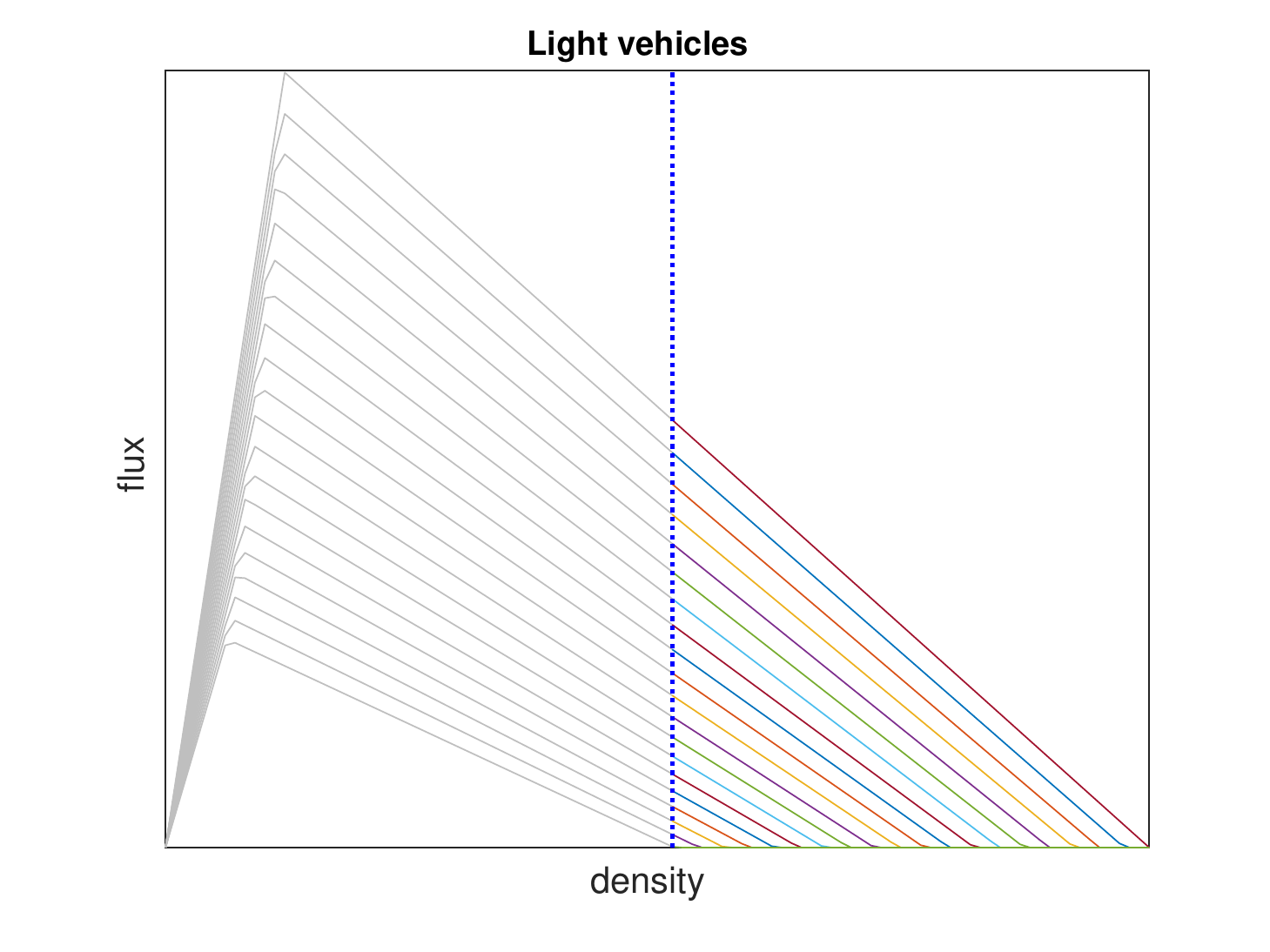}
				\put(50,47){\rotatebox{90}{\tiny transition level}}
			\end{overpic}
			\subfig{b}\includegraphics[width=6cm]{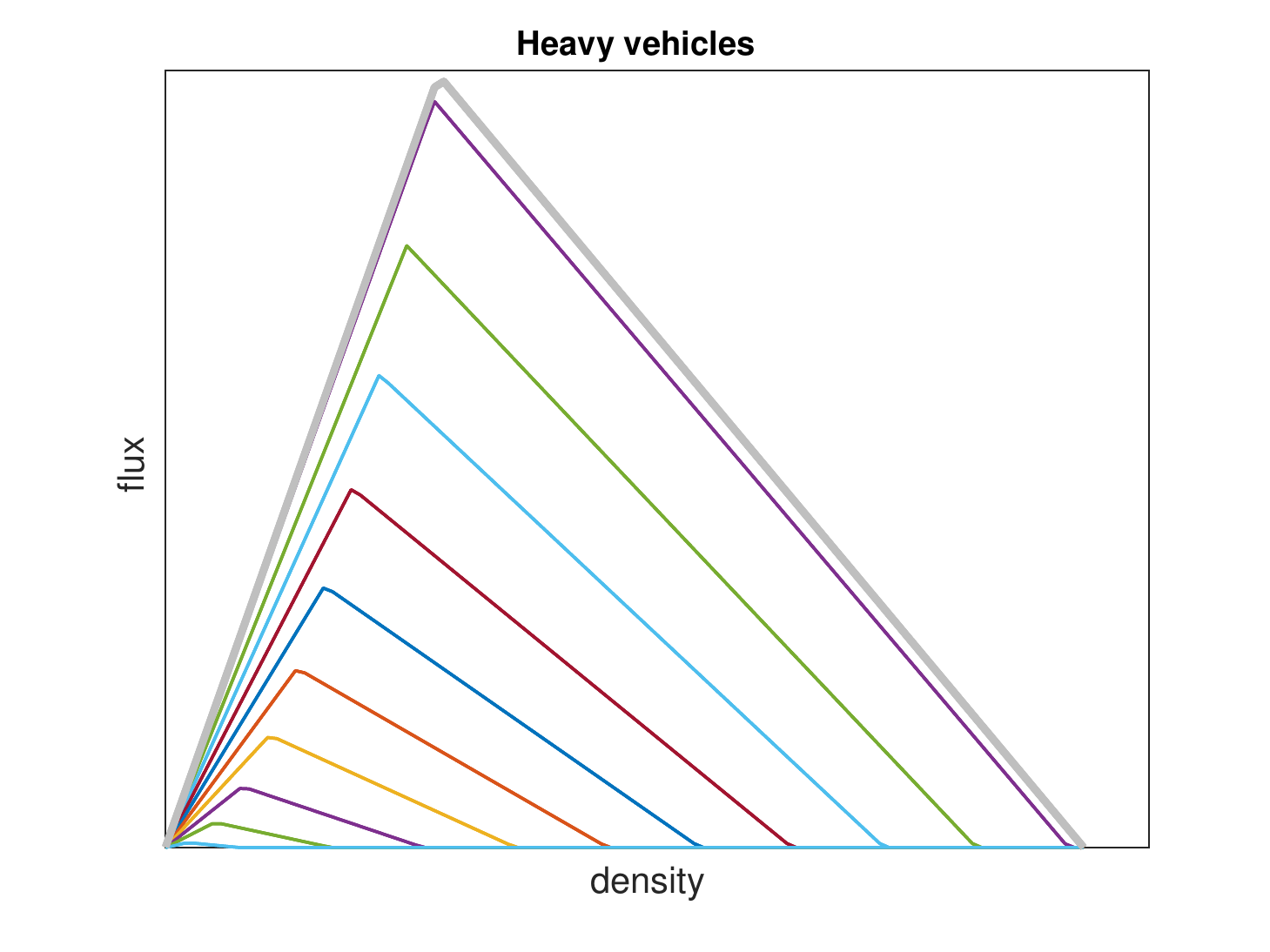}
			\caption{Fundamental diagrams of the macroscopic model in the fully-coupling phase, i.e.\ $(\rhoL,\rhoP)\in\Domain_2$.}
			\label{fig:Vel_FCphase}
		\end{figure} 
	\end{itemize}
\medskip 

	\noindent We define the \emph{transition level} the threshold density of light vehicles which act as interface between the two phases, see Fig.\ \ref{fig:dominio_MaMa}. 
	In our setting trucks are confined in one of the two available lanes, then the transition level is equal to $\rhomax-\mumax/\beta=\rhomax/2$.
	
	Note also that the fundamental diagrams we use in this work verify all the properties (L1)--(L5) and (H1)--(H5).

	
	\subsection{Multi-scale model}\label{sec:mima_model}
	
	In this section we describe the multi-scale model. Here cars are described by a first-order LWR model of type \eqref{LWR} and trucks are described by a second-order microscopic Follow-the-Leader model of type \eqref{FtL}. 
	Let us describe the microscopic model first, dropping for the moment the coupling with light vehicles.
	
	\subsubsection{Microscopic model for heavy vehicles}
	
	The microscopic model is the one presented in \cite{cristiani2019DCDS-B}, which is, in turn, inspired by the model originally proposed by Zhao and Zhang in \cite{zhao2017TRB}.
	
	In the following we denote by $\Delta_k$ the gap between truck $k$ and truck $k+1$ at any time $t$,
	$$
	\Delta_k(t):=X_{k+1}(t)-X_k(t).
	$$
	It is plain that this gap is inversely proportional to the density of heavy vehicles.
	
	We define in \eqref{FtL}
	\begin{equation}\label{def:Acinese}
		A\big(X_k,X_{k+1},V_k,V_{k+1}\big)=
		\left\{
		\begin{array}{ll}
		\frac{1}{\tau_\text{acc}}\Big(v^\textsc{zz}(\Delta_k)-V_k\Big), & \text{ if\ \ } v^\textsc{zz}(\Delta_k)\geq V_k \\
		\frac{1}{\tau_\text{dec}}\Big(v^\textsc{zz}(\Delta_k)-V_k\Big), & \text{ if\ \ } v^\textsc{zz}(\Delta_k)<V_k 
		\end{array}
		\right.
	\end{equation}
	where the function $v^\textsc{zz}$ represents the equilibrium velocity all drivers tend to, and depends on the gap $\Delta_k$.
	Parameters $\tau_\text{acc},\tau_\text{dec}>0$ are the relaxation times as usual, differentiated for the acceleration and the deceleration phase. Diversifying the relaxation times appeared to be crucial to fit real data.
	
\begin{figure}[h!]
	\begin{overpic}
		[width=6cm]{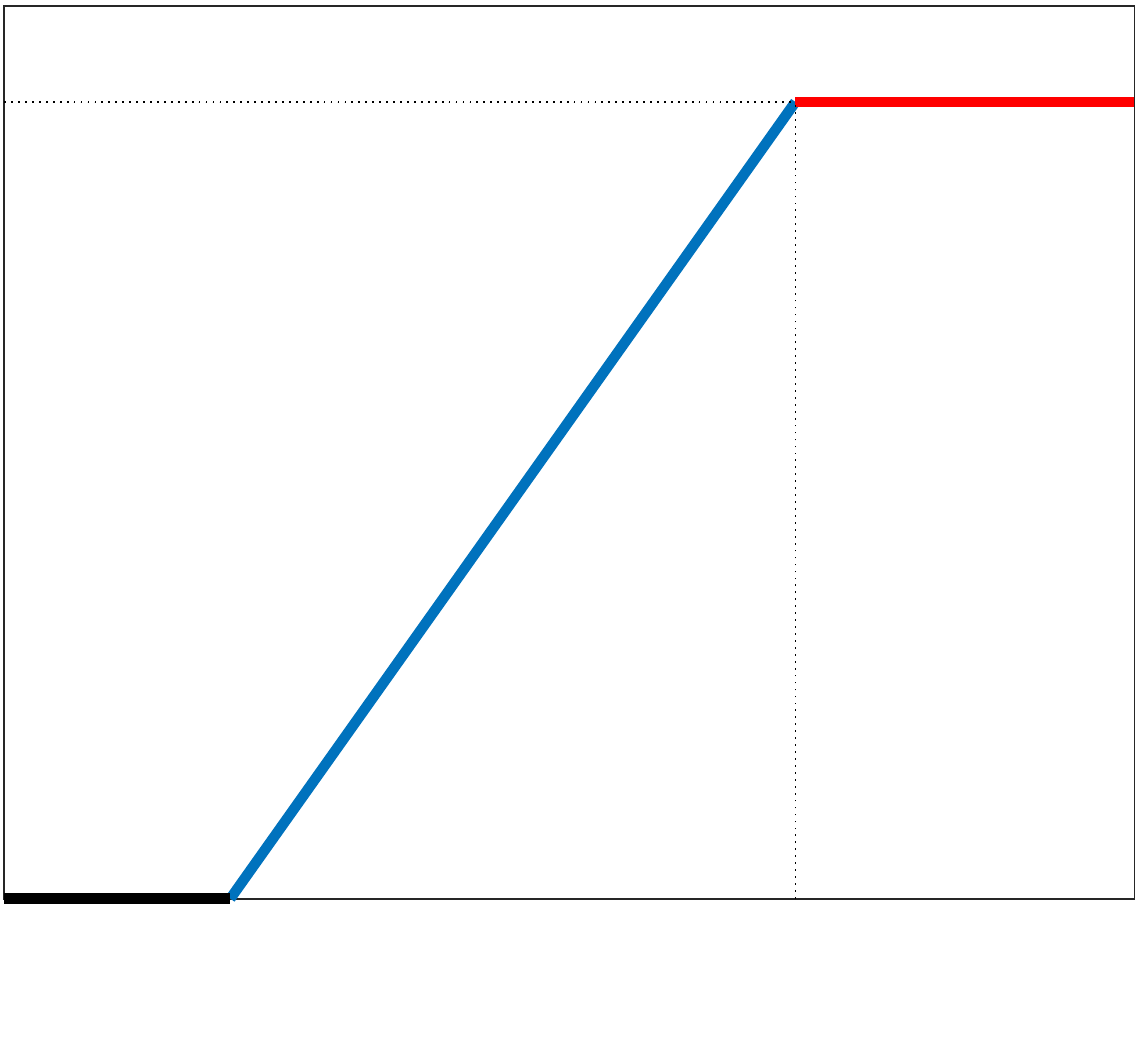}
		\put(15,5){$\Delta_\text{close}$}
		\put(65,5){$\Delta_\text{far}$}
		\put(-15,80){$V_\pes^\text{max}$}
	\end{overpic}  
	\caption{The shape of the velocity function $v^\textsc{zz}(\Delta)$ defined in \eqref{def:vZZ}.}
	\label{fig:vZZ}	
	\end{figure}

	The velocity function $v^\textsc{zz}$ is defined by
	
	\begin{equation}\label{def:vZZ}
		v^\textsc{zz}(\Delta):=\left\{
		\begin{array}{ll}
			0, & \text{ if\ \ }\Delta\leq \Delta_\text{close} \\ [2mm]
			\frac{V_\textsc{h}^\text{max}}{\Delta_\text{far}-\Delta_\text{close}}(\Delta-\Delta_\text{close}), & \text{ if\ \ } \Delta_\text{close}<\Delta< \Delta_\text{far}\\ [3mm]
			V_\pes^\text{max}, & \text{ if\ \ }\Delta\geq \Delta_\text{far}
		\end{array}
		\right.
	\end{equation}
	where $\Delta_\text{close}$, $\Delta_\text{far}$, $V_\pes^\text{max}$ are positive parameters, see Fig.\ \ref{fig:vZZ}. 
	%
	The plateau in $\Delta\in[0,\Delta_\text{close}]$ is crucial for correctly reproducing stop \& go waves. Indeed, once the relaxation times $\tau_\text{acc}$, $\tau_\text{dec}$ are fixed, the capability of the model to trigger stop \& go waves is ruled precisely by $\Delta_\text{close}$.
	
	\subsubsection{Full model}
	First of all, given a parameter $\delta>0$, we assume that cars located at $x$ are influenced by a truck iff the distance between the truck and $x$ is less than $\delta$. 
	We denote the number of trucks falling in the road interval $[x-\delta,x+\delta)$ at any time $t$ by
	\begin{equation}\label{def:Nh}
		N^\delta_\pes(x,t):=\#\{k:X_k(t)\in[x-\delta,x+\delta)\}.
	\end{equation} 
	
	Second, we denote by $\rhoL$ the density of light vehicles and $v_\leg$ their velocity. To couple the dynamics of the two classes we assume that $v_\leg$ depends on both $\rhoL$ (as in the classical LWR model) and $N^\delta_\pes$. As usual we assume that $v_\leg$ is decreasing with respect to both arguments.
	
	Finally, we couple the dynamics of heavy vehicles with those of light vehicles. 
	The interaction is obtained by introducing the dependence on $\rhoL$ in the parameters $\Delta_\text{close}$  and $\Delta_\text{far} $.
	More precisely, we introduce the increasing functions 
	$\Delta_\text{close}  =\Delta_\text{close}(\rhoL)$ and $\Delta_\text{far}=\Delta_\text{far}(\rhoL)$ 
	and we denote by $A^\textsc{c}=A^\textsc{c}(X_k,X_{k+1},V_k,V_{k+1},\rhoL)$ the coupled acceleration defined as $A$ in \eqref{def:Acinese}--\eqref{def:vZZ}, with the new dependence on $\rhoL$.
		
	We are now ready to present the fully coupled multi-scale model which reads as
	\begin{equation}\label{multiscalemodel}
	\left\{
	\begin{array}{l}
	\left\{
	\begin{array}{l}
	\dot X_k=V_k \\ [2mm]
	\dot V_k=A^\textsc{c}(X_k,X_{k+1},V_k,V_{k+1},\rhoL) 
	\end{array}, \qquad  k=1,\ldots,N-1 
	\right. 
	\\ [5mm]
	\partial_t\rhoL+\partial_x\Big(\rhoL v_\leg(\rhoL,N^\delta_\pes)\Big)=0,\qquad x\in\R,\quad t>0.
	\end{array}
	\right. 
	\end{equation}
	To be coherent with our modeling assumptions, the functions $\Delta_\text{close}$ and $\Delta_\text{far}$ are constant for cars densities below the transition level, i.e.\ $\rhoL\leq\rhomax/2$. In this case the dynamics of trucks is independent from those of cars. 
	Conversely, for $\rhoL>\rhomax/2$, we assume that the distances $\Delta_\text{close}$ and $\Delta_\text{far}$ increase linearly with respect to the average number of cars which are positioned between two trucks. This number can be easily computed considering the average number of cars in a road segment of length $\ell$ (equal to $\rhoL\ell$) and the number of trucks in the same road segment (assuming that all vehicles are uniformly distributed).
	We were unable to precisely calibrate the shape of the functions $\Delta_\text{close}$ and $\Delta_\text{far}$ from real data because it happens rarely that many cars are found between trucks: indeed, trucks tend to ``push'' cars in the fast lane rather than reacting to their presence.

	\subsection{Extension of the models to general road networks}\label{sec:generalizations}
	In order to perform a complete simulation on a generic network of highways, some important generalizations are needed.
	
	\subsubsection{Any number of lanes}\label{sec:generalizations-lanes}
	Highways have often more than two lanes. Consider a road with $n$ lanes of which $n_\pes$ can be occupied by trucks. To allow the creeping phenomenon, we assume that $n_\pes<n$, which corresponds to $\rhomax\lL-\mumax\lP>0$ in terms of space occupied, cf.\ \eqref{eq:DensCondstretta}.

	In the macroscopic approach the model is easy generalized. Fundamental diagrams are modified in such a way that trucks start interacting with cars when the density of cars becomes greater than $\frac{n_\pes}{n} \rhomax$. 
	
	In the microscopic model instead, an important modification is needed if $n_\pes>1$. Indeed, in this case trucks can overtake and the microscopic model must be able to handle this. Typically, some new parameters are introduced in order to establish when a truck decides to overtake and if the truck can actually overtake, considering suitable safety constraints. 
	From the computational point of view, an additional difficulty arises when one has to find the truck in front of any other truck, since the ordering is lost whenever a truck overtakes.
	To make the search of the preceding vehicle computationally feasible, one can keep track, in a specific list, of all trucks located in each numerical cell, and then update the list whenever a truck leaves or enter the cell. 
	
	\subsubsection{Junctions}\label{sec:generalizations-junctions}
	In order to perform a full simulation on a network of highways both theoretical and numerical treatment of junctions are needed. Typically highways have not roundabouts, traffic lights or complex junctions, so we can limit ourself to handle simple merge (2 incoming roads and 1 outgoing road) and diverge (1 incoming road and 2 outgoing roads). 
	We have adopted the approach detailed in \cite{briani2018CMS}, in which the dynamics is reformulated along paths and junctions ``disappears''. 
	The price to pay is that the number of equations is multiplied by the number of possible paths the drivers can follow at junctions. In both merge and diverge, we have only two possible paths: for example, in the case of the diverge, one can choose among the first and the second outgoing road, while in a merge one can decide to come from the first or the second incoming road.
	
	Following this approach in the macroscopic model, the densities of each class of vehicles are split around every junction, ending up with a system of four conservation laws (two paths for each of the two classes of vehicles) with discontinuous flux. 
	After the junctions, densities are gathered together again and the two-equation system \eqref{eq:model} is restored.  
	
	In the multi-scale model instead, the path-based approach is applied only for cars dynamics since managing trucks is much simpler. Indeed, in the microscopic model one can just move vehicles from one road to another on the basis of their destination, see \cite{cristiani2016NHM}.
	Unfortunately, the ordering of trucks is lost every time a change of road takes place. 
	In order to reduce the computational effort need for the computation of the preceding truck of every truck, the same solution proposed in Sect.\ \ref{sec:generalizations-lanes} can be applied.
	

	\section{Numerical approximation and calibration}\label{sec:implementation}
	In this section we describe how the models introduced above can be actually implemented. First, we briefly recall the numerical methods we have adopted, then we describe how we have used real data to set the models' parameters. 
		
		\subsection{Macroscopic model}\label{sec:numerics_mama}
		For the numerical approximation of the macroscopic model \eqref{eq:model} we employ the extension of the \textit{cell transmission model} (CTM) to heterogeneous multi-class model proposed in \cite{fan2015SIAP}. 
		Let $\dx$ and $\dt$ be the space and time step respectively, and $(\rhoL^{n,i},\rhoP^{n,i})$ the traffic densities in the $i$th cell at the $n$th time step. 
		The finite volume numerical scheme reads
			\begin{subequations}\label{eq:num2classes}
			\begin{empheq}[left=\empheqlbrace]{align}
			\hspace{3mm} \rhoL^{n+1,i}=\rhoL^{n,i} + \frac{\dt}{\dx}\left(\FL^{n,i-1/2} -\FL^{n,i+1/2}\right) \label{eq:num2classes_cars}\\ 
			\hspace{3mm} \rhoP^{n+1,i}=\rhoP^{n,i} + \frac{\dt}{\dx}\left(\FP^{n,i-1/2} -\FP^{n,i+1/2}\right) \label{eq:num2classes_trucks}
			\end{empheq}
		\end{subequations}
		where
		\begin{equation}\label{eq:num_flux2classes_cars}
		\mathcal{F}^{n,i+1/2}_{\leg} := \min\Big\{S_{\leg}(\rhoL^{n,i},\rhoP^{n,i}),R_{\leg}(\rhoL^{n,i+1},\rhoP^{n,i+1})\Big\},
		\end{equation}
		\begin{equation}\label{eq:num_flux2classes_trucks}
		\mathcal{F}^{n,i+1/2}_{\pes} := \min\Big\{S_{\pes}(\rhoL^{n,i},\rhoP^{n,i}),R_{\pes}(\rhoL^{n,i+1},\rhoP^{n,i+1})\Big\},
		\end{equation}
		and $(\SL, \RL)$, $(\SP, \RP)$ represent the sending and receiving functions of the two vehicles classes respectively, defined by
		\begin{equation}\label{eq:SR1}
		\begin{array}{c}
		\SL(\rhoL,\rhoP):=\left\{\begin{array}{lr}
		f_\leg(\rhoL,\rhoP), &\mbox{if } \rhoL \leq \sigma_\leg(\rhoP),
		\smallskip\\
		f_\leg(\sigma_\leg(\rhoP),\rhoP), &\mbox{if } \rhoL > \sigma_\leg(\rhoP),
		\end{array}\right.
		\medskip\\
		\RL(\rhoL,\rhoP):=\left\{\begin{array}{lr}
		f_\leg(\sigma_\leg(\rhoP),\rhoP), &\mbox{if } \rhoL \leq \sigma_\leg(\rhoP),
		\smallskip\\
		f_\leg(\rhoL,\rhoP), &\mbox{if } \rhoL > \sigma_\leg(\rhoP),
		\end{array}\right.
		\end{array}
		\end{equation}
		and similarly for $(\SP,\RP)$.
		
		
		The numerical grid is chosen as $\Delta x$=100 m and $\Delta t$=2.6 s. 
		\rev{The choice of the space step comes from the fact that the company Autovie Venete finds such granularity convenient to share traffic information to drivers, while the time step is dictated by the CFL condition.} 
		
		\medskip
		
\begin{figure}[h!]
		\subfig{a}\includegraphics[width=6cm]{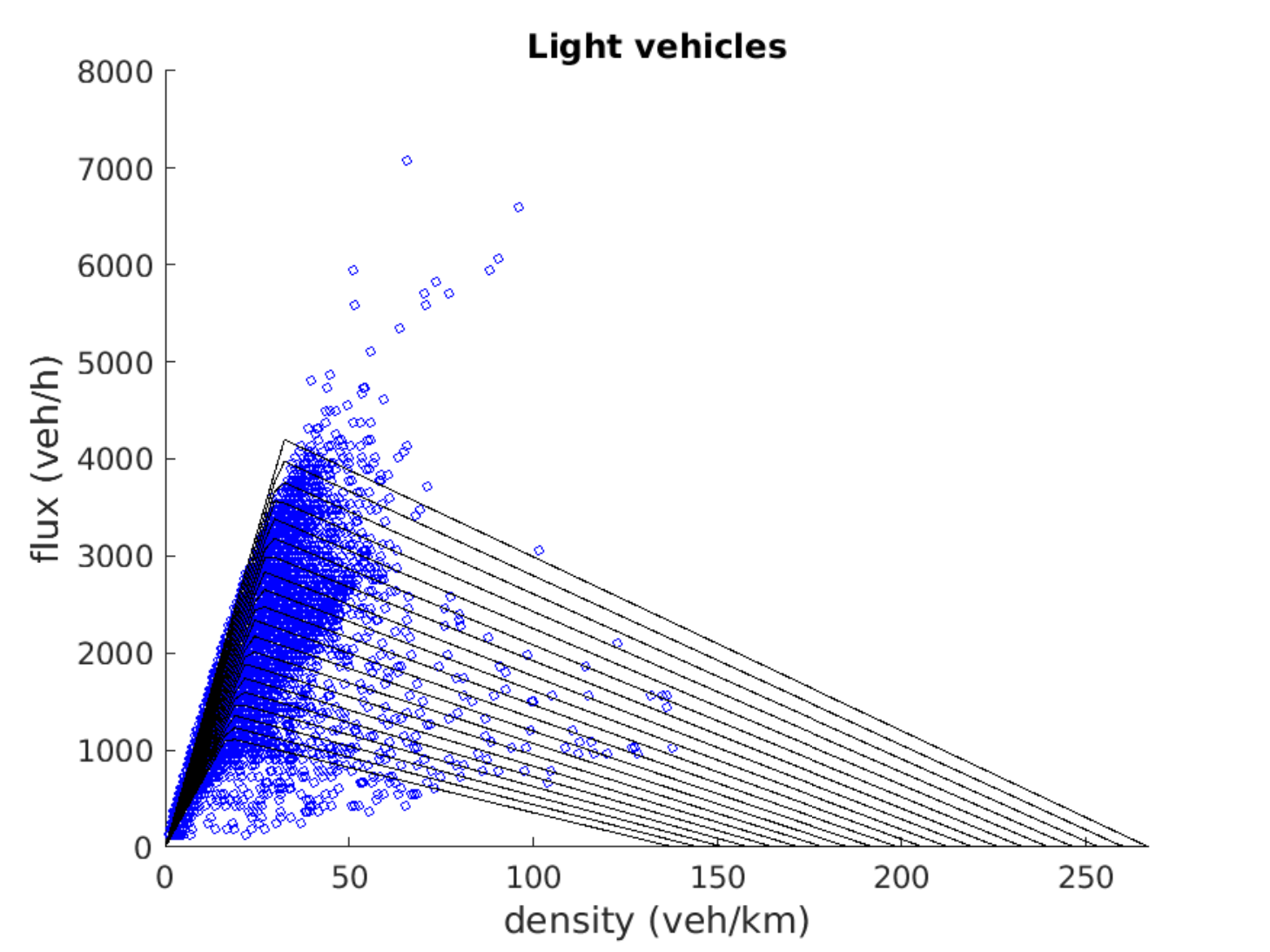}
		\subfig{b}\includegraphics[width=6cm]{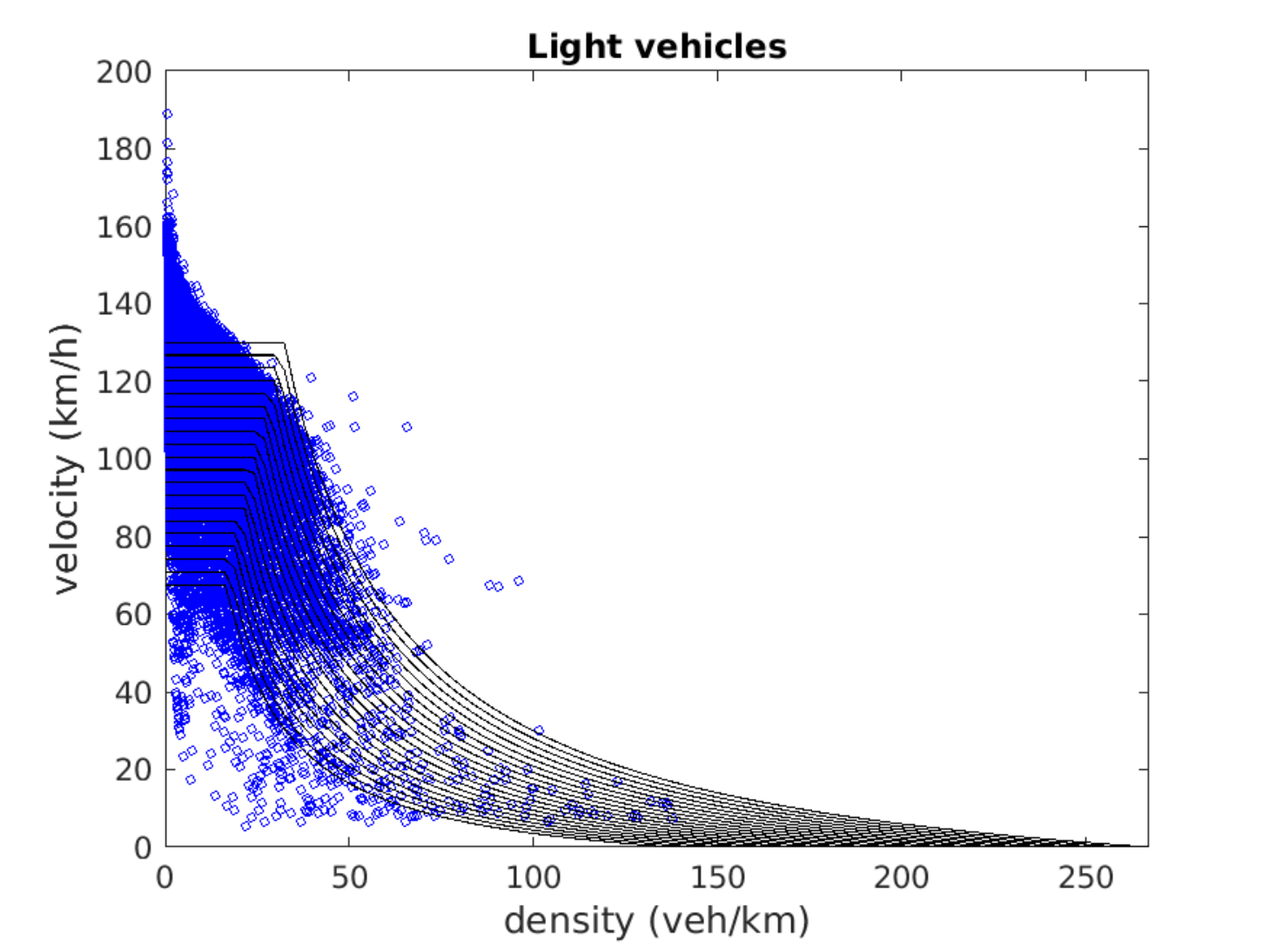}
		\caption{\subfig{a} Flux-density and \subfig{b} velocity-density relationships for cars with real data superimposed.}
		\label{fig:Calib_MaMa_cars}
	\end{figure}

\begin{figure}[h!]
		\subfig{a}\includegraphics[width=6cm]{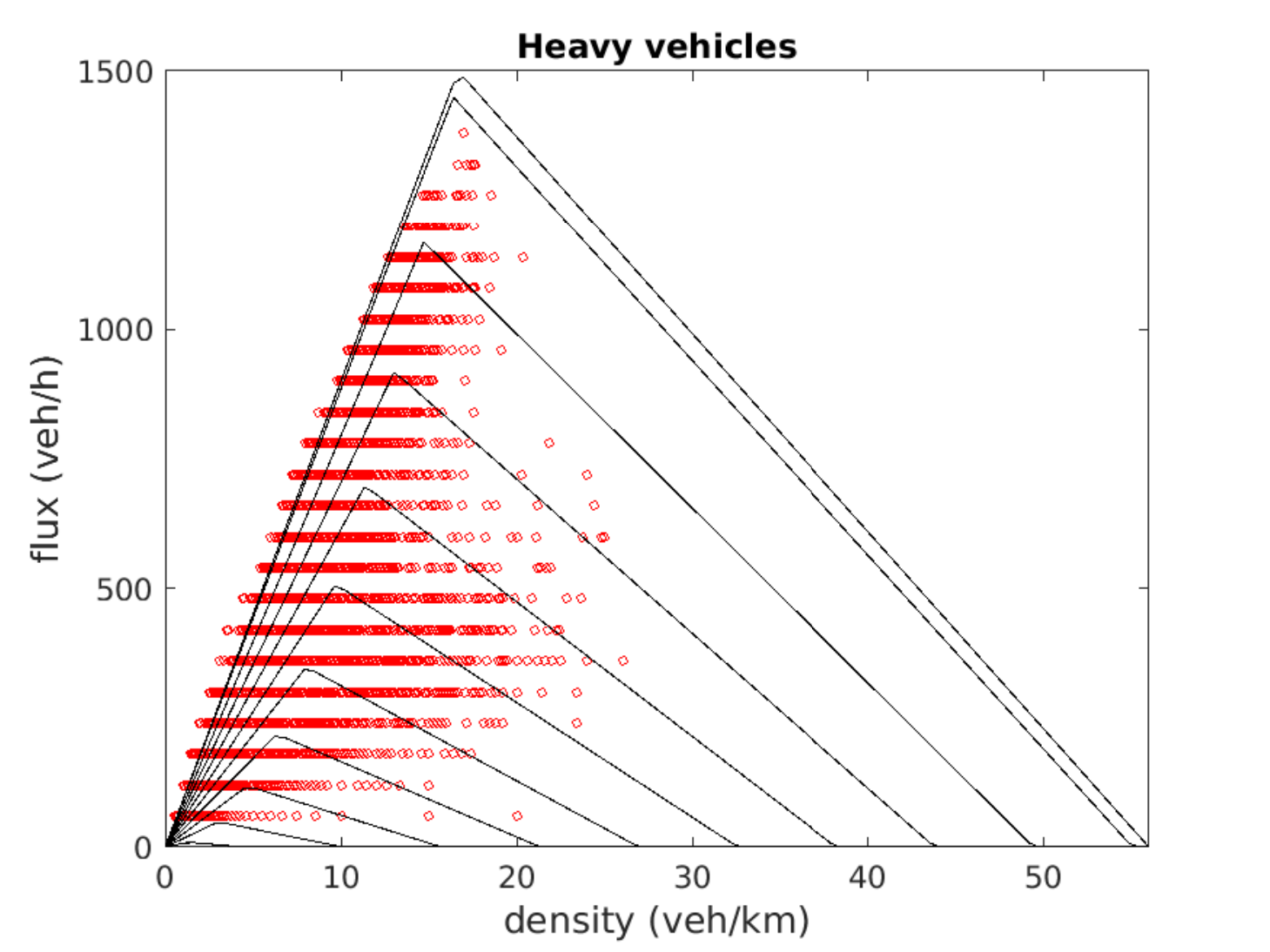}
		\subfig{b}\includegraphics[width=6cm]{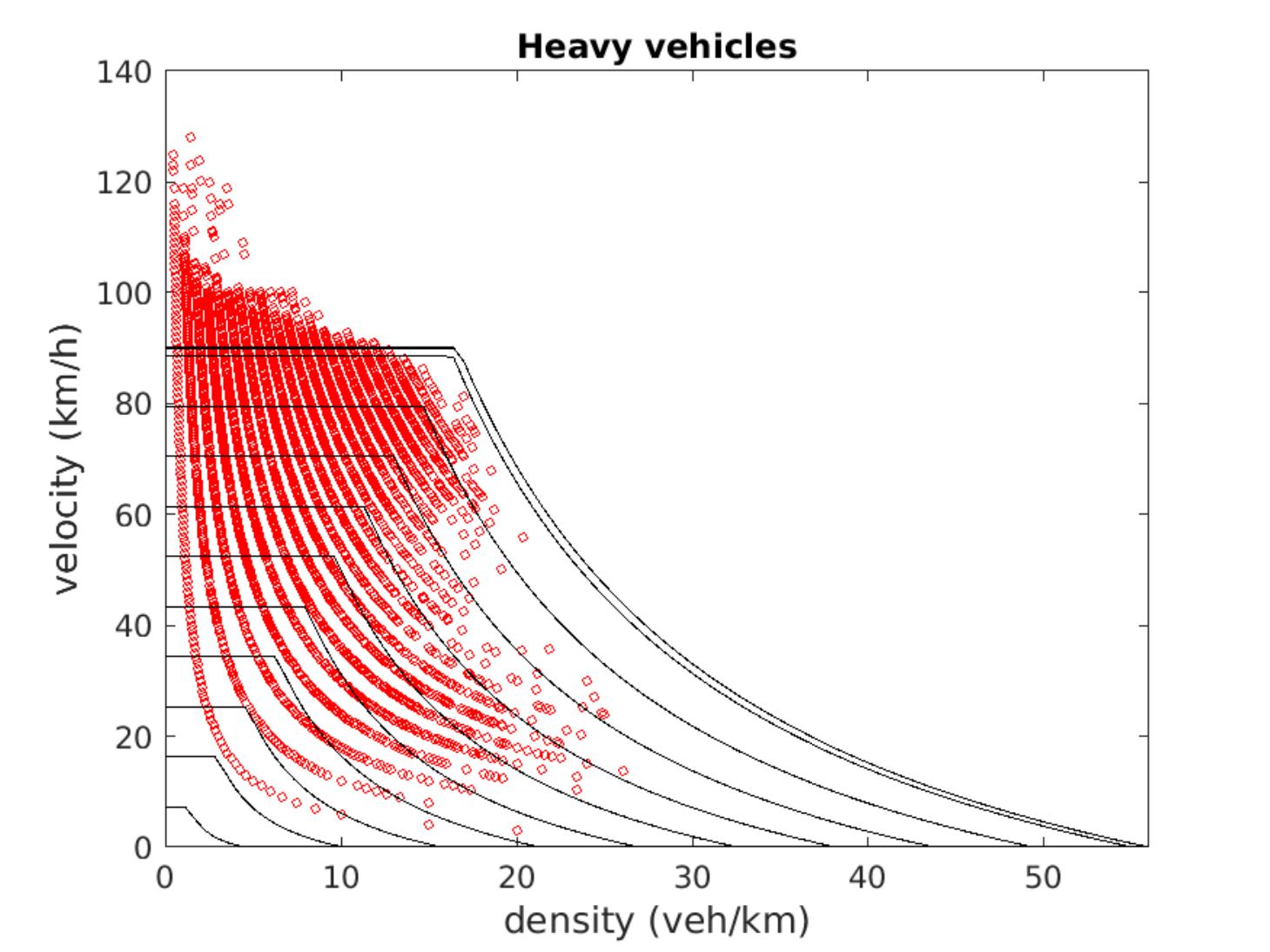}
		\caption{\subfig{a} Flux-density and \subfig{b} velocity-density relationships for trucks with real data superimposed.}
		\label{fig:Calib_MaMa_trucks}
	\end{figure}

		Calibration of the fundamental diagrams was performed by fitting real data. We used all data measured in 2019 by one fixed sensor located near Cessalto, see Figs.\ \ref{fig:Calib_MaMa_cars}-\ref{fig:Calib_MaMa_trucks}.
		Note that for high densities, the velocities drop rapidly to zero. Since we have no data for completely stationary vehicles under the sensor, we are not able to reconstruct data on high traffic density.
		For this reason, the maximal densities $\rhomax$ and $\mumax$ are estimated by simply computing the ratio between the number of available lanes for the class and the average length of vehicles of that class, see Eq.\ \eqref{def_rhomax}.

		Model parameters are summarized in Table \ref{tab:Calib_param_MaMa}. All functions which rule the dependence of $\rho, v, f$ on the density of the other class are linear.
	\begin{specialtable}[h!]
	\caption{Parameters for the macroscopic model.}
	\label{tab:Calib_param_MaMa}
	\small
	\begin{tabular}{lll}\toprule
		& Light vehicles & Heavy vehicles  \\ \midrule
		veh.\ length + safety dist.\ (km) &  $7.5\cdot 10^{-3}$ & $18\cdot 10^{-3}$  \\ \hline
		max max density (veh/km) & $\rhoLstar(0) = 267$ & $\rhoPstar(0) = 56$ \\ \hline
		min max density (veh/km)  & $\rhoLstar(\mumax) = 133$ & $\rhoPstar(\rhomax) = 0$ \\ \hline
		max max speed (km/h)  & $v_\leg(0,0) = 130$ & $v_\pes(0,0) = 90$\\ \hline
		min max speed (km/h) & $v_\leg(0,\mumax) = 65$ & $v_\pes(\rhomax,0) = 0$\\ \hline
		max max flux (veh/h) & $f_\leg(\sigma_\leg(0),0) = 4200$ & $f_\pes(0,\sigma_\pes(0)) = 1500$ \\ \hline
		min max flux (veh/h) & $f_\leg(\sigma_\leg(\mumax),\mumax) = 1200$ & $f_\pes(\rhomax,\sigma_\pes(\rhomax)) = 0$
		\\\bottomrule 
	\end{tabular}
	\end{specialtable}
	
	
	\subsection{Multi-scale model}
	For the numerical approximation of the macroscopic part of the multi-scale model \eqref{multiscalemodel} we employ again the scheme \eqref{eq:num2classes_cars}, where $N_\pes^\delta$ plays the role of $\rho_\pes$ in the obvious manner.
	Numerical grid is chosen as $\Delta x$=100 m and $\Delta t$=2 s.
	
	The dependence of the flux on $N_\pes^\delta$ can generate some issues. For example, consider the case of no trucks and  a car density $\hat\rho_\leg$ close to $\rhomax$. 
	When a truck enters the road, the maximal density allowed in the cell occupied by the truck drops to $\rhoLstar(N^\delta_\pes)$ according to \eqref{eq:RhoMax}. 
	Now, if $\rhoLstar(N^\delta_\pes)<\hat\rho_\leg$, the current density $\hat\rho_\leg$ is found not to be compatible with the new maximal density. 
	Although the entering truck perceives the cars, it is not guaranteed that the compatibility with the maximal density is respected at any time. 
	To avoid this problem, trucks must be prevented from entering cells if the new maximal density caused by the presence of the truck itself is not compatible with current traffic conditions.   

	For the numerical approximation of the microscopic part we used a standard Euler scheme with a time step $\delta t$= 0.1 s. Note that this time step is much smaller than the time step $\Delta t$ used for the Godunov scheme, meaning that the updates of the trucks and cars are asynchronous.
	
	\medskip
	
	Regarding the parameters, the macroscopic part of the model is treated as in Section \ref{sec:numerics_mama} (Table \ref{tab:Calib_param_MaMa}). 
	For the microscopic model, some parameters are easily calibrated by using real data and considering physical constraints. 
	For example $V_\pes^\text{max}$ was defined as in the macroscopic model. 
	$\Delta_\text{far}$ was set in order to guarantee that trucks do not collide even in the event that a truck suddenly brakes with full power until it stops (note that our model allows in principle collisions since deceleration is bounded).
	$\Delta_\text{close}$, instead, was set to a distance which guarantees to catch the maximal observed density of trucks. In other words, when a queue of trucks is formed, the model predicts the correct maximal density.
	
\begin{figure}[h!]
		\includegraphics[width=11cm]{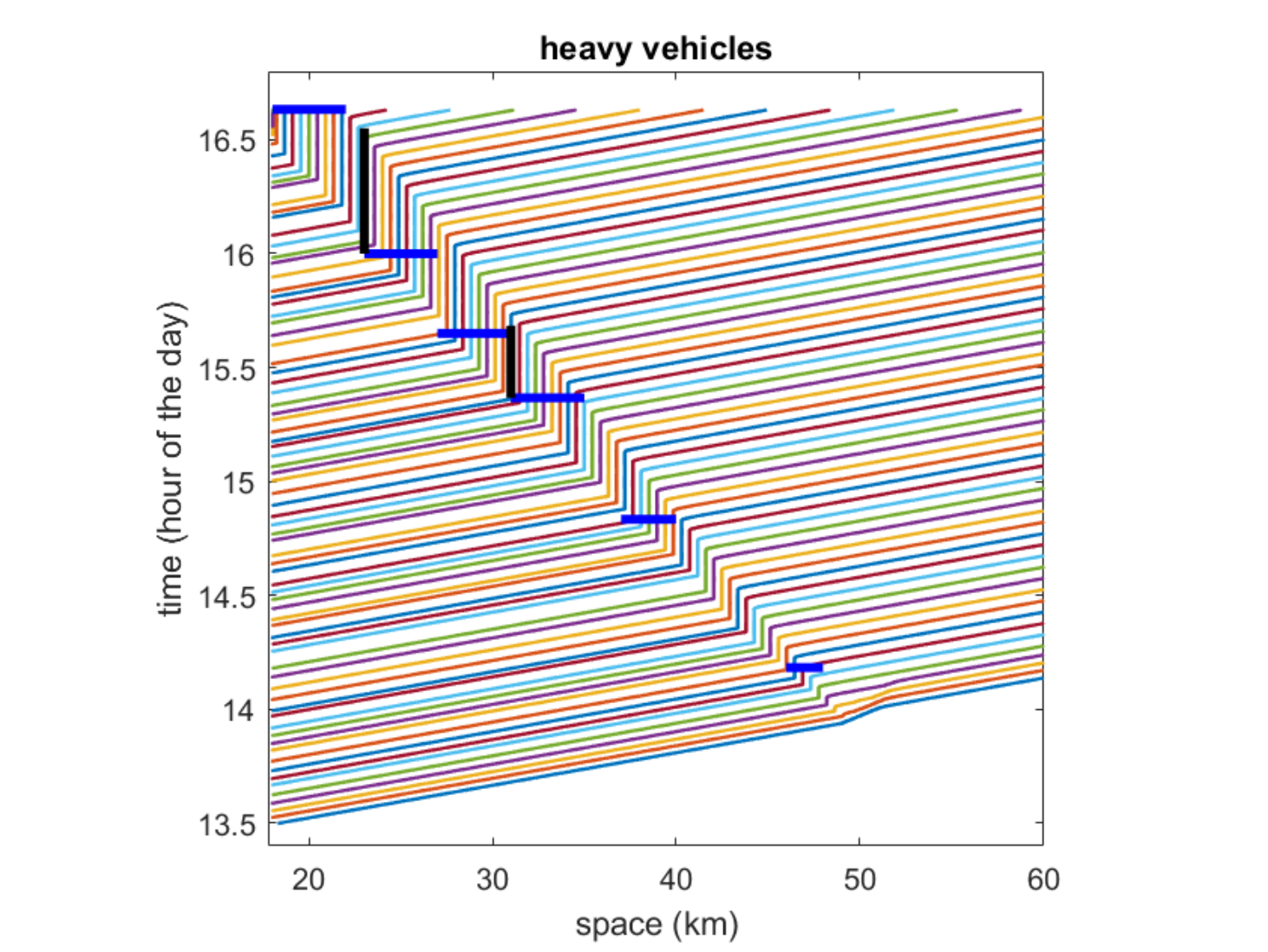}
		\caption{\rev{Simulated trajectories obtained with real inflow data as left boundary conditions} (not all vehicles are plotted for visualization purposes). Horizontal blue lines and vertical black lines indicate, respectively, the position and the duration of the real queue as measured on the field.}
		\label{fig:calibrazione_tau}
	\end{figure}

	Parameters $\tau_\text{acc}$ and $\tau_\text{dec}$ are instead more difficult to calibrate since they are not easily measurable.  
	For those values, we considered a real stop \& go wave observed by the company staff on June 12, 2017, generated by the slowdown of a truck near a bottleneck. \rev{The initial perturbation (slowdown) was amplified and in short time generated a queue which propagated backwards. We have run the microscopic  model using real inflow data as left boundary conditions, then we have fitted the parameters in order to catch the real queue as measured on the field}, see Fig.\ \ref{fig:calibrazione_tau}.
    The role of the parameters $\tau_\text{acc}$ and $\tau_\text{dec}$ is to adjust the points/times of start and end of the queue.
   \rev{We have noted a strong sensitivity of the model to those parameters. As a consequence, it is quite difficult to catch the correct speed of the backward propagation of a queue when the inertia comes into play.}
	We summarize the values of the parameters in Table \ref{tab:parametri_mima}.
	\begin{specialtable}[h!]
		\caption{Parameters for the microscopic model.}
		\label{tab:parametri_mima}
		\begin{tabular}{lcl}\toprule
			$\delta$ & $50\times 10^{-3}$ & km \\ \hline
			$\Delta_\text{close}$ & $25\times 10^{-3}$ & km \\ \hline
			$\Delta_\text{far}$ & $50\times 10^{-3}$ & km \\ \hline
			$V_\pes^\text{max}$ & $90$ & km/h \\ \hline
			$\tau_\text{dec}$ & $2\times 10^{-4}$ & h \\ \hline
			$\tau_\text{acc}$ & $1.4\times 10^{-2}$ & h 
			\\\bottomrule 
		\end{tabular}
	\end{specialtable}
	%

	%
	
	\section{Numerical results}\label{sec:numericaltests}
	In this section we present the numerical results obtained with the models \eqref{eq:model} and \eqref{multiscalemodel}.

	\subsection{Macroscopic model}\label{sec:mamatests}
	Here we present three tests which highlight how the macroscopic model reproduces some interesting phenomena arising from the coupled dynamics of cars and trucks. In particular we focus on the creeping phenomenon, the shared occupancy and the stop \& go waves.

		\subsubsection{Test 1A: creeping}\label{sec:T1mama}
\begin{figure}[h!]
		\subfig{a}\includegraphics[width=.45\linewidth]{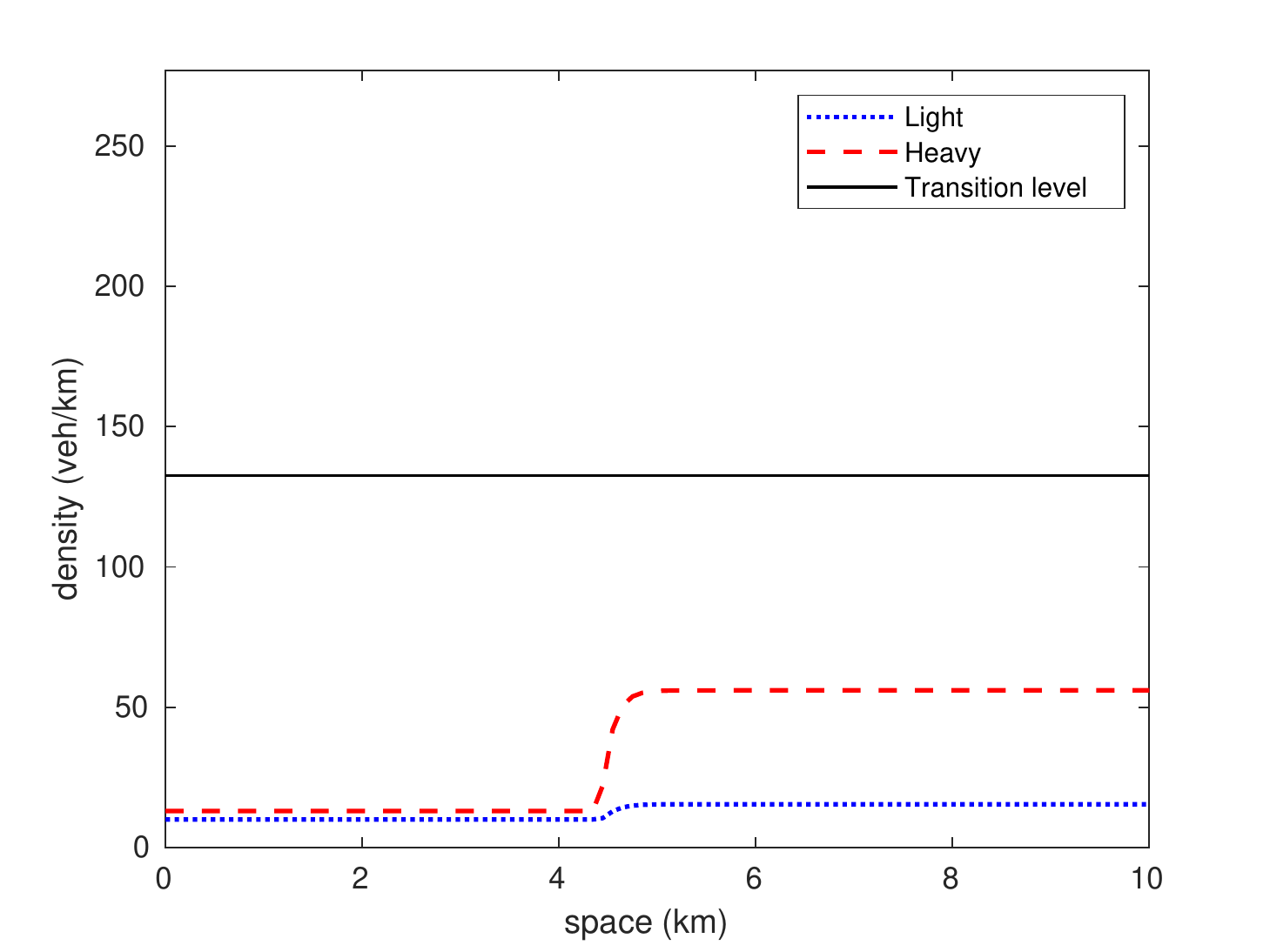}
		\subfig{b}\includegraphics[width=.45\linewidth]{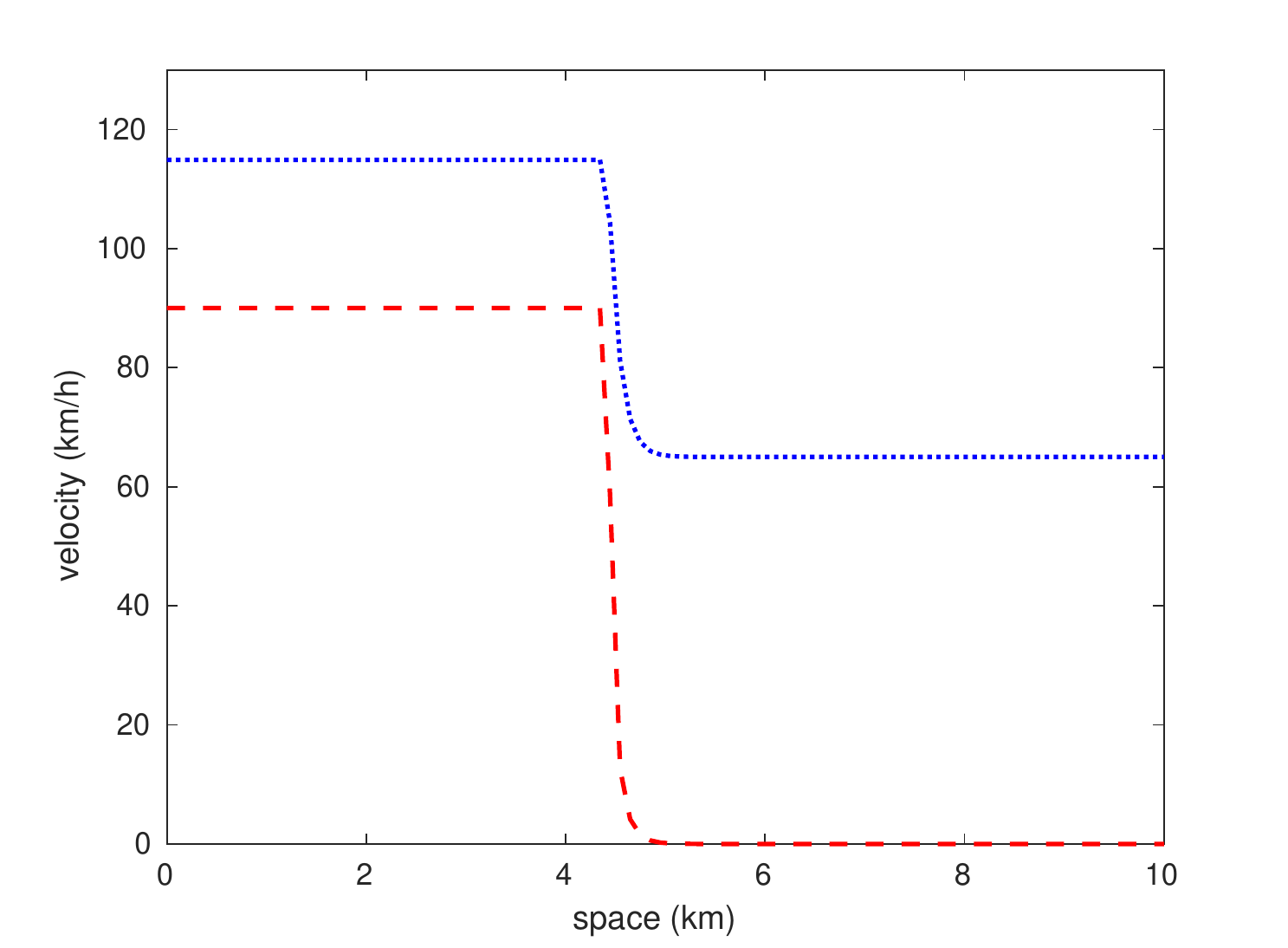}		\subfig{c}\includegraphics[width=.48\linewidth]{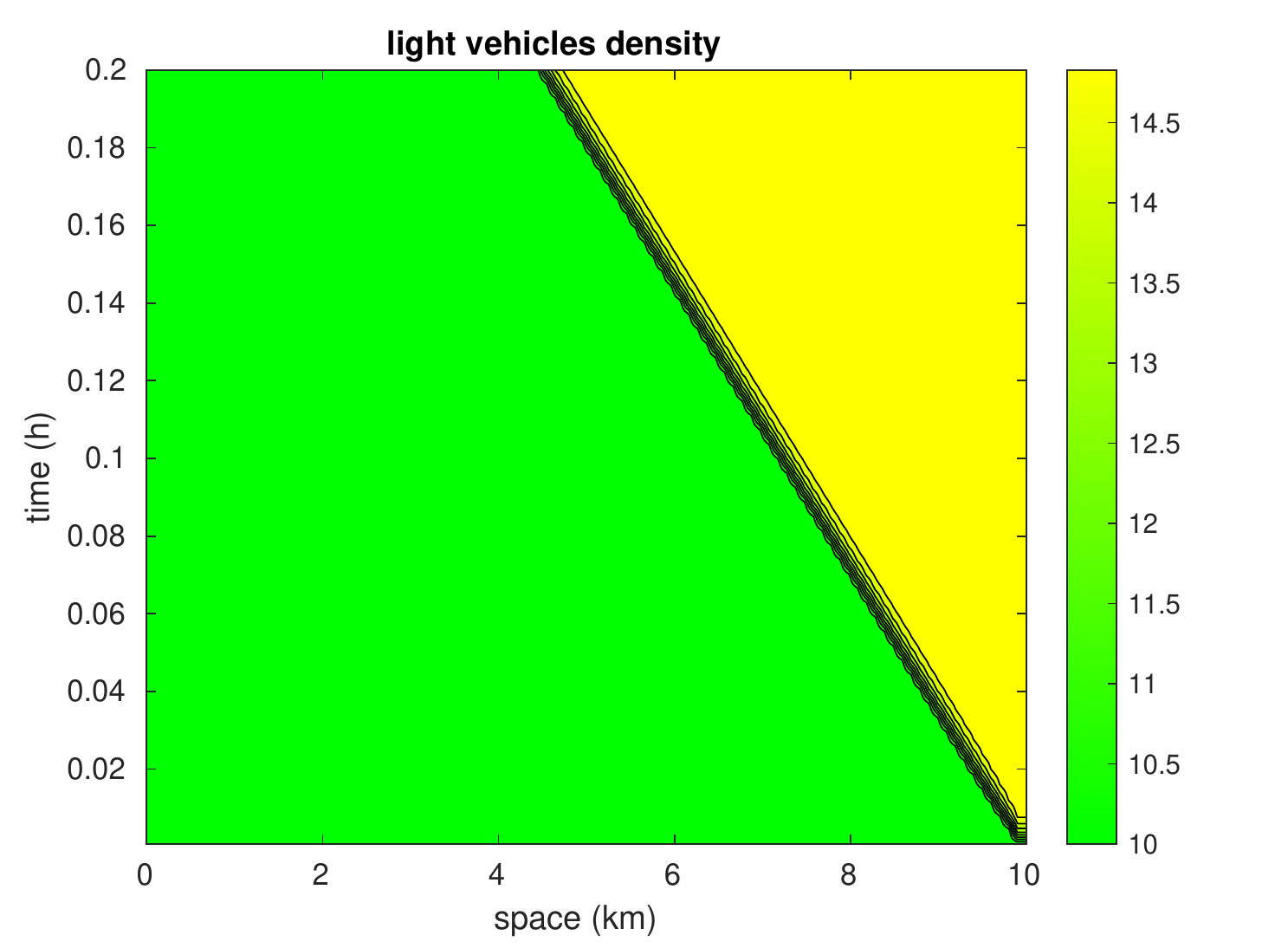}
		\subfig{d}\includegraphics[width=.48\linewidth]{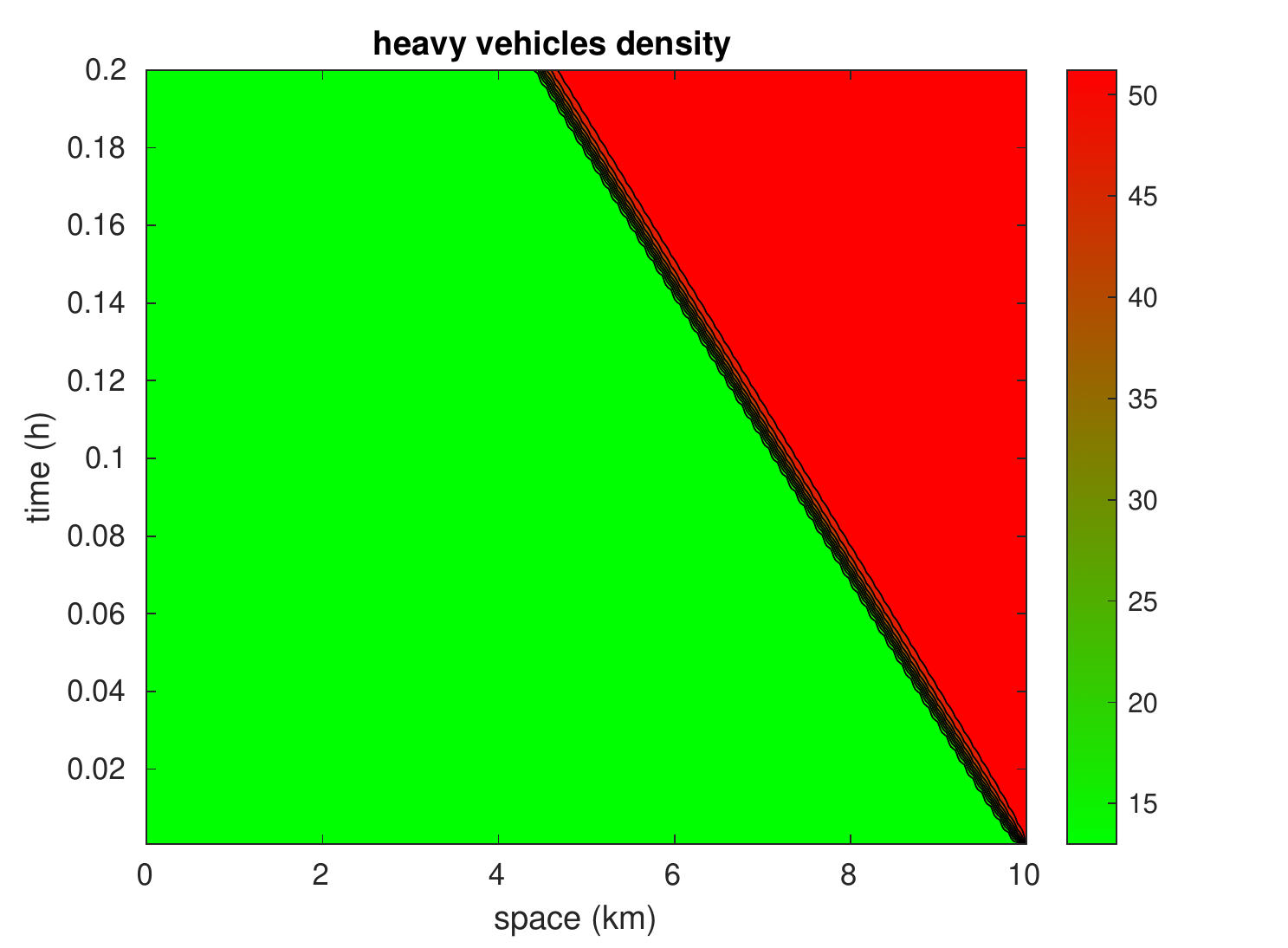}
		\caption{Test 1A: \subfig{a} Density and \subfig{b} velocity of light and heavy vehicles as a function of space at final time. \subfig{c} Density of light and \subfig{d} heavy vehicles in the space-time.}
		\label{fig:T1maya}
		\end{figure}

		In this simple test we observe the creeping phenomenon, see Fig.\ \ref{fig:T1maya}.
		The simulation starts with a constant density $(\rhoL,\rhoP)=(10,13)$ veh/km all along the road. At the end of the road (right boundary), trucks are stopped by fixing their density at its maximum value $\mumax=56$:  
		a queue of trucks is propagating backward from the end of the road while a constant flux of cars approaches the beginning of the queue.
		Once cars reach the trucks' queue they have to slow down but do not stop completely. More precisely, the cars' velocity drops to 65 km/h. 
		Note that cars density remains under the transition level and then the dynamics is in the partial-coupling phase all the time. Moreover, cars are always in the freeflow regime and then move at maximal speed, but the maximal speed changes as a function of the trucks' density.


		\subsubsection{Test 2A: cars congestion affects truck dynamics}\label{sec:T2mama}
\begin{figure}[h!]
	\subfig{a}\includegraphics[width=.45\linewidth]{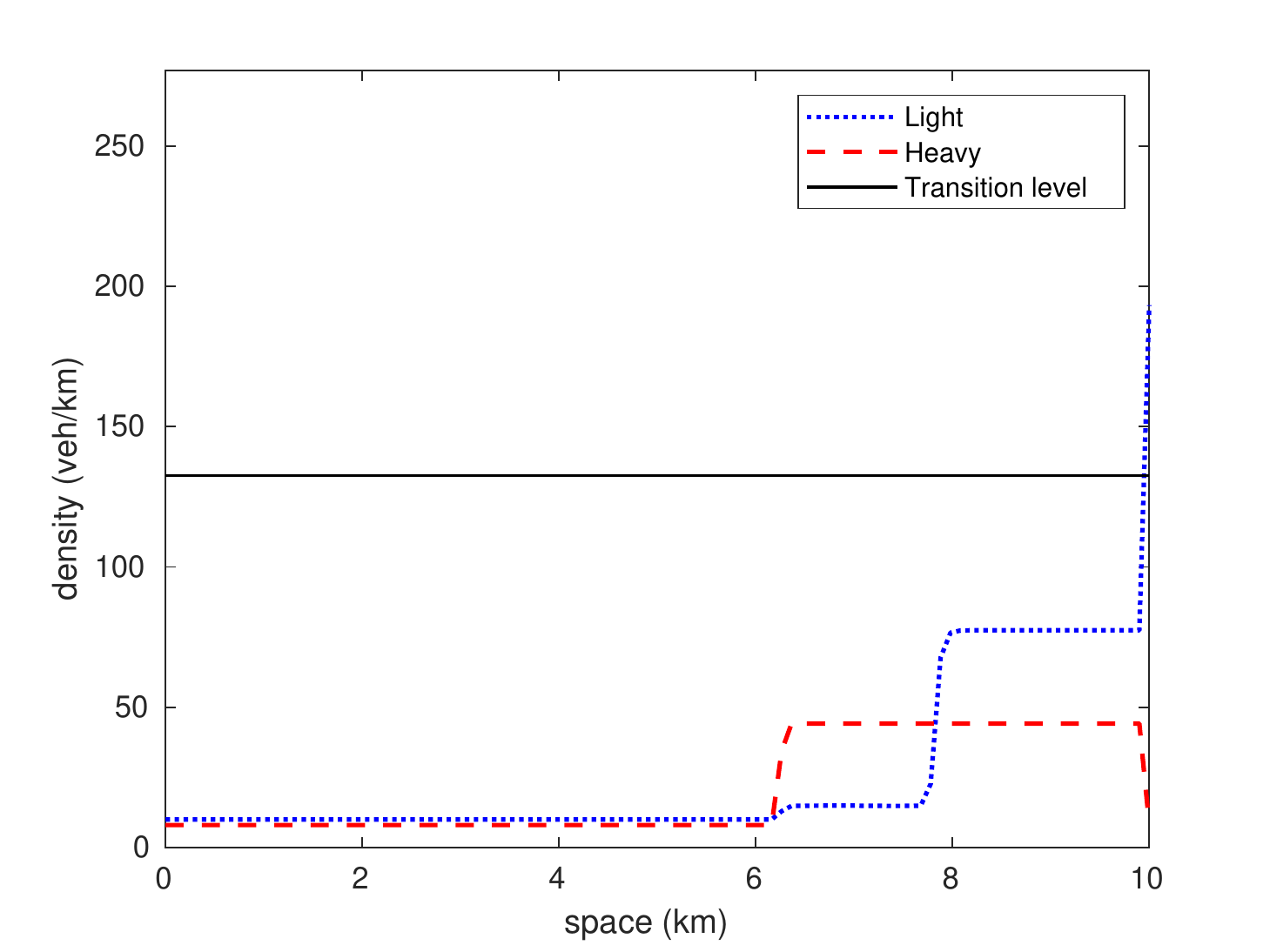}
	\subfig{b}\includegraphics[width=.45\linewidth]{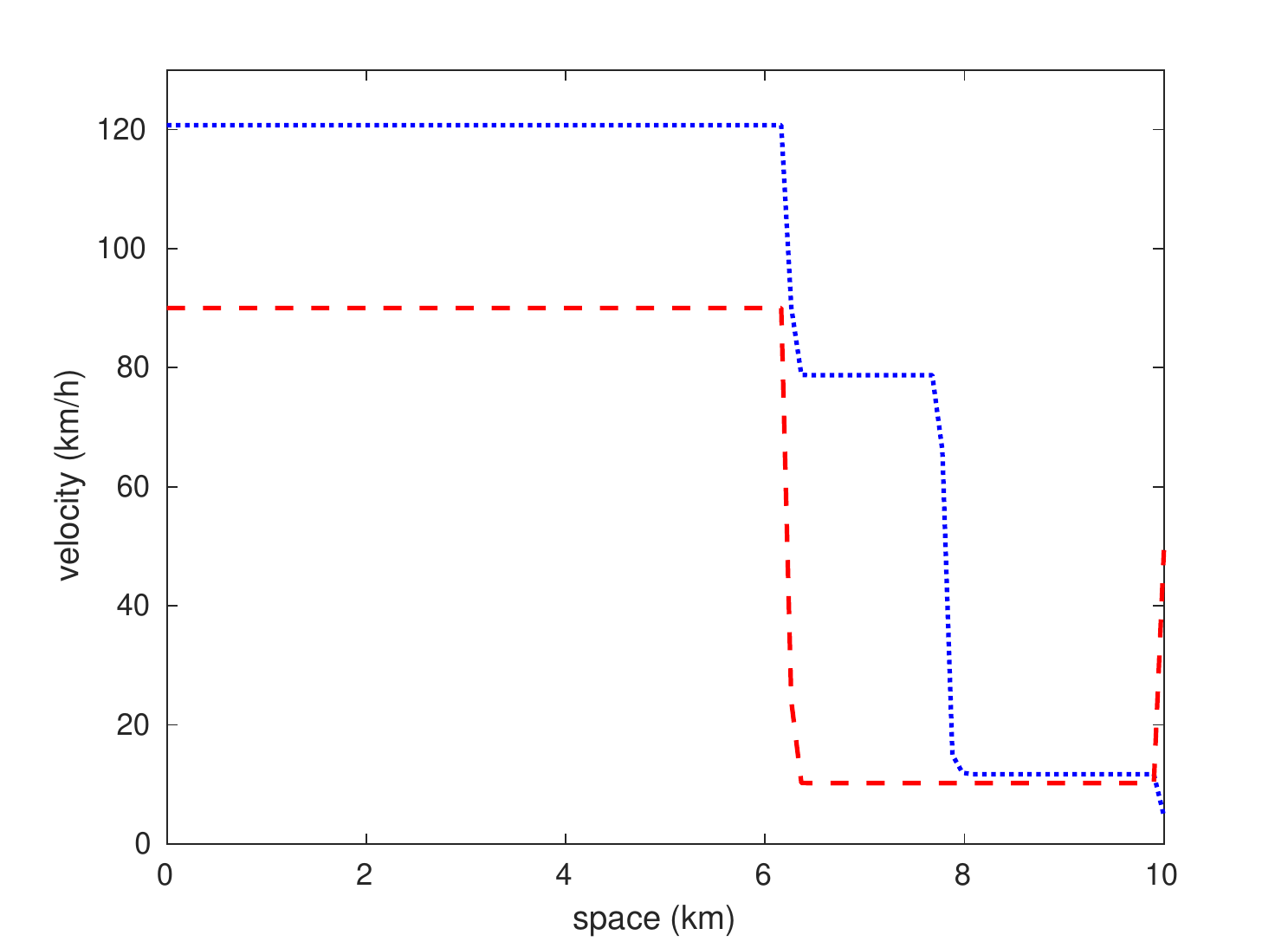}
	\subfig{c}\includegraphics[width=.48\linewidth]{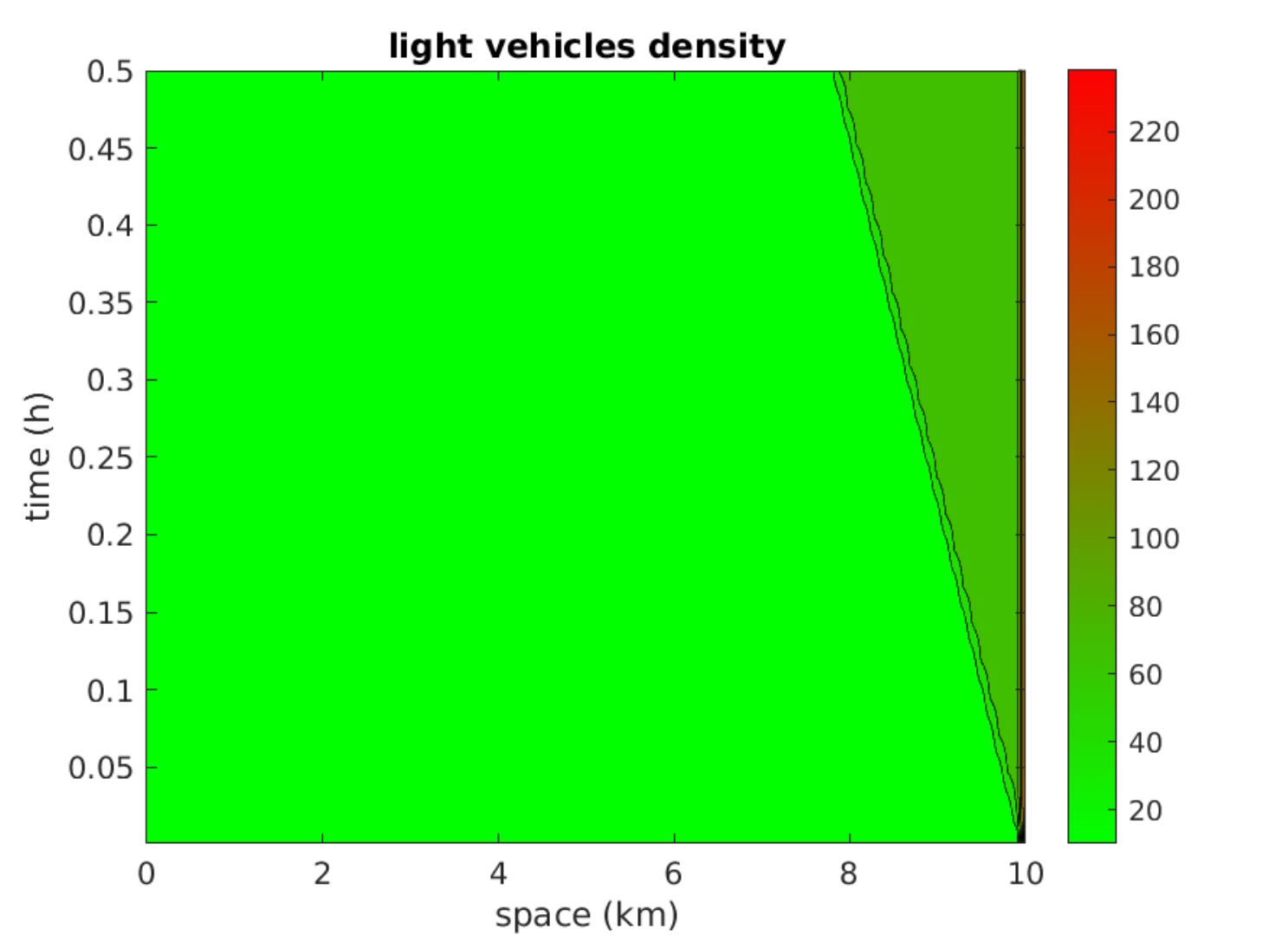}
	\subfig{d}\includegraphics[width=.48\linewidth]{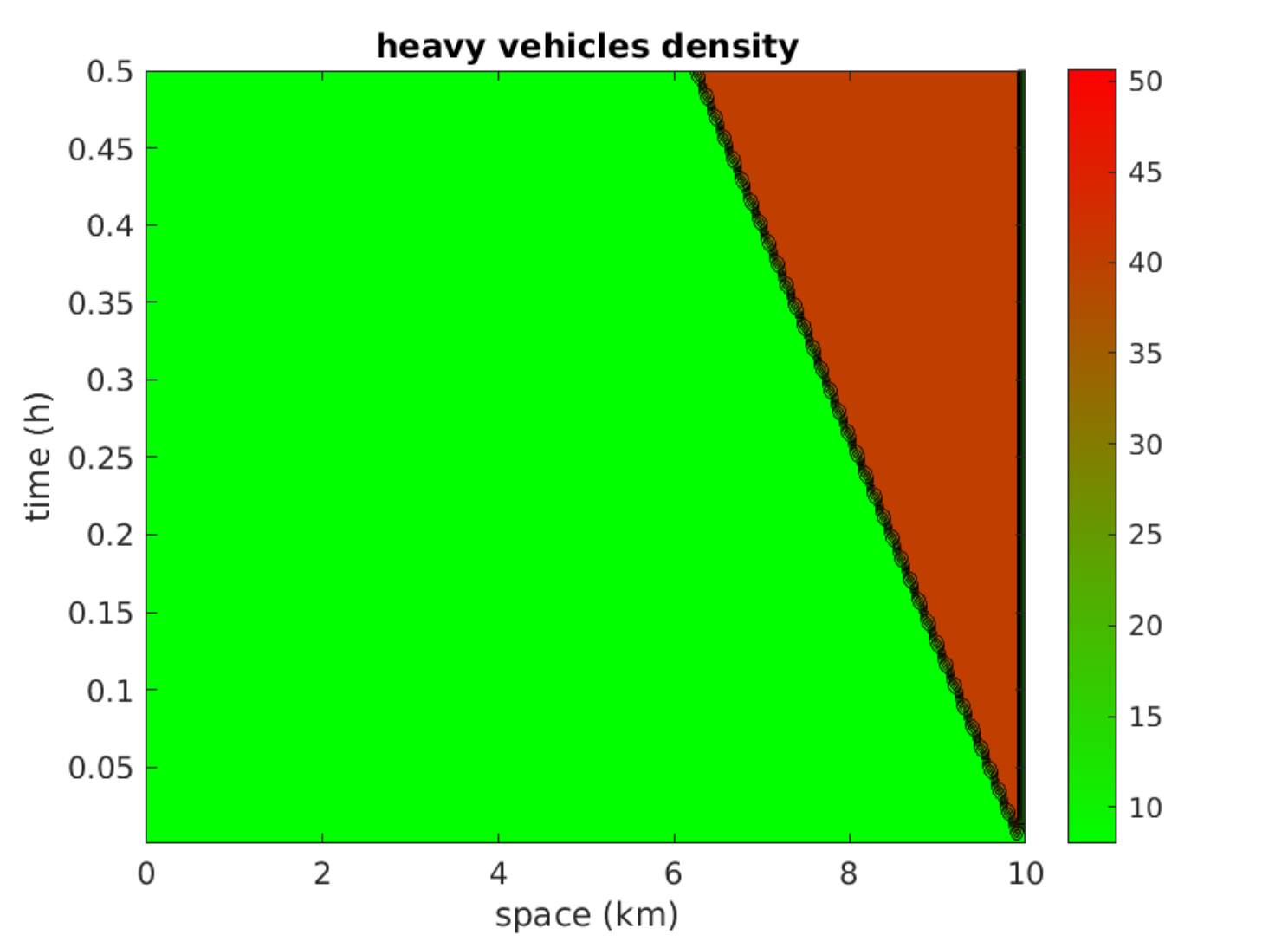}
	\caption{Test 2A: \subfig{a} Density and \subfig{b} velocity of light and heavy vehicles as a function of space at final time. \subfig{c} Density of light and \subfig{d} heavy vehicles in the space-time.}
	\label{fig:T2maya}	
\end{figure}

		In this test we observe the effect of a congestion of cars, see Fig.\ \ref{fig:T2maya}. 
		The simulation starts with a constant density $(\rhoL,\rhoP)=(10,8)$ veh/km all along the road. At the end of the road (right boundary), the density of cars is fixed to 186 veh/km to create the slowdown.
		Car density is larger than the transition level, then cars have to invade the slow lane. 
		Trucks facing cars congestion slow down but do not just occupy the space left to them by cars; rather, they conquer some extra space, thus decreasing the car density. As a result, both cars and trucks proceed slowly without stopping, and the initial cars congestion propagates backward with a density lower than the transition level. 

\subsubsection{Test 3A: stop \& go wave}\label{sec:T3mama}
\begin{figure}[h!]
	\subfig{a}\includegraphics[width=.45\linewidth]{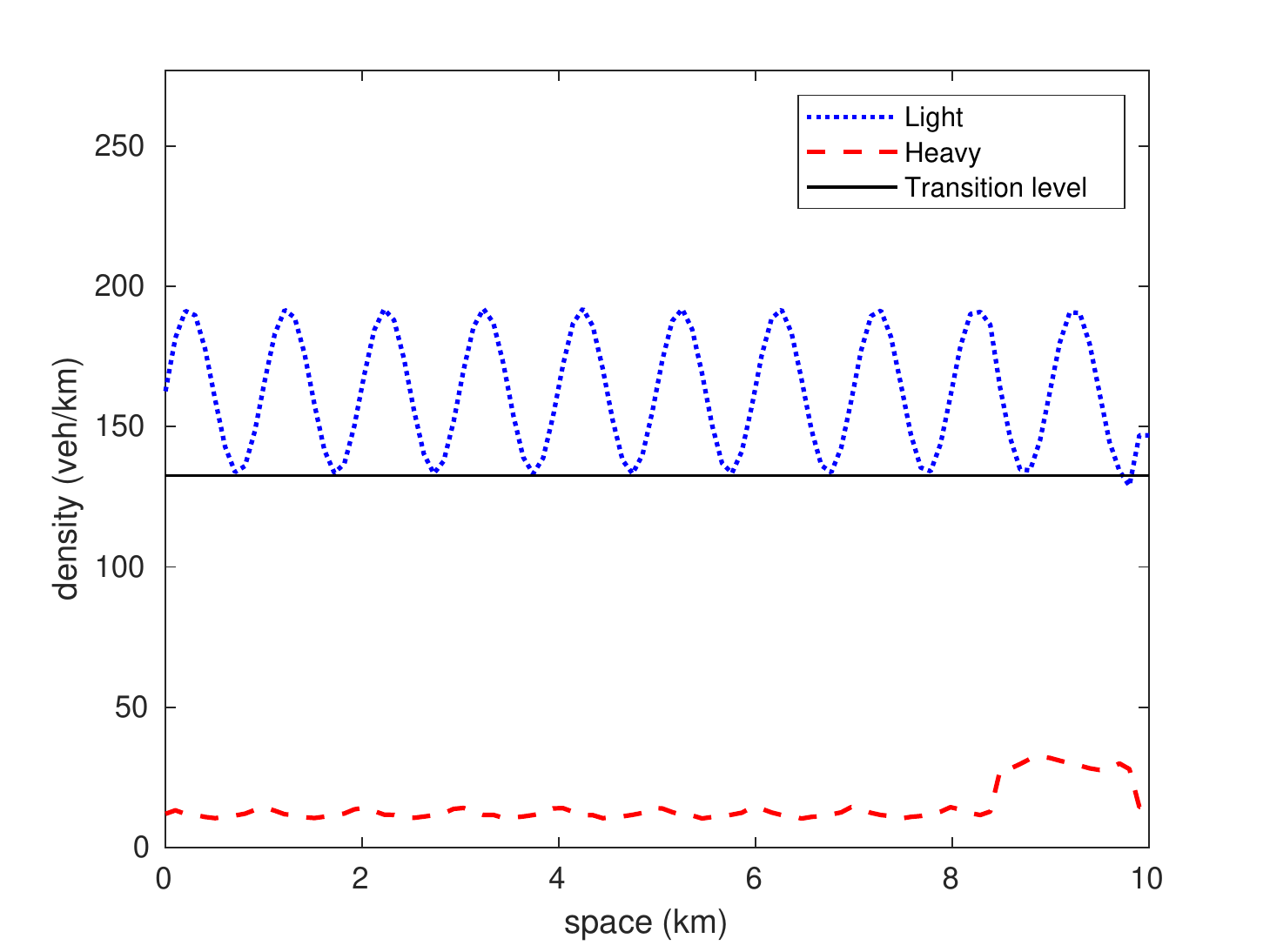}
	\subfig{b}\includegraphics[width=.45\linewidth]{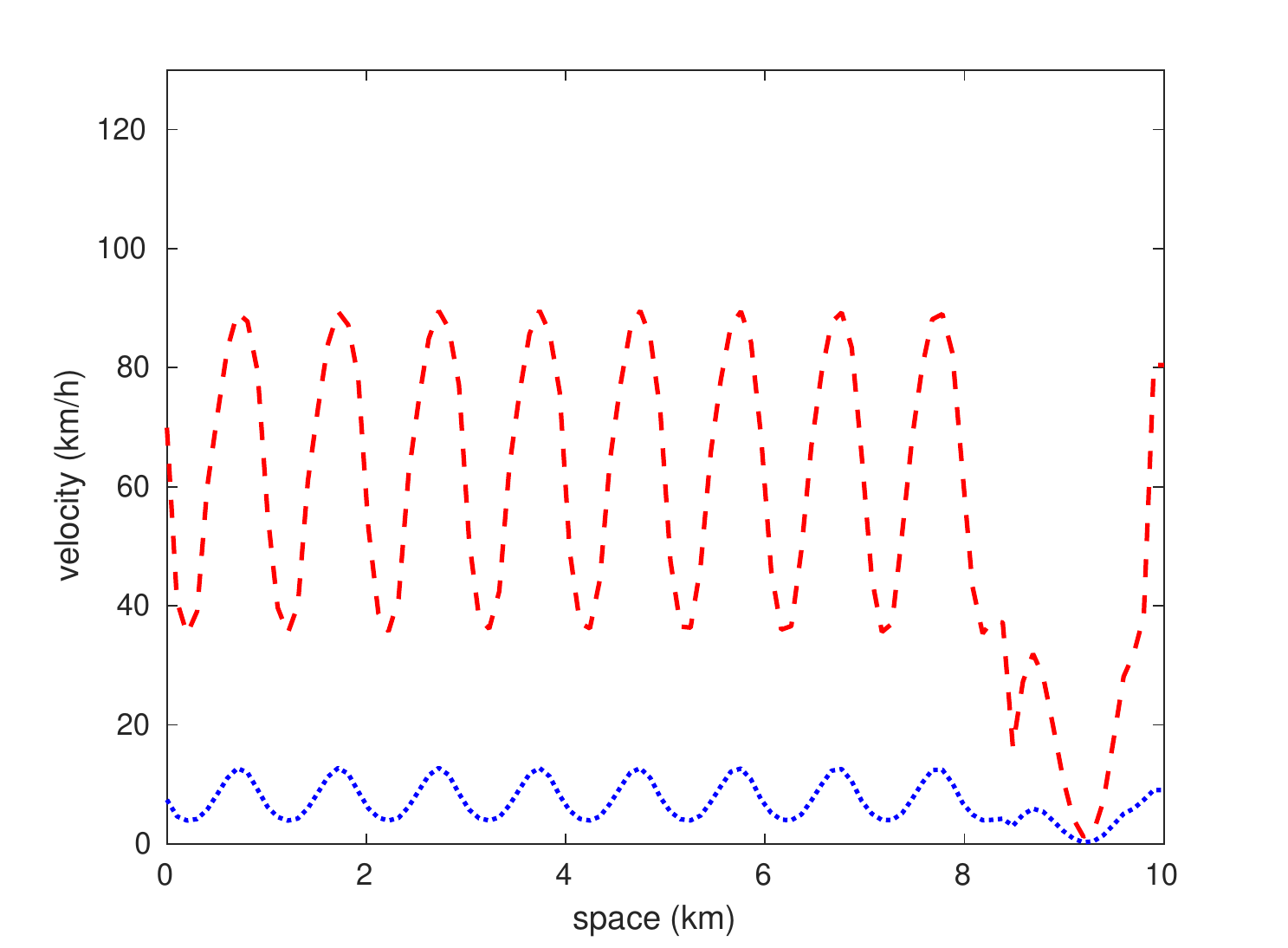}
	\subfig{c}\includegraphics[width=.48\linewidth]{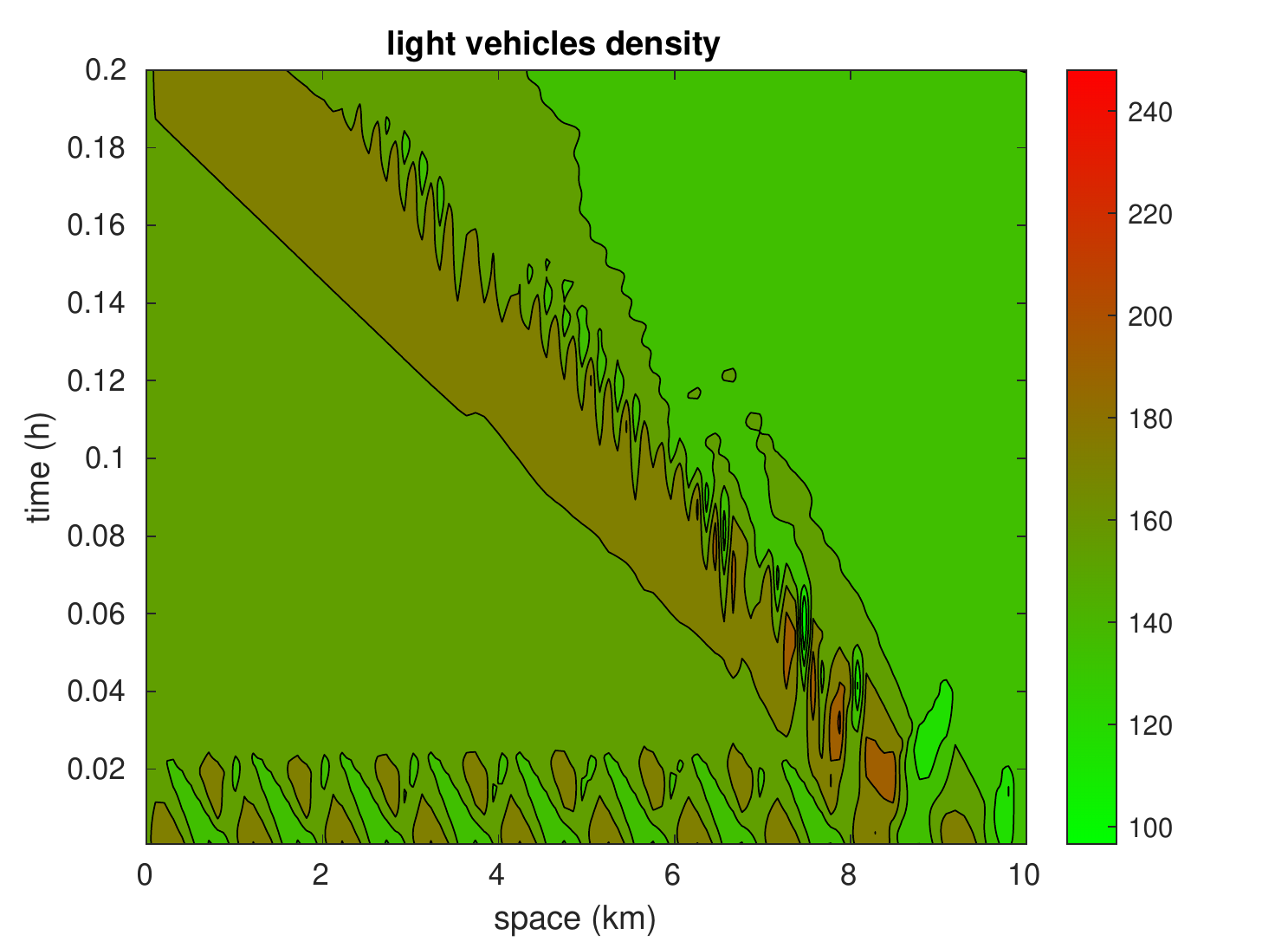}
	\subfig{d}\includegraphics[width=.48\linewidth]{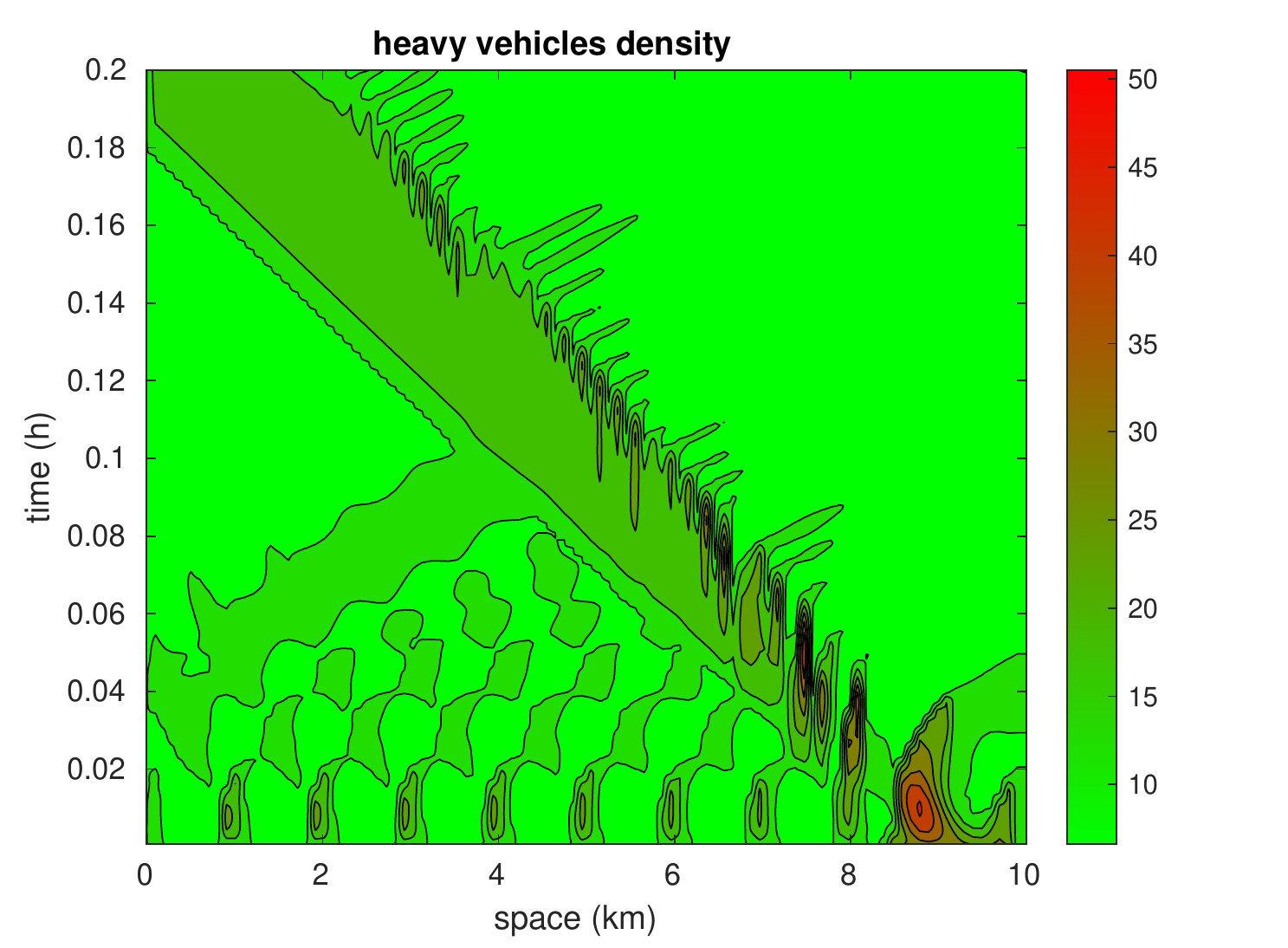}
	\caption{Test 3A: \subfig{a} Density and \subfig{b} velocity of light and heavy vehicles as a function of space at $t=\Delta t$ (i.e.\ just after the initial time). 
		\subfig{c} Density of light and \subfig{d} heavy vehicles in the space-time. It is perfectly visible the evolution of the initial perturbation of trucks density starting at km 9, which creates, in turn, a perturbation of cars density.
	}
	\label{fig:T3maya}	
\end{figure}

In this test we study the evolution of a small perturbation in the trucks density, see Fig.\ \ref{fig:T3maya}. 
At initial time the trucks density is constant and equal to 12 veh/km except for a small perturbation at the end of the road where the density is equal to 30 veh/km. Cars density instead oscillates just above the transition level. 
It is plain that a single-class LWR model for trucks only would flatten the perturbation in short time. Conversely, in this case the coupling with cars dynamics makes the perturbation propagate backward without vanishing. This second-order-like effect is obtained thanks to the fact that the fundamental diagram of trucks is continuously modified by the oscillating cars density.
%
%
%
%
	\subsection{Multi-scale model}\label{sec:mimatests}
	Here we replicate, with the multi-scale model, the first two scenarios already investigated in Section \ref{sec:mamatests}. The third scenario was already considered in Fig.\ \ref{fig:calibrazione_tau} where the second-order microscopic model is able to reproduce stop \& go waves alone, without  the need to couple cars dynamics. 
	Finally we consider the case of a merge.
	
	\subsubsection{Test 1B: creeping effect}\label{sec:T1mima}
\begin{figure}[h!]
		\subfig{a}\includegraphics[width=.45\linewidth]{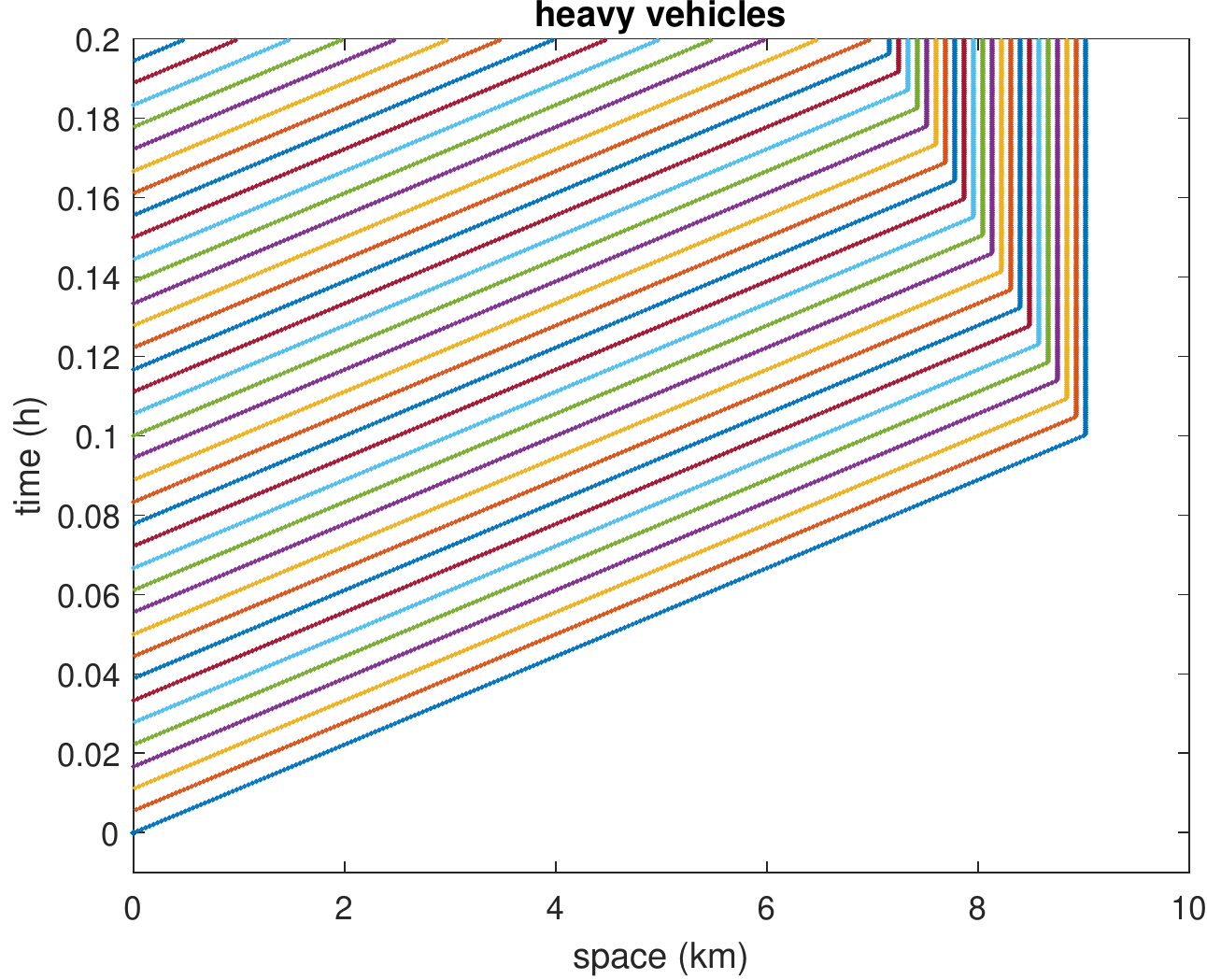}
		\subfig{b}\includegraphics[width=.482\linewidth]{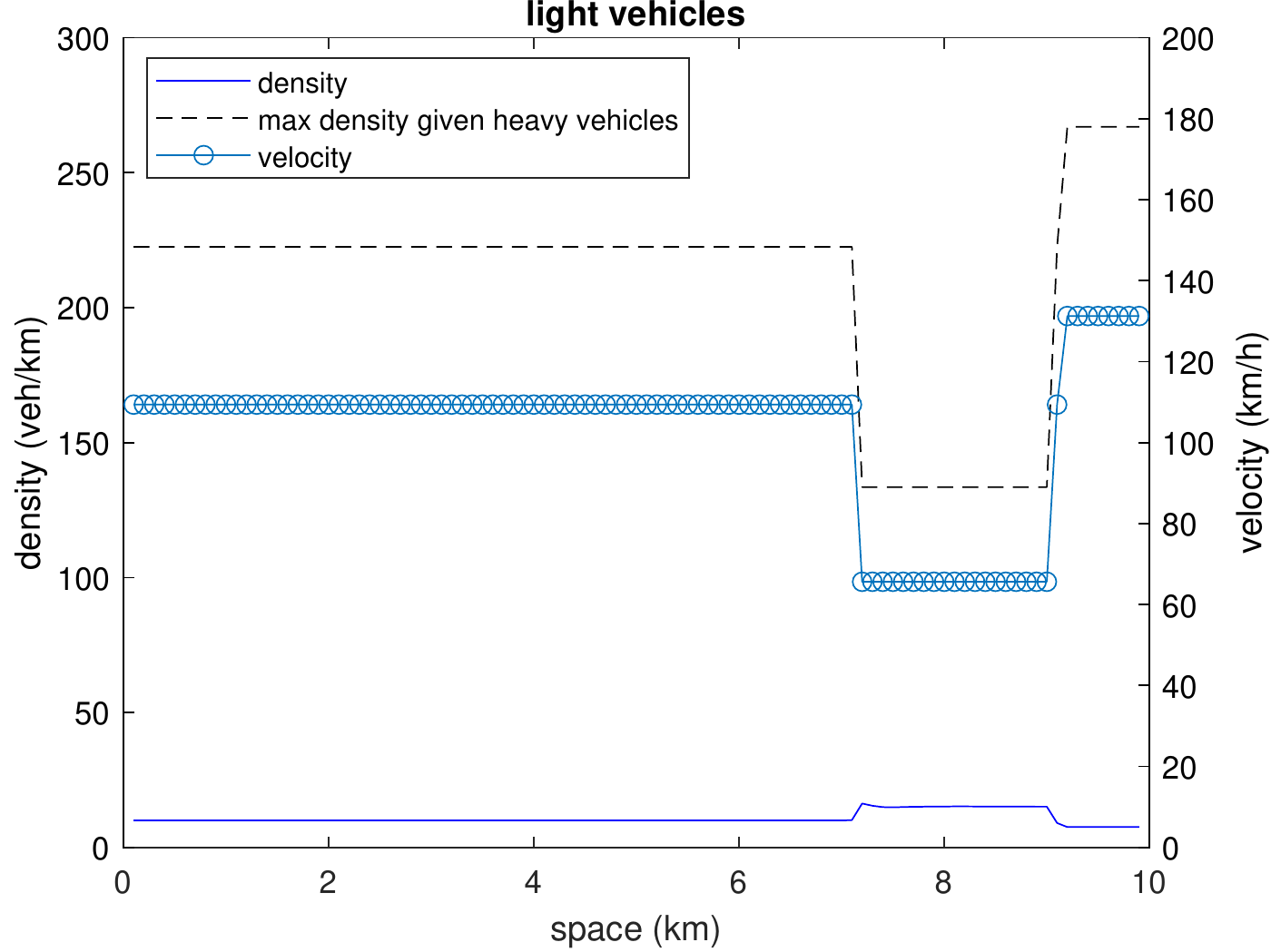}		
		\caption{Test 1B: 
			\subfig{a} Trajectories of trucks in the space-time (for visualization purposes not all trucks are actually plotted). When the first truck stops, a queue is formed behind. 
			\subfig{b} Cars density, cars velocity, and cars maximal density given the number of trucks at final time. Creeping is visible between km 7 and km 9.
		}
		\label{fig:T1emi}
	\end{figure}

	As in Test 1A in Section \ref{sec:T1mama}, here one truck stops completely and creates a long queue of trucks behind, which saturates the slow lane. When cars reach the trucks queue, they all move to the fast lane keeping moving at (the new, reduced) maximal velocity of 65 km/h, see Fig.\ \ref{fig:T1emi}.
	
	
	\subsubsection{Test 2B: cars congestion affects truck dynamics}\label{sec:T2mima}
\begin{figure}[h!]
		\subfig{a}\includegraphics[width=.45\linewidth]{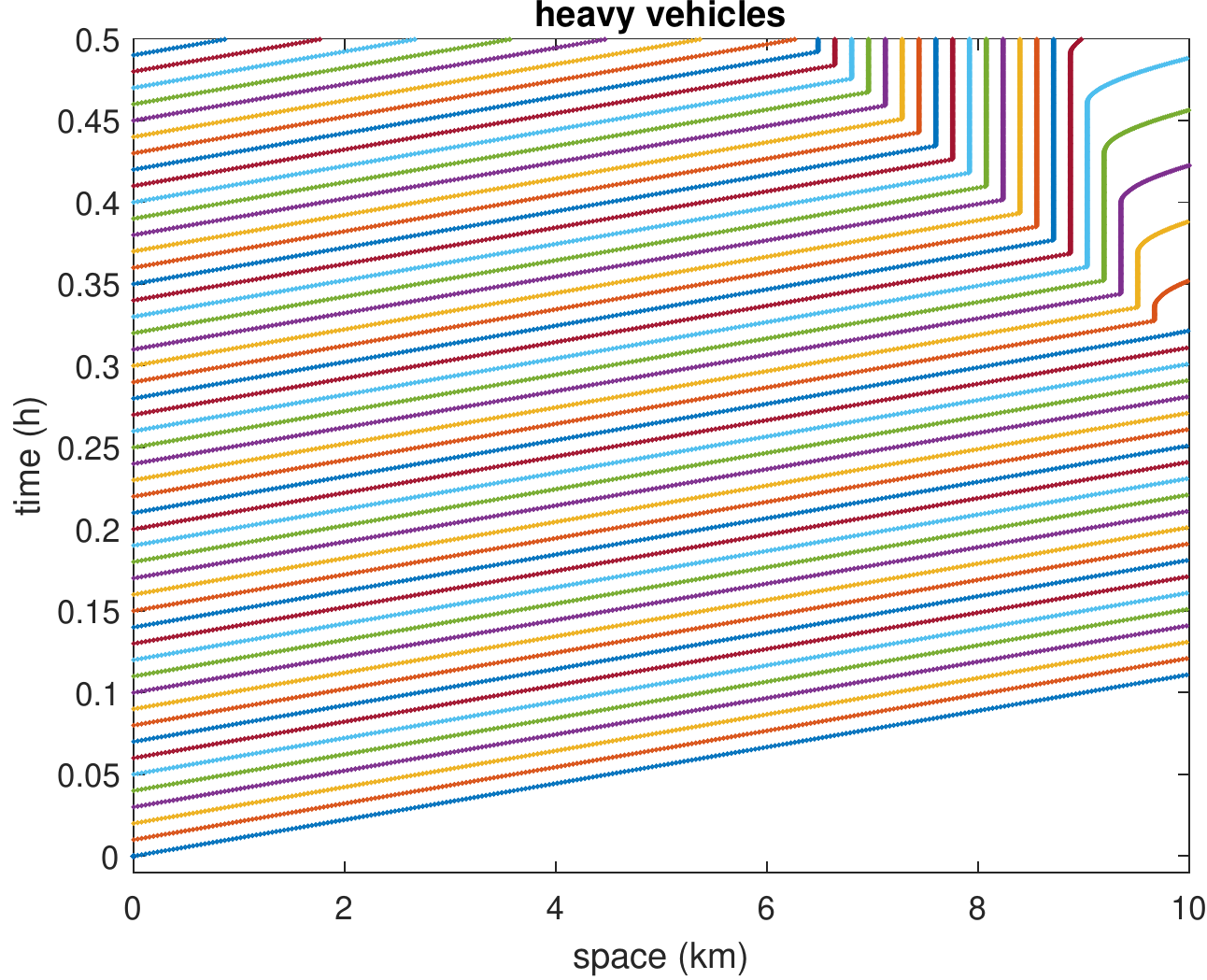}
		\subfig{b}\includegraphics[width=.482\linewidth]{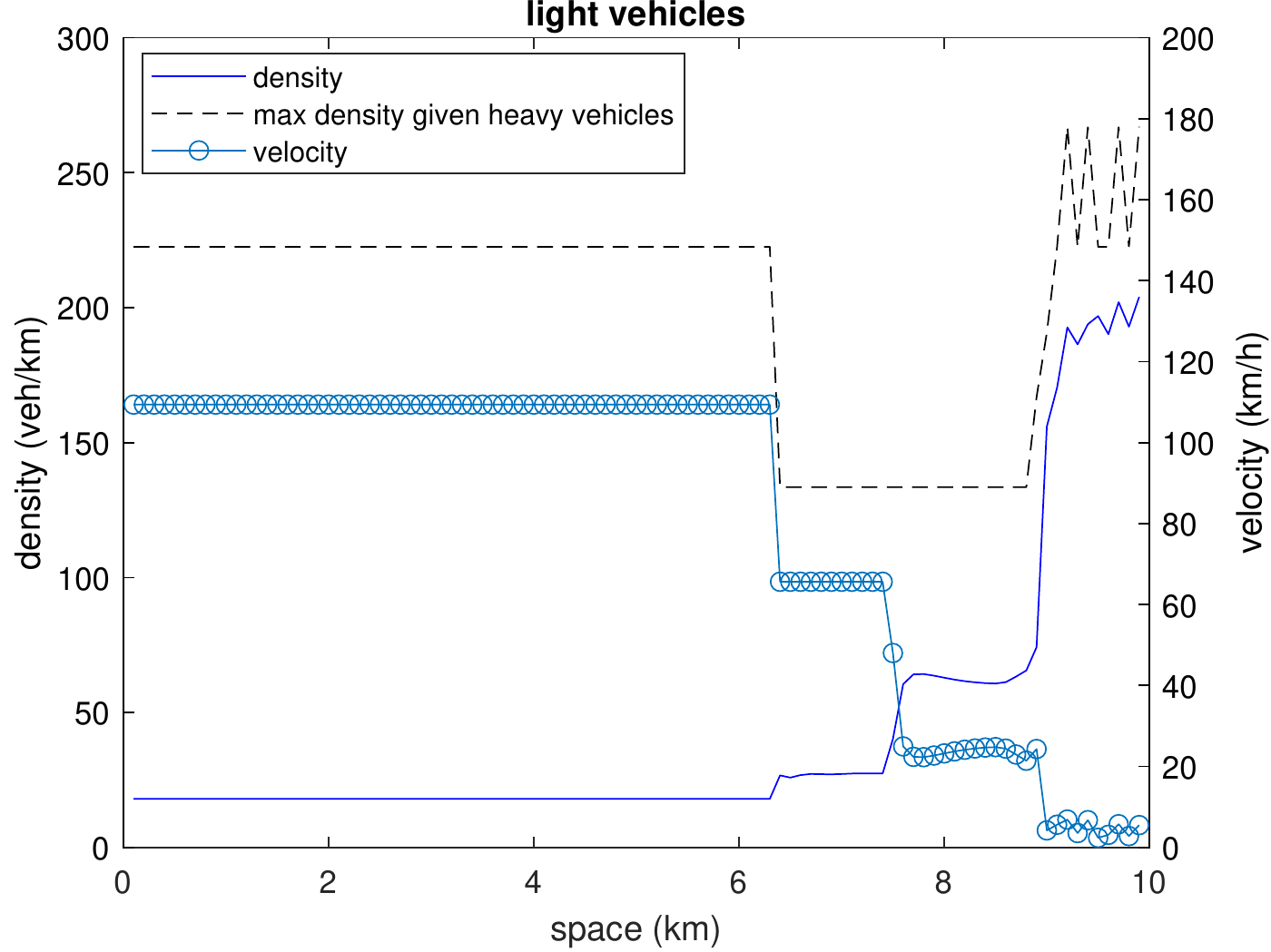}		
		\caption{Test 2B: 
			\subfig{a} Trajectories of trucks in the space-time (for visualization purposes not all trucks are actually plotted). They stop for a while then accelerate. 
			\subfig{b} Cars density, cars velocity, and cars maximal density given the number of trucks at final time. 
		}
		\label{fig:T2emi}
	\end{figure}

	As in Test 2A in Section \ref{sec:T2mama}, a congestion of cars at the end of the road slows down trucks, see Fig.\ \ref{fig:T2emi}. 
	Results are similar to those obtained by the macroscopic but here trucks stop completely, forming a queue. 

	
	\subsubsection{\rev{Test 3B: merge}}\label{sec:T3mima}
\begin{figure}[h!]
		\subfig{a}\includegraphics[width=.295\linewidth]{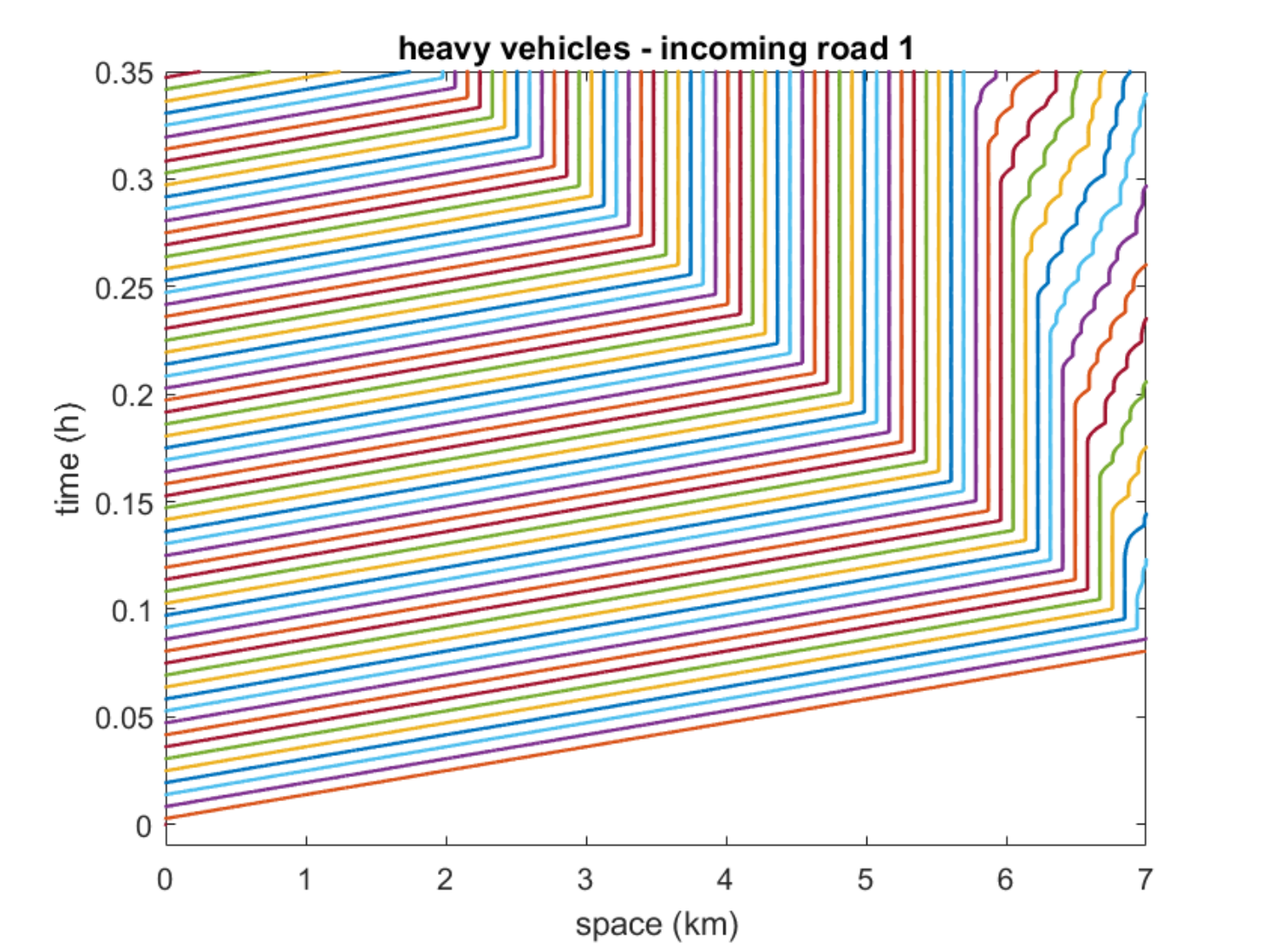}
		\subfig{b}\includegraphics[width=.295\linewidth]{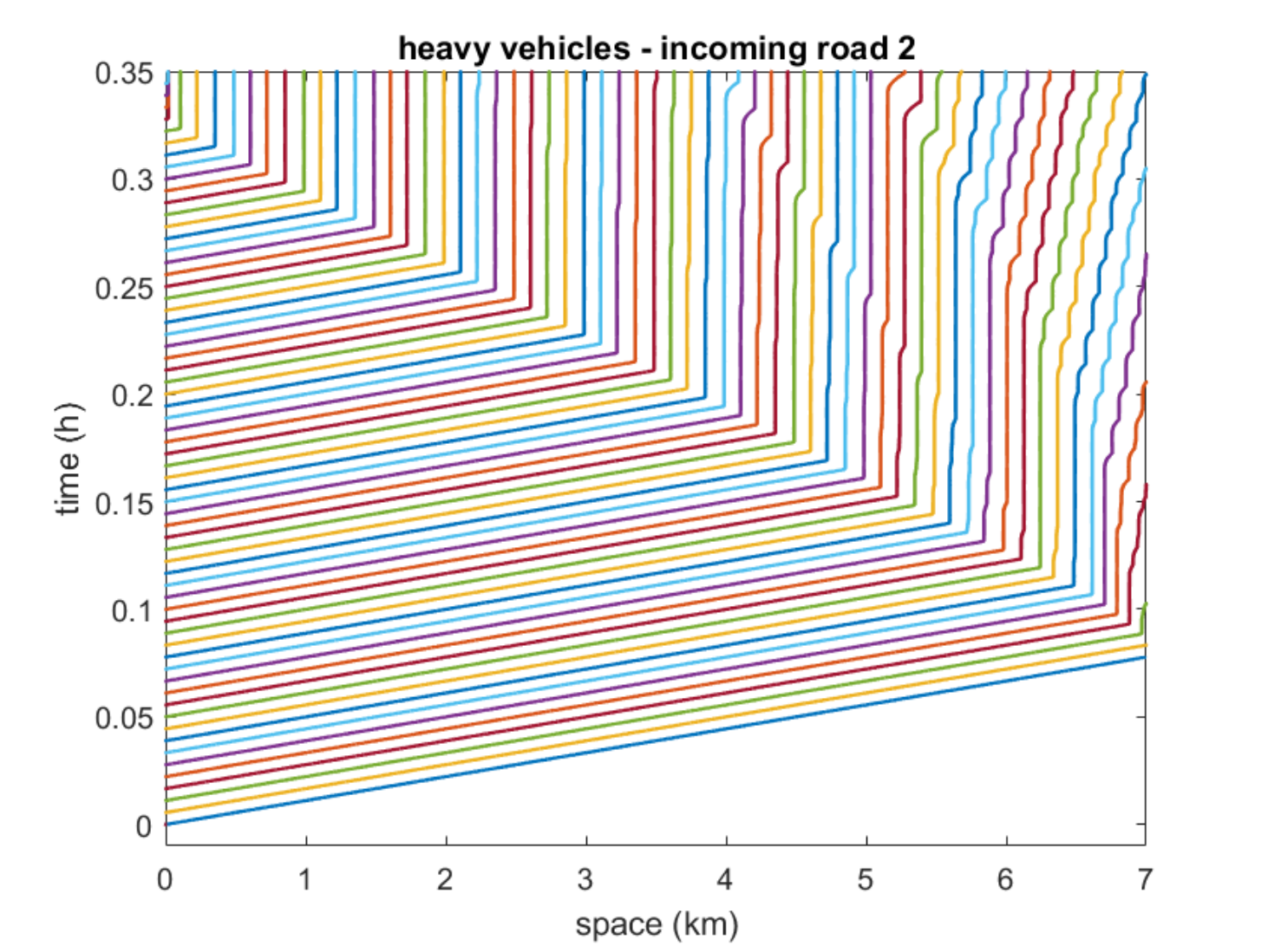}
		\subfig{c}\includegraphics[width=.295\linewidth]{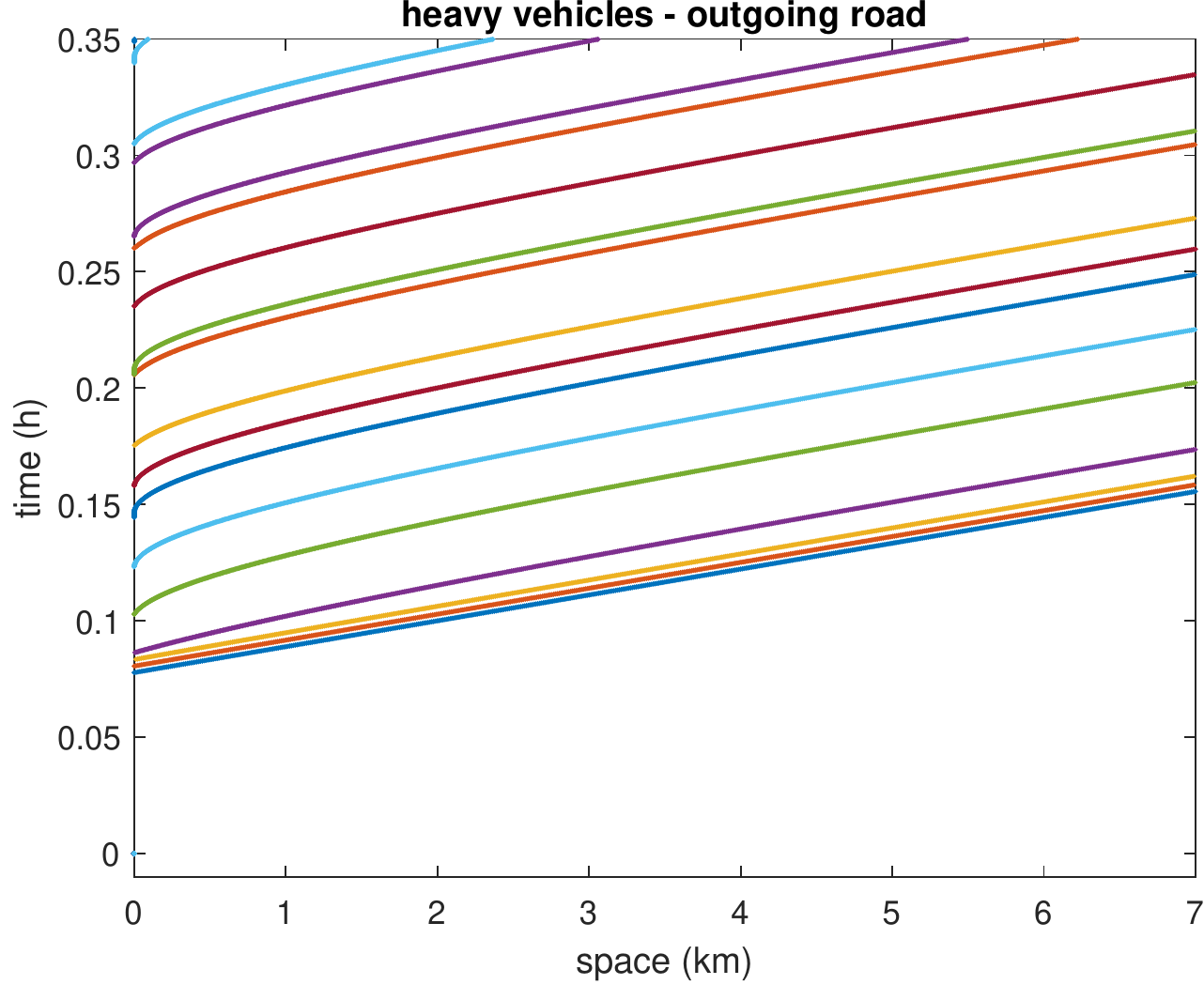}\\
		\subfig{d}\includegraphics[width=.47\linewidth]{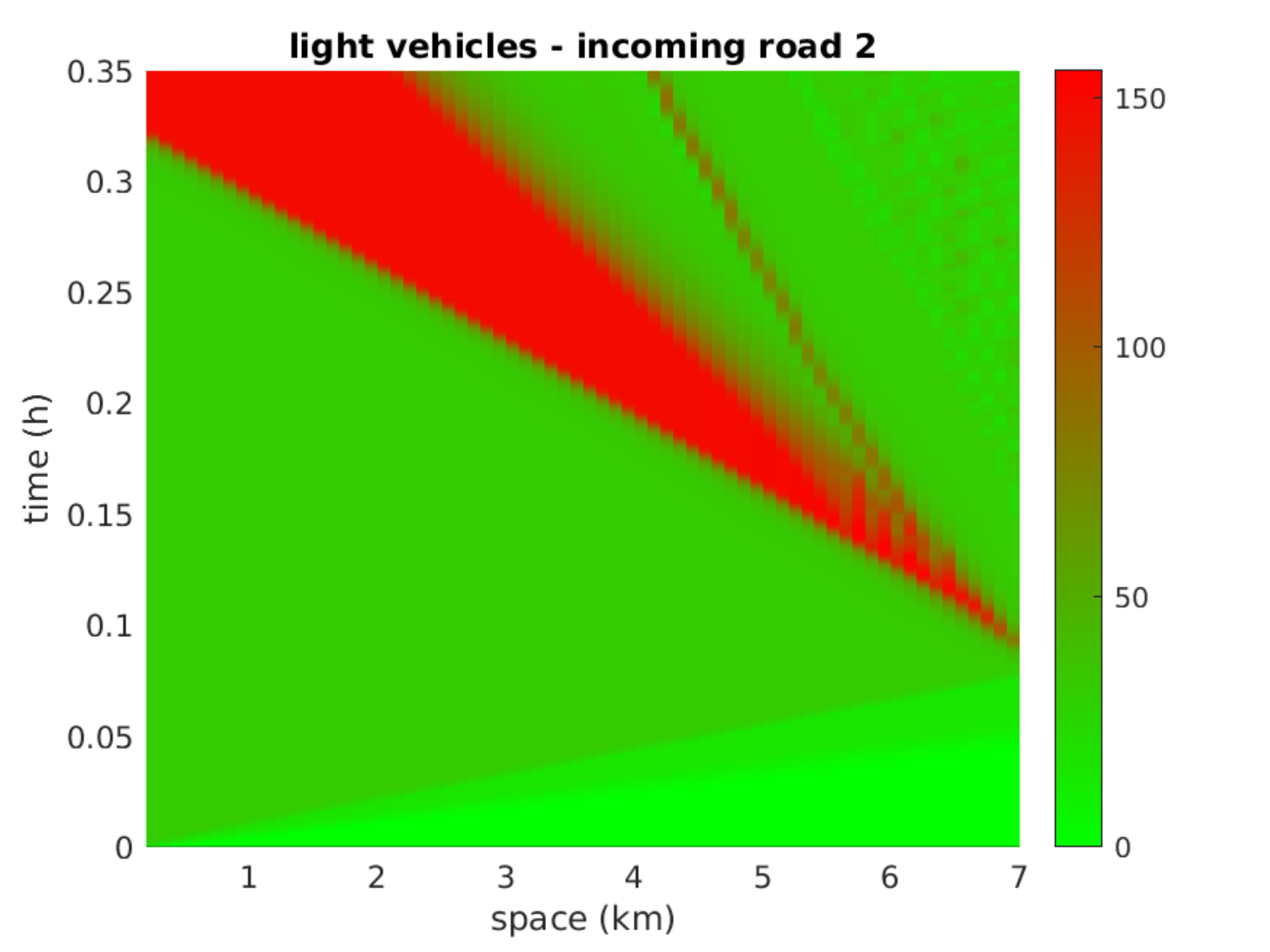}	 		\subfig{e}\includegraphics[width=.47\linewidth]{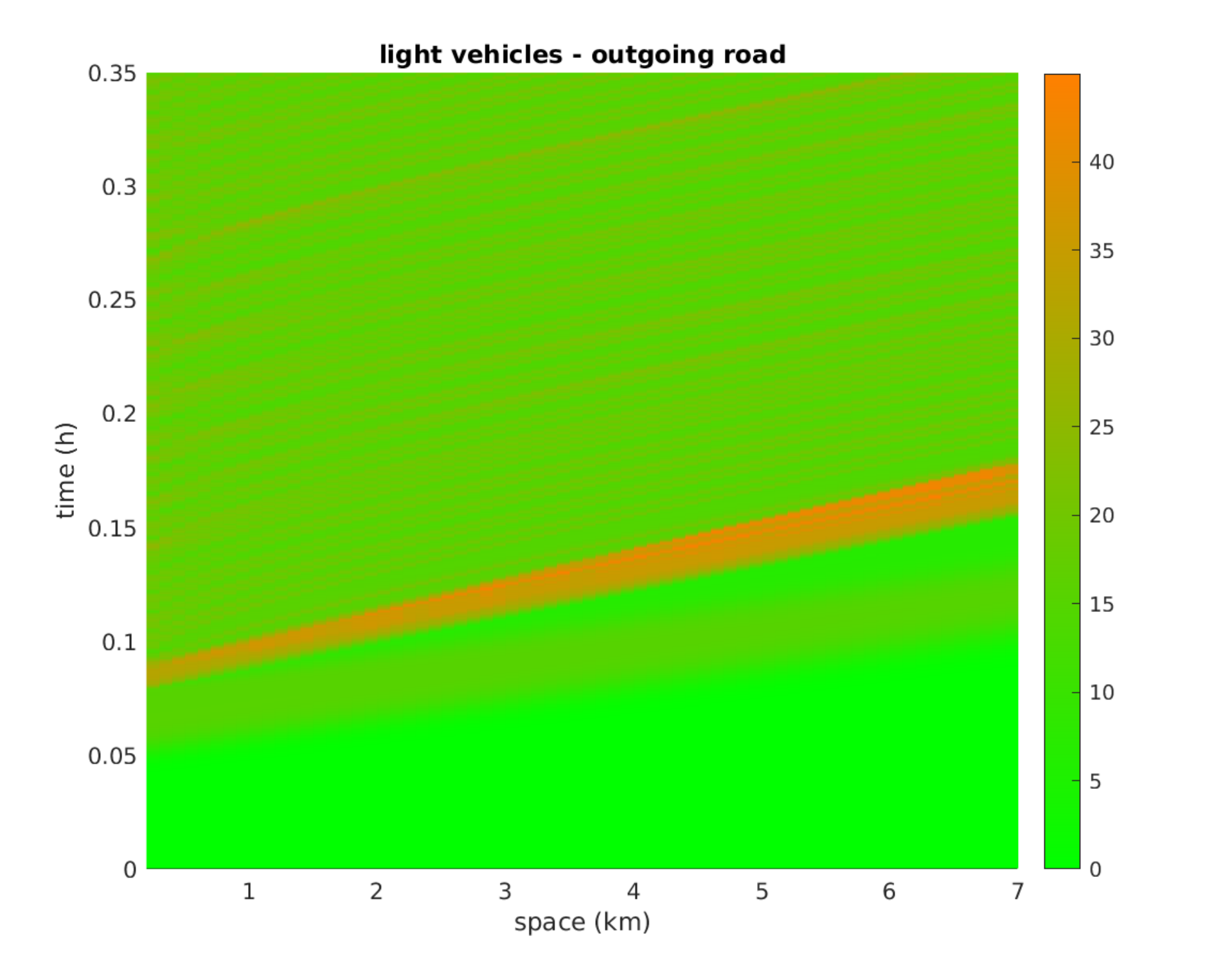}			
		\caption{Test 3B: 
			(\textbf{a},\textbf{b},\textbf{c}) Trajectories of trucks in the space-time on the first incoming road, second incoming road and outgoing road, respectively (for visualization purposes not all trucks are actually plotted).  
			(\textbf{d},\textbf{e}) Cars density on  second incoming road and outgoing road, respectively.
		}
		\label{fig:T3emi}
	\end{figure}

	In this test we consider a merge (2 incoming roads and 1 outgoing road). At time $t=0$ the three roads are empty. A constant inflow of trucks (1 every 4 s) comes from the left boundary of both incoming roads, while a constant density of cars ($\rhoL=32$)  is imposed as Dirichlet left boundary condition on the second incoming road only. The first incoming road has no cars. 
	When trucks reach the junction and merge, they suddenly break and rapidly form a queue which propagates backward along both incoming roads, see Fig.\ \ref{fig:T3emi}(\textbf{a},\textbf{b}). 
  Queues are not identical due to the presence of cars along the second incoming road. 
  One can note that when the trucks downstream of the queue start moving again, their flux is not maximal: indeed, if the flow were maximum, a queue at the junction would immediately reform as it happened in the first place. 
  This is the well known \emph{capacity drop} phenomenon, ruled by $\tau_\text{acc}$, cf.\ \cite{calvert2018JTTE}. 
  As a consequence, trucks are able to cross the junction without spillback.
  Cars, instead, move at maximal flux until they encounter the trucks queue. The queue acts as a moving bottleneck and drops the road capacity, therefore the car traffic immediately enters the congested state and the density increases. Downstream, the density remains in the free flow state and cars cross the junction without spillback, see Fig.\ \ref{fig:T3emi}(\textbf{d},\textbf{e}).


	\section{Conclusions and future work}\label{sec:conclusions}
	In this paper we have presented two models for two-class traffic flow. 
	Although models are tailored for a specific case study, they are sufficiently general to be useful in other motorways.
	Moreover, both models can be easily generalized to more than two classes of vehicles and a different ratio between the number of lanes used by trucks and  the number of lanes used by cars. 
	
	We have shown that the models are able to reproduce both qualitatively and quantitatively some notable traffic phenomena arising from the interactions of the two classes. 
	Interestingly, the macroscopic model, although purely first-order, is able to reproduce stop \& go waves thanks to the coupling of the two classes. 
	
	After this preliminary analysis, it is possible to sketch some conclusions about advantages and drawbacks of the two models:
	the multi-scale model has a greater potential since the second-order microscopic part makes it more realistic and then suitable for quantitative predictions. Nevertheless, the macroscopic model appears to be simpler and more manageable, thus representing a valid alternative if one wants to avoid tracking all single vehicles, especially for saving computational time.
	
	\rev{In conclusions, we think that both the proposed models represent the best compromise between accuracy and implementability. In fact, decoupling the dynamics of different classes simplifies excessively the problem description and does not allow to get accurate forecast; 
	conversely, moving to second-order macroscopic models or including multi-lane features to the models increases notably the complexity of the code as well as the number of parameters to be tuned. These generalizations would allow in principle to catch easily inertia-based phenomena in all classes of vehicles and to track the density of \emph{each class} of vehicles in \emph{each lane} but, in our opinion, they make that model unfeasible for practical applications.}
	
	\medskip
	
	In the next future we plan to improve the models including the possibility that they are fed by both Lagrangian (GPS-like) and Eulerian data coming from mobile and fixed sensors, respectively, cf.\ \cite{colombo2016M3AS}. 
	Moreover, we plan to estimate in real time the difference between predicted and measured  densities using the machinery developed in \cite{briani2018CMS}, hopefully creating an algorithm for the auto-calibration of the models in real time.

	\authorcontributions{Conceptualization, Maya Briani and Emiliano Cristiani; Data curation, Paolo Ranut; Funding acquisition, Maya Briani, Emiliano Cristiani and Paolo Ranut; Investigation, Maya Briani, Emiliano Cristiani and Paolo Ranut; Methodology, Maya Briani and Emiliano Cristiani; Visualization, Maya Briani and Emiliano Cristiani; Writing – original draft, Maya Briani and Emiliano Cristiani; Writing – review \& editing, Maya Briani, Emiliano Cristiani and Paolo Ranut.}

	\funding{
	This work was partially funded by the company Autovie Venete S.p.A. 
	
	The work was also carried out within the research project ``SMARTOUR: Intelligent Platform for Tourism'' (No.\ SCN\_00166) funded by the Ministry of University and Research with the Regional Development Fund of European Union (PON Research and Competitiveness 2007-2013).
	
	Authors also acknowledge the Italian Minister of Instruction, University and Research to support this research with funds coming from the project entitled \textit{Innovative numerical methods for evolutionary partial differential equations and applications} (PRIN Project 2017, No.\ 2017KKJP4X).
	
	M.B.\ and E.C.\ are members of the INdAM Research group GNCS.
}

	\acknowledgments{
		The authors want to thank all the Autovie Venete staff as well as Gabriella Bretti, Matteo Piu, Elisa Iacomini, Caterina Balzotti, and Elia Onofri for valuable help.
	}

	\conflictsofinterest{The authors declare no conflict of interest.}

\end{paracol}

\reftitle{References}
\externalbibliography{yes}
\bibliography{biblio}

\begin{thebibliography}{999}

\bibitem[Lighthill and Whitham(1955)]{lighthill1955PRSLA}
Lighthill, M.J.; Whitham, G.B.
\newblock On kinematic waves {II}. {A} theory of traffic flow on long crowded
  roads.
\newblock {\em Proc. R. Soc. Lond. Ser. A} {\bf 1955}, {\em 229},~317--345.
\newblock
  doi:{\changeurlcolor{black}\href{https://doi.org/10.1098/rspa.1955.0089}{\detokenize{10.1098/rspa.1955.0089}}}.

\bibitem[Richards(1956)]{richards1956OR}
Richards, P.I.
\newblock Shock waves on the highway.
\newblock {\em Oper. Res.} {\bf 1956}, {\em 4},~42--51.
\newblock
  doi:{\changeurlcolor{black}\href{https://doi.org/10.1287/opre.4.1.42}{\detokenize{10.1287/opre.4.1.42}}}.

\bibitem[Cristiani and Iacomini(2019)]{cristiani2019DCDS-B}
Cristiani, E.; Iacomini, E.
\newblock An interface-free multi-scale multi-order model for traffic flow.
\newblock {\em Discrete Contin. Dyn. Syst. Ser. B} {\bf 2019}, {\em
  24},~6189--6207.
\newblock
  doi:{\changeurlcolor{black}\href{https://doi.org/10.3934/dcdsb.2019135}{\detokenize{10.3934/dcdsb.2019135}}}.

\bibitem[Garavello and Piccoli(2006)]{piccolibook}
Garavello, M.; Piccoli, B.
\newblock {\em Traffic Flow on Networks}; American Institute of Mathematical
  Sciences,  2006.

\bibitem[Coclite \em{et~al.}(2005)Coclite, Garavello, and
  Piccoli]{coclite2005SIMA}
Coclite, G.M.; Garavello, M.; Piccoli, B.
\newblock Traffic flow on a road network.
\newblock {\em SIAM J. Math. Anal.} {\bf 2005}, {\em 36},~1862--1886.

\bibitem[Holden and Risebro(1995)]{holden1995SIMA}
Holden, H.; Risebro, H.
\newblock A mathematical model of traffic flow on a network of unidirectional
  roads.
\newblock {\em SIAM J. Math. Anal.} {\bf 1995}, {\em 26},~999--1017.
\newblock
  doi:{\changeurlcolor{black}\href{https://doi.org/10.1137/S0036141093243289}{\detokenize{10.1137/S0036141093243289}}}.

\bibitem[Bressan and Nguyen(2015)]{bressan2015NHM}
Bressan, A.; Nguyen, K.T.
\newblock Conservation law models for traffic flow on a network of roads.
\newblock {\em Netw. Heterog. Media} {\bf 2015}, {\em 10},~255--293.

\bibitem[Garavello and Goatin(2012)]{garavello2012DCDS-A}
Garavello, M.; Goatin, P.
\newblock The {C}auchy problem at a node with buffer.
\newblock {\em Discrete Contin. Dyn. Syst. Ser. A} {\bf 2012}, {\em
  32},~1915--1938.

\bibitem[Herty \em{et~al.}(2009)Herty, Lebacque, and Moutari]{herty2009NHM}
Herty, M.; Lebacque, J.P.; Moutari, S.
\newblock A novel model for intersections of vehicular traffic flow.
\newblock {\em Netw. Heterog. Media} {\bf 2009}, {\em 4},~813--826.

\bibitem[Bretti \em{et~al.}(2014)Bretti, Briani, and
  Cristiani]{bretti2014DCDS-S}
Bretti, G.; Briani, M.; Cristiani, E.
\newblock An easy-to-use algorithm for simulating traffic flow on networks:
  Numerical experiments.
\newblock {\em Discrete Contin. Dyn. Syst. Ser. S} {\bf 2014}, {\em
  7},~379--394.
\newblock
  doi:{\changeurlcolor{black}\href{https://doi.org/10.3934/dcdss.2014.7.379}{\detokenize{10.3934/dcdss.2014.7.379}}}.

\bibitem[Briani and Cristiani(2014)]{briani2014NHM}
Briani, M.; Cristiani, E.
\newblock An easy-to-use algorithm for simulating traffic flow on networks:
  Theoretical study.
\newblock {\em Netw. Heterog. Media} {\bf 2014}, {\em 9},~519--552.
\newblock
  doi:{\changeurlcolor{black}\href{https://doi.org/10.3934/nhm.2014.9.519}{\detokenize{10.3934/nhm.2014.9.519}}}.

\bibitem[Hilliges and Weidlich(1995)]{hilliges1995TRB}
Hilliges, M.; Weidlich, W.
\newblock A phenomenological model for dynamic traffic flow in networks.
\newblock {\em Transportation Res. Part B} {\bf 1995}, {\em 29},~407--431.

\bibitem[Briani \em{et~al.}(2018)Briani, Cristiani, and
  Iacomini]{briani2018CMS}
Briani, M.; Cristiani, E.; Iacomini, E.
\newblock Sensitivity analysis of the {LWR} model for traffic forecast on large
  networks using {W}asserstein distance.
\newblock {\em Commun. Math. Sci.} {\bf 2018}, {\em 16},~123--144.
\newblock
  doi:{\changeurlcolor{black}\href{https://doi.org/10.4310/CMS.2018.v16.n1.a6}{\detokenize{10.4310/CMS.2018.v16.n1.a6}}}.

\bibitem[Fan and Work(2015)]{fan2015SIAP}
Fan, S.; Work, D.B.
\newblock A heterogeneous multiclass traffic flow model with creeping.
\newblock {\em SIAM J. Appl. Math.} {\bf 2015}, {\em 75},~813--835.

\bibitem[{van Wageningen-Kessels}(2016)]{kessels2016TRR}
{van Wageningen-Kessels}, F.
\newblock Framework to assess multiclass continuum traffic flow models.
\newblock {\em Transportation Research Record} {\bf 2016}, {\em
  2553},~150--160.
\newblock
  doi:{\changeurlcolor{black}\href{https://doi.org/10.3141/2553-16}{\detokenize{10.3141/2553-16}}}.

\bibitem[{(Sean) Qian} \em{et~al.}(2017){(Sean) Qian}, Li, Li, Zhang, and
  Wang]{qian2017TRB}
{(Sean) Qian}, Z.; Li, J.; Li, X.; Zhang, M.; Wang, H.
\newblock Modeling heterogeneous traffic flow: {A} pragmatic approach.
\newblock {\em Transportation Res. Part B} {\bf 2017}, {\em 99},~183--204.
\newblock
  doi:{\changeurlcolor{black}\href{https://doi.org/10.1016/j.trb.2017.01.011}{\detokenize{10.1016/j.trb.2017.01.011}}}.

\bibitem[Ferrara \em{et~al.}(2018)Ferrara, Sacone, and Siri]{ferrarabook}
Ferrara, A.; Sacone, S.; Siri, S.
\newblock {\em Freeway Traffic Modelling and Control}; Springer,  2018.

\bibitem[Kessels(2019)]{kesselsbook}
Kessels, F.
\newblock {\em Traffic Flow Modelling}; Springer,  2019.

\bibitem[Agarwal and L\"ammel(2016)]{agarwal2016TDE}
Agarwal, A.; L\"ammel, G.
\newblock Modeling \emph{seepage} behavior of smaller vehicles in mixed traffic
  conditions using an agent based simulation.
\newblock {\em Transp. in Dev. Econ.} {\bf 2016}, {\em 2},~8.
\newblock
  doi:{\changeurlcolor{black}\href{https://doi.org/10.1007/s40890-016-0014-9}{\detokenize{10.1007/s40890-016-0014-9}}}.

\bibitem[Benzoni-Gavage and Colombo(2003)]{benzoni2003EJAM}
Benzoni-Gavage, S.; Colombo, R.M.
\newblock An $n$-populations model for traffic flow.
\newblock {\em Euro. Jnl of Applied Mathematics} {\bf 2003}, {\em
  14},~587--612.

\bibitem[Balzotti and G\"ottlich(2021)]{balzotti2021NHM}
Balzotti, C.; G\"ottlich, S.
\newblock A two-dimensional multi-class traffic flow model.
\newblock {\em Netw. Heterog. Media} {\bf 2021}, {\em 16},~69--90.
\newblock
  doi:{\changeurlcolor{black}\href{https://doi.org/10.3934/nhm.2020034}{\detokenize{10.3934/nhm.2020034}}}.

\bibitem[Fan and Seibold(2013)]{fan2013TRR}
Fan, S.; Seibold, B.
\newblock Data-fitted first-order traffic models and their second-order
  generalizations. {C}omparison by trajectory and sensor data.
\newblock {\em Transportation Research Record} {\bf 2013}, {\em 2391},~32--43.
\newblock
  doi:{\changeurlcolor{black}\href{https://doi.org/10.3141/2391-04}{\detokenize{10.3141/2391-04}}}.

\bibitem[Fan \em{et~al.}(2014)Fan, Herty, and Seibold]{fan2014NHM}
Fan, S.; Herty, M.; Seibold, B.
\newblock Comparative model accuracy of a data-fitted generalized
  {Aw-Rascle-Zhang} model.
\newblock {\em Netw. Heterog. Media} {\bf 2014}, {\em 9},~239--268.
\newblock
  doi:{\changeurlcolor{black}\href{https://doi.org/10.3934/nhm.2014.9.239}{\detokenize{10.3934/nhm.2014.9.239}}}.

\bibitem[Klar \em{et~al.}(2004)Klar, G\"unther, Wegener, and
  Materne]{klar2004SIAP}
Klar, A.; G\"unther, M.; Wegener, R.; Materne, T.
\newblock Multivalued fundamental diagrams and stop and go waves for continuum
  traffic flow equations.
\newblock {\em SIAM J. Appl. Math.} {\bf 2004}, {\em 64},~468--483.
\newblock
  doi:{\changeurlcolor{black}\href{https://doi.org/10.1137/S0036139902404700}{\detokenize{10.1137/S0036139902404700}}}.

\bibitem[Herty and Illner(2012)]{herty2012KRM}
Herty, M.; Illner, R.
\newblock Coupling of non-local driving behaviour with fundamental diagrams.
\newblock {\em Kinetic \& Related Models} {\bf 2012}, {\em 5},~843--855.
\newblock
  doi:{\changeurlcolor{black}\href{https://doi.org/10.3934/krm.2012.5.843}{\detokenize{10.3934/krm.2012.5.843}}}.

\bibitem[Ni \em{et~al.}(2018)Ni, Hsieh, and Jiang]{ni2018AMM}
Ni, D.; Hsieh, H.K.; Jiang, T.
\newblock Modeling phase diagrams as stochastic processes with application in
  vehicular traffic flow.
\newblock {\em Appl. Math. Model.} {\bf 2018}, {\em 53},~106--117.

\bibitem[Paipuri and Leclercq(2020)]{paipuri2020TRB}
Paipuri, M.; Leclercq, L.
\newblock Bi-modal macroscopic traffic dynamics in a single region.
\newblock {\em Transportation Res. Part B} {\bf 2020}, {\em 133},~257--290.
\newblock
  doi:{\changeurlcolor{black}\href{https://doi.org/10.1016/j.trb.2020.01.007}{\detokenize{10.1016/j.trb.2020.01.007}}}.

\bibitem[Puppo \em{et~al.}(2016)Puppo, Semplice, Tosin, and
  Visconti]{puppo2016CMS}
Puppo, G.; Semplice, M.; Tosin, A.; Visconti, G.
\newblock Fundamental diagrams in traffic flow: the case of heterogeneous
  kinetic models.
\newblock {\em Commun. Math. Sci.} {\bf 2016}, {\em 14},~643--669.

\bibitem[Visconti \em{et~al.}(2017)Visconti, Herty, Puppo, and
  Tosin]{tosin2017MMS}
Visconti, G.; Herty, M.; Puppo, G.; Tosin, A.
\newblock Multivalued fundamental diagrams of traffic flow in the kinetic
  {Fokker--Planck} limit.
\newblock {\em Multiscale Model. Simul.} {\bf 2017}, {\em 15},~1267--1293.
\newblock
  doi:{\changeurlcolor{black}\href{https://doi.org/10.1137/16M1087035}{\detokenize{10.1137/16M1087035}}}.

\bibitem[Wang \em{et~al.}(2013)Wang, Ni, Chen, and Li]{wang2013JAT}
Wang, H.; Ni, D.; Chen, Q.Y.; Li, J.
\newblock Stochastic modeling of the equilibrium speed-density relationship.
\newblock {\em J. Adv. Transp.} {\bf 2013}, {\em 47},~126--150.
\newblock
  doi:{\changeurlcolor{black}\href{https://doi.org/10.1002/atr.172}{\detokenize{10.1002/atr.172}}}.

\bibitem[Fan \em{et~al.}()Fan, Sun, Piccoli, Seibold, and Work]{CGARZ}
Fan, S.; Sun, Y.; Piccoli, B.; Seibold, B.; Work, D.B.
\newblock A collapsed generalized {Aw-Rascle-Zhang} model and its model
  accuracy.
\newblock ArXiv:1702.03624.

\bibitem[Bretti \em{et~al.}(2018)Bretti, Cristiani, Lattanzio, Maurizi, and
  Piccoli]{bretti2018MiE}
Bretti, G.; Cristiani, E.; Lattanzio, C.; Maurizi, A.; Piccoli, B.
\newblock Two algorithms for a fully coupled and consistently macroscopic
  {PDE-ODE} system modeling a moving bottleneck on a road.
\newblock {\em Mathematics in Engineering} {\bf 2018}, {\em 1},~55--83.
\newblock
  doi:{\changeurlcolor{black}\href{https://doi.org/10.3934/Mine.2018.1.55}{\detokenize{10.3934/Mine.2018.1.55}}}.

\bibitem[Zhao and Zhang(2017)]{zhao2017TRB}
Zhao, Y.; Zhang, H.M.
\newblock A unified follow-the-leader model for vehicle, bicycle and pedestrian
  traffic.
\newblock {\em Transportation Res. Part B} {\bf 2017}, {\em 105},~315--327.

\bibitem[Colombo(2002)]{colombo2002SIAP}
Colombo, R.M.
\newblock Hyperbolic phase transitions in traffic flow.
\newblock {\em SIAM J. Appl. Math.} {\bf 2002}, {\em 63},~708--721.
\newblock
  doi:{\changeurlcolor{black}\href{https://doi.org/10.1137/S0036139901393184}{\detokenize{10.1137/S0036139901393184}}}.

\bibitem[Colombo \em{et~al.}(2010)Colombo, Goatin, and
  Piccoli]{colombo2010JHDE}
Colombo, R.M.; Goatin, P.; Piccoli, B.
\newblock Road networks with phase transitions.
\newblock {\em J. Hyperbolic Differ. Eq.} {\bf 2010}, {\em 7},~85--106.
\newblock
  doi:{\changeurlcolor{black}\href{https://doi.org/10.1142/S0219891610002025}{\detokenize{10.1142/S0219891610002025}}}.

\bibitem[Delle~Monache \em{et~al.}(2021)Delle~Monache, Chi, Chen, Goatin, Han,
  Qiu, and Piccoli]{dellemonache2021AX}
Delle~Monache, M.L.; Chi, K.; Chen, Y.; Goatin, P.; Han, K.; Qiu, J.; Piccoli,
  B.
\newblock Three-phase fundamental diagram from three-dimensional traffic data.
\newblock {\em Axioms} {\bf 2021}, {\em 10},~17.
\newblock
  doi:{\changeurlcolor{black}\href{https://doi.org/10.3390/axioms10010017}{\detokenize{10.3390/axioms10010017}}}.

\bibitem[Cristiani and Sahu(2016)]{cristiani2016NHM}
Cristiani, E.; Sahu, S.
\newblock On the micro-to-macro limit for first-order traffic flow models on
  networks.
\newblock {\em Netw. Heterog. Media} {\bf 2016}, {\em 11},~395--413.
\newblock
  doi:{\changeurlcolor{black}\href{https://doi.org/10.3934/nhm.2016002}{\detokenize{10.3934/nhm.2016002}}}.

\bibitem[Calvert \em{et~al.}(2018)Calvert, {van Wageningen-Kessels}, and
  Hoogendoorn]{calvert2018JTTE}
Calvert, S.C.; {van Wageningen-Kessels}, F.L.M.; Hoogendoorn, S.P.
\newblock Capacity drop through reaction times in heterogeneous traffic.
\newblock {\em Journal of Traffic and Transportation Engineering} {\bf 2018},
  {\em 5},~96--104.

\bibitem[Colombo and Marcellini(2016)]{colombo2016M3AS}
Colombo, R.M.; Marcellini, F.
\newblock A traffic model aware of real time data.
\newblock {\em Math. Models Methods Appl. Sci.} {\bf 2016}, {\em 26},~445--467.
\newblock
  doi:{\changeurlcolor{black}\href{https://doi.org/10.1142/S0218202516500081}{\detokenize{10.1142/S0218202516500081}}}.

\end{thebibliography}
\end{document}